\documentclass[11pt,twoside,a4paper]{article}
\usepackage{mathrsfs}
\usepackage{mathrsfs}
\usepackage{mathrsfs}
\usepackage{mathrsfs}
\usepackage{mathrsfs}
\usepackage{mathrsfs}
\usepackage{mathrsfs}

\usepackage{amsfonts}
\usepackage{amsmath,amssymb}
\usepackage{mathrsfs}
\usepackage{hyperref}

\newtheorem{thm}{Theorem}[section]

\newtheorem{lem}[thm]{Lemma}
\newtheorem{prop}[thm]{Proposition}

\numberwithin{equation}{section}\allowdisplaybreaks

\def\leq{\leqslant}
\def\geq{\geqslant}

\begin{document}

\title{ {\bf \large ABSENCE OF SHOCKS FOR 1D EULER-POISSON SYSTEM}}

{ {\author{\footnotesize YAN GUO,\ \ LIJIA HAN,\ \ JINGJUN ZHANG}}
\date{}
}
\maketitle

\begin{abstract}

 It is shown that smooth solutions with small amplitude to the 1D Euler-Poisson system for electrons persist forever  with no  shock  formation.

\end{abstract}

\section{Introduction}
In this paper, we  consider the 1D Euler-Poisson system in plasma physics:
\begin{align}\label{ep0}
\begin{split}
&n_t  + (nv)_x  = 0, \\
&v_t  +  v v_x + \frac{1}{m_e n}p(n)_x = \frac{e}{m_e} \psi_x
\end{split}
\end{align}
with the electric field  $\psi_{x}$ which satisfies the Poisson equation
\begin{align*}
\psi_{xx}= 4 \pi e (n-n_0), \text{ with }|\psi|\rightarrow 0, \text{ when } x\rightarrow \infty.
\end{align*}
Here, the electrons of charge $e$ and mass $m_e$ are described by a density $n(t, x)$ and an
average velocity $v(t, x)$. The constant equilibrium-charged density of ions and electrons
is $\pm e n_0$. $p$ denotes the  pressure.

Euler-Poisson system \eqref{ep0}  describes the simplest  two-fluid model in plasma physics.
In this model, the ions are treated as immobile and only form a constant charged
background $n_0$. The ¡°two-fluid¡± models
describe dynamics of two separate compressible fluids of ions and electrons interacting
with their self-consistent electromagnetic field. As pointed in the classical book of Jackson \cite[P. 337]{Jackson}, ``The
adiabatic law $p=p_{0}(n/n_{0})$ can be assumed, but the customary accoustic
value $\gamma =\frac{5}{3}$ for a gas of particles with 3 external, but no
internal, degrees of freedom is not valid. The reason is that the frequency
of the present density oscillations is much higher than the collision
frequency, contrary to the acoustical limit. Consequently the
one-dimensional nature of the density oscillations is maintained. A value of
$\gamma $ appropriate to $1$ translational degree of freedom must be used.
Since $\gamma =(m+2)/m,$ where $m$ is the number of degrees of freedom, we
have in this case $\gamma =3.$''  We therefore concentrate in this paper on this
most significant physical case,  and assume the pressure is given by
\begin{align}\label{pressure}
p(n)=\frac{1}{3}n^3.
\end{align}

 In 1998, Guo in \cite{Guo}  first  studied  Euler-Poisson system in three dimensional case. He observed that the linearized Euler-Poisson system for the electron fluid is the Klein-Gordon equation, due to plasma oscillations
created by the electric field, and constructed the smooth irrotational solutions with small amplitude for all time
 (never develop shocks).  This is a very surprising result compared to the work of Sideris \cite{ST} for
pure Euler equations, where the solutions will blow up even under small perturbations. It is the dispersive effect of
the electric field that enhances the linear decay rate and prevents shock formation. Note that the decay rate in the
$L^\infty - L^1$ decay estimate for the linear Klein-Gordon equation is $t^{-\frac{d}{2}}$, which is integrable when
$d=3$.

In lower dimension case (1D and 2D case), as the decay rate for the linearized Euler-Poisson equations is worse than 3D
case, so that the construction of global smooth solution is much more challenging. In 2D case, the decay rate in the $L^\infty - L^1$ decay estimate is the borderline $t^{-1}$, so the main obstructions in the 2D Euler-Poisson system are slow (non-integrable) dispersion. Recently,
smooth irrotational solutions for the 2D Euler-Poisson system \eqref{ep0} are constructed independently in \cite{IP,LW}. See also \cite{Jang,JLZ} for related results on two dimensional case.

Such an unexpected and subtle dispersive effect has been discovered
and exploited in other two-fluid models, which leads to persistence of
global smooth solutions and absence of shock formations. Among the results, we refer to \cite{Deng,GP,GMP,GIP,GIP2,GuoPausader,IP2,HZG}.

It has remained as an outstanding question about whether or not shock
formations can be suppressed in 1D for any two-fluid model. For the
Euler-Poisson system \eqref{ep0}, the linear time decay rate is merely of $%
t^{-1/2},$ and even for general 1D scalar nonlinear Klein-Gordon equation,
singularity (shock waves) might develop for small initial data \cite{Hormander}.
Nevertheless, we settle this question in affirmative for the Euler-Poisson
system \eqref{ep0} with $\gamma =3$ by constructing global smooth solutions with
small amplitude$.$  To state precisely our result, we set all the physical constants $m_e$, $e$, $n_0$ and $4\pi$ to be one. From \eqref{pressure},  system \eqref{ep0} reduces to
\begin{align}\label{ep001}
\begin{split}
&n_t  + (nv)_x  = 0, \\
&v_t  +  v v_x + nn_x =  \psi_x,\\
& \psi_{xx}=n-1.
\end{split}
\end{align}
Moreover, if $E:=\psi_x$, system \eqref{ep001} can be further rewritten as
\begin{align}\label{ep2}
\begin{split}
&E_t  + v + v E_x  = 0, \\
&v_t  + E_{xx} -  E + v v_x + E_x E_{xx} = 0.
\end{split}
\end{align}
From now on, we mainly focus on the above system. Let
\begin{align}\label{trans1}
 r:= \frac{E}{2}, \quad u:= -\frac{v}{2\langle\partial_x\rangle},
 \end{align}
where $\langle\partial_x\rangle:=\sqrt{1-\partial_{x}^2}$, then \eqref{ep2} can be written in an equivalent form
\begin{align}\label{ep3}
\begin{pmatrix}
           r \\
           u \\
         \end{pmatrix}_t + \begin{pmatrix}
                0  & -\langle\partial_x\rangle \\
               \langle\partial_x\rangle  & 0 \\
               \end{pmatrix}
                \begin{pmatrix}
           r \\
           u \\
         \end{pmatrix} = \begin{pmatrix}
           2 \langle\partial_x\rangle u  \ r_x \\
           \frac{\partial_x}{\langle\partial_x\rangle}[ (\langle\partial_x\rangle u)^2 + (r_x)^2 ] \\
         \end{pmatrix}.
\end{align}

Once we obtain global smooth solutions $(r,u)$  for system \eqref{ep3}, then we also obtain smooth solutions $(E, v)$ for system \eqref{ep2} by the relation  \eqref{trans1}, and thus the density $n$ in \eqref{ep001} is $n=1+\psi_{xx}= 1 + E_{x}$.

The main result of the paper is stated in the following theorem.

\begin{thm}\label{mainthm}
Let $N=300,\ N_1=15,\ 0<p_0<10^{-3}$,  $U:=(r,u)^T$ and $\Gamma:=t\partial_x+x\partial_t$. Then there exists $\epsilon_0=\epsilon_0(p_0)>0$ sufficiently small such that if
\begin{align}\label{state1}
\|U(0)\|_{H^N}+\|xU(0)\|_{H^{N_1+1}}+\|\langle\xi\rangle^{N_1+10}\widehat{U(0)}\|_{L^\infty}\leq  \epsilon_0,
\end{align}
the system \eqref{ep3} admits a global solution $U\in C(\mathbb{R}^+;H^N)$ satisfying
\begin{align*}
\sup_{t>0}\big[(1+t)^{-p_0}\|U(t)\|_{H^N}&+(1+t)^{-p_0}\|\Gamma U(t)\|_{H^{N_1}}
+(1+t)^{1/2}\|U(t)\|_{W^{N_1+10,\infty}}\big]\lesssim  \epsilon_0.
\end{align*}
\end{thm}

We remark that \eqref{state1} implies the  neutrality condition
$$
\int_{\mathbb{R}}(n(0,x)-1)dx=0,
$$
which is conserved for all time. The above theorem shows that under small perturbations around the equilibrium, system \eqref{ep001} still has a global
smooth solution. However, unlike the 2D or 3D case,  we can not obtain the usual scattering result for 1D Euler-Poisson
system. Instead,  we will see that solutions approach to a nonlinear asymptotic state. To show this phenomenon,  we
set
 \begin{align}\label{defofh}
 h=\frac{1}{2}E-\frac{i}{2\langle\partial_x\rangle}v=r+iu,
 \end{align}
then system \eqref{ep3} is equivalent to the following complex-valued Klein-Gordon  equation
\begin{align}\label{ep002}
h_t+i\langle\partial_x\rangle h=&\frac{1}{2i}(h+\overline{h})_x\langle\partial_x\rangle(h-\overline{h})
+\frac{\partial_x}{4i\langle\partial_x\rangle}[\langle\partial_x\rangle(h-\overline{h})]^2
-\frac{\partial_x}{4i\langle\partial_x\rangle}[(h+\overline{h})_x]^2.
\end{align}
By Shatah's normal form transformation \cite{Shatah}, we may make a change
of new unknown $g$ (see \eqref{ps-a1}) such that
\begin{align}\label{introps3}
g_t+i\langle\partial_x\rangle g=\mathcal{N}(h),
\end{align}
where the cubic term $\mathcal{N}(h)$ is given in  \eqref{ps4}.  In this work, we show that there exists a unique $w_\infty(\xi)\in L^\infty$ such that
\begin{align*}
\sup_{t\geq 0}[(1+t)^\delta\|\langle\xi\rangle^{N_1+10} e^{i\vartheta(t,\xi)}\widehat{w}(t,\xi)-w_\infty(\xi)\|_{L^\infty}]\lesssim \epsilon_0
\end{align*}
for some $\delta>0$, where $w:=e^{i\langle \partial _{x}\rangle }g$ is the linear profile of $g$, and $\vartheta$ is a real-valued function defined by  \eqref{definition of theta}. This result says the solution of the equation \eqref{ep002} tends to a nonlinear asymptotic state as time goes to infinity, thus  such equation possesses a modified scattering behavior. Therefore, we extend the previous work on electron type Euler-Poisson system to one dimensional case. Together with the work \cite{Guo,IP,LW}, our result provides a complete picture of Klein-Gordon effect which prevents formation of small shocks
in all physical dimensions for the Euler-Poisson system \eqref{ep0}. Moreover, even though \eqref{ep0} is a hyperbolic
system of conservation (balance) law \cite{Dafermous}, the construction of its global BV solutions with
small amplitude (hence uniqueness) has remained outstanding. Our result also
demonstrates that the standard 1D BV theory is not needed for small smooth initial
data for \eqref{ep0}, and it is ill-suited to capture the delicate dispersive
Klein-Gordon effect which prevents the shock formation.

Our work is inspired by recent work of \cite{AD,IPu2,IPu3,IPu4} on water waves system,
which depends on a delicate interplay between higher energy estimates and a low
order $L^{\infty }$ estimate.  It is well-known that due to poor decay rate
of $t^{-1/2}$, the classical energy estimate with quadratic nonlinearity is impossible to close, and it is necessary to perform the energy estimate in a new
system with a \textit{cubic} nonlinearity. In other words, one would wish to
make an ``energy normal form'' transformation in the energy estimate.
Unfortunately, Shatah's normal form transformation introduces ``loss'' of
derivatives. Even though it is sufficient for lower order $L^{\infty }$
decay estimates, it is in general not compatible for high order energy
estimate. As a matter of fact, such an ``energy normal form'' may not exist
for general 1D quasi-linear Klein-Gordon equations.

Our first important step is
the construction of an ``energy normal form'' transformation in Section 2. We
follow the procedure in \cite{AD}, and  the special structure with $\gamma =3$
enables us to discover subtle cancelations for the  part of the quadratic terms $Q-B$ (see Proposition \ref{p3.1}), during the Sobolev
energy estimates. Meanwhile, we construct normal form transformations without ``loss'' of
derivatives, which
eliminates the other part of quadratic terms $B$ (see Proposition \ref{p4}). In Section 2.4, we complete the whole process of higher order energy estimates (Proposition \ref{eprop}).

For the $L^{\infty} $ decay estimate, we employ the following refined
linear decay estimate for the solution $g$ (see Lemma A.1),
\begin{align}\label{introdecay}
\|g\|_{L^\infty}\lesssim (1+t)^{-\frac{1}{2}}\|\widehat{w} \|_{L^\infty}
+(1+t)^{-\frac{5}{8}}(\|w\|_{H^2}+\|xw\|_{H^1}),\ \forall\ t\geq 0,
\end{align}
where $w=e^{it\langle\partial_x\rangle}g$. It is important to note that $-\frac{5}{8}<-\frac{1}{2},$ so there is room
for mild growth for $||w||_{H^{2}}$ and $||xw||_{H^{1}}$. Then it reduces to low order estimates for $\|xw\|_{H^{N_1}}$ and $\|\langle \xi \rangle^{N_1+10}\widehat{w} \|_{L^\infty}$, respectively.

The second important step is to estimate $||xw||_{H^{N_1}}$. In the work \cite{AD}
and \cite{IPu2}, the crucial homogeneous scaling operator $S = \frac{1}{2} t \partial _{t}+ x\partial _{x}$ for the gravity water waves system is employed.
Unfortunately, in our problem, the natural  operator\footnote{The operator $\widetilde{\Gamma}$ was used in \cite{HN1,HN3} to study the scattering behavior for cubic and quadratic nonlinear Klein-Gordon equation without derivatives.} $
\widetilde{\Gamma}= t\partial_x-i\langle\partial_x\rangle x$ for the Klein-Gordon case, is not homogeneous. So the
energy estimate fails for $\widetilde{\Gamma} U$, as $\widetilde{\Gamma}$ could not commute with the nonlinear terms. Instead, we use the homogeneous vector field operator $\Gamma = t \partial_x + x \partial_t$  to perform energy estimate for $\Gamma U$. The key observation is the  following relation between $\Gamma$ and $\widetilde{\Gamma}$,
 \begin{align*}
\widetilde{\Gamma} g = \Gamma g - x(\partial_t + i \langle \partial_x \rangle )g+ \frac{i\partial_x}{\langle\partial_x\rangle}g
=\Gamma g - x \mathcal{N}(h)+ \frac{i\partial_x}{\langle\partial_x\rangle}g.
\end{align*}
 In Sections 3.1-3.3, we obtain energy estimate for $\Gamma U$ by applying similar strategy as used in Section 2, see Proposition \ref{weprop}. In addition, using this modified normal form process,  we also control the low  energy estimate for $x U$ in Section 3.4, which is necessary when estimating the difference between $\Gamma g$ and $\widetilde{\Gamma} g$. We establish that  $\|xU\|_{H^{N_1}}$ grows almost linearly
$$
\|xU(t)\|_{H^{N_1}}\lesssim (1+t)^{1+p_0}
$$
for $p_0\ll1$, see Proposition \ref{energyxuprop}. Thanks again to the \emph{cubic} structure of $\mathcal{N}(h)$, it yields that $\|x\mathcal{N}(h)\|_{H^{N_1}}$ can be bounded by $(1+t)^{p_0}$, which is sufficient  for our argument. In virtue of the identity
$$
\langle\partial_x\rangle (xw)= ie^{it\langle\partial_x\rangle}\widetilde{\Gamma} g,
$$
 we finally can able to control  $||xw||_{H^{N_1}}$ via the estimates of $g$, $\Gamma g$ and $x\mathcal{N}(h)$. The details are presented in  Section 4.1.

The estimate for $||\langle\xi\rangle^{N_1+10}\widehat{w}||_{L^\infty }$ is carried out in Section 4.2 as an
adaptation of the proof in \cite{IPu1,IPu2,IPu3,IPu4}. Through precise frequency decompositions and stationary phase analysis, we notice that a phase correction is needed to the leading order term and thus leads to the modified scattering behavior (Proposition \ref{scproposition}). Using the above norms and \eqref{introdecay}, we close our decay argument in Section 4.3.

Finally,  the  global existence result follows from \eqref{introdecay}, Proposition \ref{eprop}, Propositions \ref{weprop}--\ref{energyxuprop} and Proposition \ref{decaye}.

\textbf{Notations:}

$\bullet$ The Fourier transform and Fourier inverse transform are defined by
\begin{align*}
(\mathscr{F}f)(\xi)&=\int_{\mathbb{R}^d}e^{-ix\cdot \xi}f(x)dx=\widehat{f}(\xi),\\
(\mathscr{F}^{-1}g)(x)&=(2\pi)^{-d}\int_{\mathbb{R}^d}e^{ix\cdot \xi}g(\xi)d\xi.
\end{align*}

$\bullet$ Assume $f$ is a scalar  function, $V$  is a vector-valued function (or scalar function) and $ M(\xi_1, \xi_2)$ is a matrix symbol (or scalar symbol). Define the bilinear operator
\begin{align}\label{defofO}
\mathcal{O}[f,M]V:=\frac{1}{(2 \pi)^2} \int_{\mathbb{R}^2}  e^{i x (\xi_1 +\xi_2)}  \widehat{f}(\xi_1) M(\xi_1, \xi_2)\widehat{V}(\xi_2)d\xi_1 d\xi_2.
\end{align}

$\bullet$  Let
$\varphi\in C_c^\infty(\mathbb{R})$ be a radial function with the properties such that $0\leq \varphi\leq 1$, $\varphi(\xi)=1$ for $|\xi|\leq 5/4$ and $\rm {supp \varphi \subset [-8/5, 8/5]}$. Then for $k\in \mathbb{\mathbb{Z}}$, we write $\varphi_k(\xi):=\varphi(\xi/2^k)-\varphi(\xi/2^{k-1})$. The
dyadic frequency localization operator $P_k$ is defined by
$$\widehat{P_kf}(\xi):= \varphi_k(|\xi|)\widehat{f}(\xi).$$
Moreover, for $a>0$, we denote by $P_\leq a$ the projector with symbol $\varphi(|\xi|/a)$.

$\bullet$  For any $\rho \in \{0\}\cup \mathbb{ N}$, we denote by $C^0(\mathbb{ R})$ the space of bounded continuous functions,
by $C^\rho(\mathbb{ R})$ the space of $C^0(\mathbb{ R})$ functions whose derivatives of order less or equal to $\rho$ are in $C^0(\mathbb{ R})$.

$\bullet$ $\langle\partial_x\rangle:=\sqrt{1-\partial_{x}^2}$, $\langle\xi\rangle:=\sqrt{1+\xi^2}$.

\section{Energy  estimate}

Our aim in  this section is to prove the following energy estimate.
\begin{prop}\label{eprop} Let $U(t)\in C([0,T];H^N)$ be the solution of system \eqref{ep3}. Assume \eqref{state1} holds and \begin{align}\label{main a prior bound1}
\sup_{t\in [0,T]}\big[(1+t)^{-p_0}\|U(t)\|_{H^N}+(1+t)^{1/2}\|U(t)\|_{W^{N_1+10,\infty}}\big]\lesssim  \epsilon_1,
\end{align}
where $0<\epsilon_0\ll\epsilon_1\ll 1$, $N=300$, $N_1=15$ and $0<p_0<10^{-3}$. Then we have
\begin{align}\label{main energy bound1}
\sup_{t\in [0,T]}\big[(1+t)^{-p_0}\|U(t)\|_{H^N}\big]\lesssim \epsilon_0+ \epsilon_1^2,
\end{align}
where the implicit constant depends only on $p_0$.
\end{prop}

\subsection{Decomposition of the nonlinear terms}

Fix a cut off function $\theta\in C^\infty(\mathbb{R} \times \mathbb{R})$ satisfying\\
(1)  There exist $\tilde{\epsilon}_1 $, $\tilde{\epsilon}_2$ such that  $0< 2 \tilde{\epsilon}_1 < \tilde{\epsilon}_2 <1/2$ and \begin{align}\label{property1}
\begin{split}
&\theta(\xi_1, \xi_2)=1, \quad |\xi_1|\leq \tilde{\epsilon}_1 |\xi_2| ,\\
&\theta(\xi_1, \xi_2)=0, \quad |\xi_1| \geq \tilde{\epsilon}_2 |\xi_2|.
 \end{split}
 \end{align}
(2) For any $\alpha$, $\beta\in \mathbb{N}\cup\{0\}$, there holds
\begin{align}\label{property2}
|\partial_{\xi_1}^{\alpha}\partial_{\xi_2}^\beta\theta(\xi_1, \xi_2)|\leq C_{\alpha,\beta}\langle\xi_2\rangle^{-\alpha-\beta},\quad \forall\ \xi_1,\xi_2\in \mathbb{R}.
\end{align}
(3) $\theta$ satisfies the symmetry condition
\begin{align}\label{property3}
\theta(\xi_1, \xi_2) = \theta(-\xi_1, -\xi_2) = \theta(-\xi_1, \xi_2).
\end{align}
Then we define the paraproduct $T_fg$ and the remainder $R_B(f, g)$ as
\begin{align}
T_f g&: = \frac{1}{(2 \pi)^2} \int_{\mathbb{R}^2}  e^{i x (\xi_1 +\xi_2)} \theta(\xi_1, \xi_2) \widehat{f}(\xi_1)\widehat{g}(\xi_2)d\xi_1 d\xi_2, \nonumber\\
R_B(f, g)& : = \frac{1}{(2 \pi)^2} \int_{\mathbb{R}^2} e^{i x (\xi_1 +\xi_2)} (1- \theta(\xi_1, \xi_2)-\theta(\xi_2, \xi_1))\widehat{f}(\xi_1)\widehat{g}(\xi_2)d\xi_1 d\xi_2. \nonumber
\end{align}
With this definition, for any $f$ and $g$, we have the following Bony decomposition
\begin{align}
&f g = T_f g + T_g f + R_B(f, g). \label{bony}
\end{align}

Let
\begin{align}\label{DU}
&U = \begin{pmatrix}
           r \\
           u \\
         \end{pmatrix}, \quad\quad
D= \begin{pmatrix}
                0  & -\langle\partial_x\rangle \\
               \langle\partial_x\rangle  & 0 \\
               \end{pmatrix},
\end{align}
and
\begin{align}
& Q_1(\xi_1, \xi_2):=\begin{pmatrix}
   0 & q_1(\xi_1,\xi_2) \\
     q_4(\xi_1,\xi_2) & 0 \\
     \end{pmatrix}=
\begin{pmatrix}
   0 & 2i\xi_1\langle\xi_2\rangle \\
     \frac{-2i (\xi_1 +\xi_2)\xi_1\xi_2}{\langle\xi_1+\xi_2\rangle } & 0 \\
     \end{pmatrix}
\theta(\xi_1, \xi_2),\label{Q1}\\
& Q_2(\xi_1, \xi_2):=\begin{pmatrix}
   q_2(\xi_1,\xi_2) &0 \\
     0 & q_3(\xi_1,\xi_2)  \\
     \end{pmatrix}=
\begin{pmatrix}
 2i \xi_2 \langle\xi_1\rangle  & 0\\
   0 &\frac{2i( \xi_1 +  \xi_2) \langle\xi_1\rangle \langle\xi_2\rangle }{\langle\xi_1+\xi_2\rangle } \\
   \end{pmatrix}
\theta(\xi_1, \xi_2),\label{Q2}\\
& S_1(\xi_1, \xi_2):=\begin{pmatrix}
   0 & s_1(\xi_1,\xi_2) \\
     s_4(\xi_1,\xi_2) & 0 \\
     \end{pmatrix}=
\begin{pmatrix}
  0 & i\xi_1\langle\xi_2\rangle  \\
   \frac{-i (\xi_1 +\xi_2)\xi_1\xi_2}{\langle\xi_1+\xi_2\rangle } & 0 \\
    \end{pmatrix}
(1- \theta(\xi_1, \xi_2)-\theta(\xi_2, \xi_1)), \label{S1}\\
& S_2(\xi_1, \xi_2):=\begin{pmatrix}
   s_2(\xi_1,\xi_2) &0 \\
     0 & s_3(\xi_1,\xi_2)  \\
     \end{pmatrix}=
\begin{pmatrix}
 i \xi_2 \langle\xi_1\rangle  & 0\\
    0 &\frac{i( \xi_1 +  \xi_2) \langle\xi_1\rangle \langle\xi_2\rangle }{\langle\xi_1+\xi_2\rangle } \\
     \end{pmatrix}
(1- \theta(\xi_1, \xi_2)-\theta(\xi_2, \xi_1)). \label{S2}
\end{align}
By \eqref{defofO} and \eqref{bony}, system \eqref{ep3} is then transformed into
\begin{align}\label{main}
 U_t + D U =  \mathcal{O}[r, Q_1]U + \mathcal{O}[u, Q_2]U +  \mathcal{O}[r, S_1]U + \mathcal{O}[u, S_2]U .
\end{align}
Here, $Q_1$, $Q_2$ are the symbols of low-high interaction terms, with one  local/global derivative on the function of  high frequency, and $S_1$, $S_2$ are the symbols of nonlinear terms with high-high interactions. The low-high terms will cause loss of derivatives when performing energy estimate, so we shall do some modifications with these terms, see Proposition \ref{p3.1} in the next subsection.

\subsection{Modifying low-high interaction terms}
In this subsection, we prove
\begin{prop}\label{p3.1} Let $Q_1$, $Q_2$ be given by \eqref{Q1}--\eqref{Q2}. Then there exist two matrices $B_1$ and $B_2$  with the form
\begin{align}\label{B12}
B_1=\begin{pmatrix}
0 & b_1(\xi_1, \xi_2)  \\
b_4(\xi_1, \xi_2) &  0 \\
\end{pmatrix}, \quad\quad
B_2=\begin{pmatrix}
 b_2(\xi_1, \xi_2) &  0 \\
 0 &  b_3(\xi_1, \xi_2) \\
 \end{pmatrix}
\end{align}
such that
\begin{align}
&\mathrm{Re} \langle\langle\partial_x\rangle^N \mathcal{O}[r, Q_1-B_1]U, \langle\partial_x\rangle^N U\rangle = 0, \label{2.1}\\
&\mathrm{Re} \langle\langle\partial_x\rangle^N \mathcal{O}[u, Q_2-B_2]U, \langle\partial_x\rangle^N U\rangle = 0, \label{2.2}
\end{align}
where $\langle\cdot, \cdot\rangle$ denotes the inner product of  $L^2$ space. Moreover, for any $ \alpha,\beta=0,1$,
\begin{align}\label{bjestimate}
|\partial_{\xi_1}^\alpha\partial_{\xi_2}^\beta b_j(\xi_1,\xi_2) | \lesssim \langle\xi_1\rangle^2, \quad j=1,2,3,4,
\end{align}
and for any $\rho\geq 3$,
\begin{align}
&\|\langle\partial_x\rangle^{N}\mathcal{O}[f,B_1]U \|_{L^2} + \|\langle\partial_x\rangle^{N}\mathcal{O}[f,B_2] U \|_{L^2} \lesssim \|f\|_{C^\rho}  \| U\|_{H^N}.\label{4.5}
\end{align}
\end{prop}

To prove this proposition, one should use the following lemma.

\begin{lem}\label{ld} \cite{AD} Assume $f$ is a real-valued function, and $M$ is a matrix. Then we have
$$ (\mathcal{O}[f, M])^* = \mathcal{O}[f, \widetilde{M}],\ \ \ \   \widetilde{M}(\xi_1, \xi_2) := \overline{M^T(-\xi_1, \xi_1 + \xi_2)},$$
where $M^T$ is the transpose of $M$.
\end{lem}

\noindent{ \emph{Proof.}} By the definition of $\mathcal{O}[f, M]W$ (see \eqref{defofO}),
\begin{align}
\mathscr{F}(\mathcal{O}[f, M]W)(\eta) = \frac{1}{2 \pi} \int_{\mathbb{R}}  \widehat{f}(\xi)M(\xi, \eta -\xi)\widehat{W}(\eta-\xi)d\xi.\nonumber
\end{align}
Using  the fact $f$ is real-valued, we have
\begin{align*}
\langle \mathcal{O}[f, M]^* V,W \rangle
&= \langle V, \mathcal{O}[f, M]W \rangle \\
& = \frac{1}{(2 \pi)^2} \int_{\mathbb{R}}  \widehat{V}(\eta) \Big( \int_{\mathbb{R}} \overline{\widehat{f}(\xi)M(\xi, \eta-\xi)\widehat{W}(\eta-\xi)} d\xi  \Big) d\eta\\
& = \frac{1}{(2 \pi)^2} \int_{\mathbb{R}^2} \widehat{V}(\xi + \eta)\widehat{f}(-\xi)\overline{M(\xi, \eta)}\ \overline{\widehat{W}(\eta)} d\xi d\eta \\
& =  \frac{1}{ 2 \pi }  \int_{\mathbb{R}}  \frac{1}{ 2 \pi } \Big( \int_{\mathbb{R}}  \widehat{f}(\xi)\overline{M^T(-\xi, \eta)}\widehat{V}(\eta - \xi) d\xi\Big) \overline{\widehat{W}(\eta)} d\eta \nonumber\\
&=  \langle \mathcal{O}[f, \widetilde{M}] V,W \rangle,
\end{align*}
from which we can obtain the desired result  $\widetilde{M}(\xi,\eta)=\overline{M^T(-\xi,\xi+\eta)}$.
$\hfill \square$

We will also need the following anisotropic multiplier estimate.
\begin{lem}\label{l5.1}
There holds
\begin{align}
\|\mathcal{O}[f, M]V\|_{L^2(\mathbb{R})} \lesssim \| M(\xi, \eta- \xi)\|_{L_\eta^\infty H_\xi^1}\|f\|_{L^\infty (\mathbb{R})}\|V\|_{L^2(\mathbb{R})}.\label{5.5}
\end{align}
\end{lem}

Similar estimates are also used in \cite{GuoPausader, GT}.
The proof of this lemma is given in Lemma B.1 of the appendix.
\vspace{5mm}

\noindent{\emph{Proof of Proposition \ref{p3.1}.}}
We rewrite equation \eqref{2.1} as
\begin{align}
0=&2 \mathrm{Re} \langle\langle\partial_x\rangle^N \mathcal{O}[r,Q_1-B_1] U, \langle\partial_x\rangle^N U\rangle  \nonumber\\
=&  \langle\langle\partial_x\rangle^N \mathcal{O}[r,Q_1-B_1] U, \langle\partial_x\rangle^N U\rangle + \overline{ \langle\langle\partial_x\rangle^N \mathcal{O}[r,Q_1-B_1] U, \langle\partial_x\rangle^N U\rangle}\nonumber\\
=&  \langle\langle\partial_x\rangle^N \mathcal{O}[r,Q_1-B_1] U,\langle\partial_x\rangle^N U\rangle + \langle\langle\partial_x\rangle^N U,  \langle\partial_x\rangle^N \mathcal{O}[r,Q_1-B_1] U\rangle \nonumber\\
=&  \langle\langle\partial_x\rangle^{2N}\mathcal{O}[r,Q_1-B_1] U,  U\rangle + \langle\langle\partial_x\rangle^{2N} U, \mathcal{O}[r,Q_1-B_1] U\rangle \nonumber\\
=&\langle\langle\partial_x\rangle^{2N}\mathcal{O}[r,Q_1-B_1] U,  U\rangle + \langle\mathcal{O}[r,Q_1-B_1]^*(\langle\partial_x\rangle^{2N} U), U\rangle. \nonumber
\end{align}
In order to prove \eqref{2.1}, we only need to verify
\begin{align}
\langle\partial_x\rangle^{2N}\mathcal{O}[r,Q_1] U +  \mathcal{O}[r,Q_1]^* (\langle\partial_x\rangle^{2N}U) =  \langle\partial_x\rangle^{2N}\mathcal{O}[r,B_1] U + \mathcal{O}[r,B_1]^* (\langle\partial_x\rangle^{2N}U). \label{2.18.1}
\end{align}
From Lemma \ref{ld}, we see
\begin{align}
 \langle\partial_x\rangle^{2N} \mathcal{O}[r,Q_1] U &= \mathcal{O}[r, \langle\xi_1+\xi_2\rangle^{2N} Q_1(\xi_1 , \xi_2)] U, \nonumber\\
  \mathcal{O}[r,Q_1]^* (\langle\partial_x\rangle^{2N}U)& = \mathcal{O}[r,  \langle\xi_2\rangle^{2N} \overline{Q_1^T(-\xi_1 , \xi_1+ \xi_2)}] U. \nonumber
\end{align}
Define
\begin{align}\label{B1}
B_1 (\xi_1, \xi_2): = \frac{\langle\xi_1 + \xi_2\rangle^{2N} Q_1(\xi_1 , \xi_2) + \langle\xi_2\rangle^{2N} \overline{ Q_1^T(-\xi_1 , \xi_1+ \xi_2)}}{ \langle\xi_1 + \xi_2\rangle^{2N} + \langle\xi_2\rangle^{2N}},
\end{align}
then $B_1(\xi_1,\xi_2)=\overline{B_1^T(-\xi_1,\xi_1+\xi_2)}$, and by Lemma \ref{ld}, we have $ \mathcal{O}[r,B_1]= \mathcal{O}[r,B_1]^* $. With such choice of $B_1(\xi_1,\xi_2)$,  the identity \eqref{2.18.1} holds, and \eqref{2.1} thus follows.

Similarly, in order to prove \eqref{2.2}, we only need to show
\begin{align}
\langle\partial_x\rangle^{2N}\mathcal{O}[u,Q_2] U +  \mathcal{O}[u,Q_2]^* (\langle\partial_x\rangle^{2N}U) =  \langle\partial_x\rangle^{2N}\mathcal{O}[u,B_2] U + \mathcal{O}[u,B_2]^* (\langle\partial_x\rangle^{2N}U). \nonumber
\end{align}
Define
\begin{align}\label{B2}
B_2 (\xi_1, \xi_2): = \frac{\langle\xi_1 + \xi_2\rangle^{2N} Q_2(\xi_1 , \xi_2) + \langle\xi_2\rangle^{2N} \overline{ Q_2^T(-\xi_1 , \xi_1+ \xi_2)}}{ \langle\xi_1 + \xi_2\rangle^{2N} + \langle\xi_2\rangle^{2N}},
\end{align}
then we can check $ \mathcal{O}[u,B_2]= \mathcal{O}[u,B_2]^* $ and  \eqref{2.2} holds.

In order to prove \eqref{bjestimate}, we should calculate $b_j (\xi_1, \xi_2)$ ($j=1,2,3,4$) carefully. From \eqref{B12}, \eqref{B1}, \eqref{Q1} and \eqref{property3}, the expression for  $b_1(\xi_1 , \xi_2)$   can be written as exactly as
\begin{align}
&b_1(\xi_1 , \xi_2) \nonumber\\
& = \frac{2i \xi_1}{\langle\xi_1 + \xi_2\rangle^{2N} + \langle\xi_2\rangle^{2N}}\Big[\langle\xi_1 + \xi_2\rangle^{2N} \langle \xi_2\rangle\theta(\xi_1,  \xi_2) -\langle\xi_2\rangle^{2N}\frac{\xi_2 (\xi_1 + \xi_2)}{\langle \xi_2\rangle} \theta(-\xi_1, \xi_1+ \xi_2)\Big]\nonumber\\
&= \frac{2i \xi_1}{(\langle\xi_1 + \xi_2\rangle^{2N} + \langle\xi_2\rangle^{2N})\langle  \xi_2\rangle}\left[\langle\xi_1 + \xi_2\rangle^{2N} \langle  \xi_2\rangle^2\theta(\xi_1,  \xi_2) -\langle\xi_2\rangle^{2N}\xi_2 (\xi_1 + \xi_2)\theta(\xi_1, \xi_1+ \xi_2) \right]\nonumber\\
&= \frac{2i \xi_1 \xi_2}{(\langle\xi_1 + \xi_2\rangle^{2N} + \langle\xi_2\rangle^{2N})\langle  \xi_2\rangle }\left[(\xi_1 + \xi_2)^{2N}   \xi_2 \theta(\xi_1,  \xi_2) -  (\xi_1 + \xi_2) \xi_2^{2N}\theta(\xi_1, \xi_1+ \xi_2)\right] +  r_1(\xi_1, \xi_2)\nonumber\\
&= \frac{2i \xi_1 \xi_2}{(\langle\xi_1 + \xi_2\rangle^{2N} + \langle\xi_2\rangle^{2N})\langle  \xi_2\rangle }\chi(\xi_1,  \xi_2) +  r_1(\xi_1, \xi_2),\label{2.5}
\end{align}
where
\begin{align}
\chi(\xi_1,\xi_2)&:=(\xi_1+\xi_2)^{2N}\xi_2\theta(\xi_1,\xi_2)-(\xi_1+\xi_2)\xi_2^{2N}\theta(\xi_1,\xi_1+\xi_2),
\label{definitionchi}\\
r_1(\xi_1,\xi_2) &:=\frac{ 2i \xi_1  \theta(\xi_1,  \xi_2) }{\langle  \xi_2\rangle (\langle\xi_1 + \xi_2\rangle^{2N} + \langle\xi_2\rangle^{2N}) }[\langle\xi_1 + \xi_2\rangle^{2N} \langle  \xi_2\rangle^2
- (\xi_1 + \xi_2)^{2N} \xi_2^2 ]\nonumber\\
&\quad-\frac{ 2i \xi_1\xi_2 (\xi_1 + \xi_2)\theta(\xi_1, \xi_1+ \xi_2)}{\langle  \xi_2\rangle (\langle\xi_1 + \xi_2\rangle^{2N} + \langle\xi_2\rangle^{2N}) }[\langle\xi_2\rangle^{2N}- \xi_2^{2N}].\nonumber
\end{align}
We decompose $\chi(\xi_1,\xi_2)$ into $I_1+I_2$ with
\begin{align*}
 I_1&:=((\xi_1+\xi_2)^{2N}\xi_2-(\xi_1+\xi_2)\xi_2^{2N})\theta(\xi_1,\xi_2),\\
 I_2&:= (\xi_1+\xi_2)\xi_2^{2N}( \theta(\xi_1,\xi_2) - \theta(\xi_1,\xi_1+\xi_2)).
\end{align*}
Recall the bound \eqref{property2}. For $I_1$, it is easy to see
\begin{align}
|\partial_{\xi_1}^\alpha\partial_{\xi_2}^\beta I_1| \lesssim  |\xi_1| |\xi_2|^{2N},\  \text{where} \ \ |\xi_1|\ll |\xi_2|, \quad \alpha,\beta=0,1. \label{3.24}
\end{align}
For $I_2$, note that $\theta(\xi_1,  \xi_2)-\theta(\xi_1, \xi_1+ \xi_2) \neq 0$ implies $  |\xi_1| \sim |\xi_2|$, then
\begin{align}
|\partial_{\xi_1}^\alpha\partial_{\xi_2}^\beta I_2| \lesssim |\xi_1| |\xi_2|^{2N} ,\  \quad \alpha,\beta=0,1.\label{3.25}
\end{align}
Also, a direct computation shows that the remainder $r_1(\xi_1,\xi_2)$ satisfies
\begin{align}
|\partial_{\xi_1}^\alpha\partial_{\xi_2}^\beta r_1(\xi_1,\xi_2)| \lesssim  \langle  \xi_1\rangle\langle  \xi_2\rangle^{-1} , \quad \alpha,\beta=0,1.\label{3.26}
\end{align}
Therefore, we conclude from  \eqref{2.5}--\eqref{3.26} that \eqref{bjestimate} holds for $j=1$.

Similarly,  using \eqref{B1}, \eqref{B2}, \eqref{Q1} and \eqref{Q2}, the  expressions for $b_4(\xi_1,\xi_2)$, $b_2(\xi_1,\xi_2)$, $b_3(\xi_1,\xi_2)$ in \eqref{B12} are
\begin{align}
b_4(\xi_1 , \xi_2) : &= \frac{ 2i \xi_1}{\langle\xi_1 + \xi_2\rangle^{2N} + \langle\xi_2\rangle^{2N}}\Big[-\langle\xi_1 + \xi_2\rangle^{2N}\frac{\xi_2 (\xi_1 + \xi_2)}{\langle \xi_1+ \xi_2\rangle} \theta(\xi_1,  \xi_2)\nonumber\\
&\qquad\qquad\qquad\qquad\qquad\quad +\langle\xi_2\rangle^{2N} \langle\xi_1 + \xi_2\rangle \theta(\xi_1, \xi_1+ \xi_2)\Big],\nonumber\\
b_2(\xi_1 , \xi_2)  : &= \frac{2i\langle\xi_1 \rangle}{\langle\xi_1 + \xi_2\rangle^{2N} + \langle\xi_2\rangle^{2N}}[\langle\xi_1 + \xi_2\rangle^{2N}\xi_2\theta(\xi_1,  \xi_2)-\langle  \xi_2\rangle^{2N} (\xi_1 + \xi_2)\theta(\xi_1, \xi_1+ \xi_2)],\nonumber\\
b_3(\xi_1 , \xi_2)  : &= \frac{2i\langle\xi_1 \rangle}{\langle\xi_1 + \xi_2\rangle^{2N} + \langle\xi_2\rangle^{2N}} \Big[\langle\xi_1 + \xi_2\rangle^{2N}\frac{(\xi_1 + \xi_2)\langle\xi_2 \rangle}{\langle\xi_1 + \xi_2\rangle}\theta(\xi_1,  \xi_2)\nonumber\\
&\qquad\qquad\qquad\qquad\qquad\quad- \langle \xi_2\rangle^{2N}\frac{ \xi_2\langle \xi_1+ \xi_2 \rangle}{\langle\xi_2\rangle}\theta(\xi_1, \xi_1+ \xi_2)\Big].\nonumber
\end{align}
Using the function $\chi(\xi_1,\xi_2)$ (see \eqref{definitionchi}), we  have
\begin{align}\label{reductionofb423}
\begin{split}
b_4(\xi_1 , \xi_2)  &= \frac{ -2i \xi_1(\xi_1+\xi_2)}{(\langle\xi_1 + \xi_2\rangle^{2N} + \langle\xi_2\rangle^{2N})\langle\xi_1+\xi_2\rangle
}\chi(\xi_1,\xi_2)+r_4(\xi_1,\xi_2),\\
b_2(\xi_1 , \xi_2)   &= \frac{2i\langle\xi_1 \rangle}{\langle\xi_1 + \xi_2\rangle^{2N} + \langle\xi_2\rangle^{2N}}\chi(\xi_1,\xi_2)+r_2(\xi_1,\xi_2),\\
b_3(\xi_1 , \xi_2)   &= \frac{2i\langle\xi_1 \rangle (\xi_1+\xi_2)\xi_2}{(\langle\xi_1 + \xi_2\rangle^{2N} + \langle\xi_2\rangle^{2N})\langle\xi_1+\xi_2\rangle\langle\xi_2\rangle}\chi(\xi_1,\xi_2)+r_3(\xi_1,\xi_2),
\end{split}
\end{align}
where
\begin{align*}
r_4(\xi_1 , \xi_2)  &= \frac{ 2i \xi_1\xi_2(\xi_1+\xi_2)\theta(\xi_1,\xi_2)}{(\langle\xi_1 + \xi_2\rangle^{2N} + \langle\xi_2\rangle^{2N})\langle\xi_1+\xi_2\rangle
}[(\xi_1+\xi_2)^{2N}-\langle\xi_1+\xi_2\rangle^{2N}]\\
&\quad+\frac{ 2i \xi_1\theta(\xi_1,\xi_1+\xi_2)}{(\langle\xi_1 + \xi_2\rangle^{2N} + \langle\xi_2\rangle^{2N})\langle\xi_1+\xi_2\rangle
}[\langle\xi_2\rangle^{2N}\langle\xi_1+\xi_2\rangle^2-\xi_2^{2N}(\xi_1+\xi_2)^2],\\
r_2(\xi_1 , \xi_2)  &= \frac{ 2i \langle\xi_1\rangle\xi_2\theta(\xi_1,\xi_2)}{\langle\xi_1 + \xi_2\rangle^{2N} + \langle\xi_2\rangle^{2N}}[\langle\xi_1+\xi_2\rangle^{2N}-(\xi_1+\xi_2)^{2N}]\\
&\quad-\frac{ 2i \langle\xi_1\rangle(\xi_1+\xi_2)\theta(\xi_1,\xi_1+\xi_2)}{\langle\xi_1 + \xi_2\rangle^{2N} + \langle\xi_2\rangle^{2N}}[\langle\xi_2\rangle^{2N}-\xi_2^{2N}],\\
r_3(\xi_1 , \xi_2)  &= \frac{ 2i\langle \xi_1\rangle(\xi_1+\xi_2)\theta(\xi_1,\xi_2)}{(\langle\xi_1 + \xi_2\rangle^{2N} + \langle\xi_2\rangle^{2N})\langle\xi_1+\xi_2\rangle\langle\xi_2\rangle
}[\langle\xi_1+\xi_2\rangle^{2N}\langle\xi_2\rangle^2-(\xi_1+\xi_2)^{2N}\xi_2^2]\\
&\quad-\frac{ 2i \langle\xi_1\rangle\xi_2\theta(\xi_1,\xi_1+\xi_2)}{(\langle\xi_1 + \xi_2\rangle^{2N} + \langle\xi_2\rangle^{2N})\langle\xi_1+\xi_2\rangle\langle\xi_2\rangle
}[\langle\xi_2\rangle^{2N}\langle\xi_1+\xi_2\rangle^2-\xi_2^{2N}(\xi_1+\xi_2)^2].
\end{align*}
With similar argument as above, we can obtain, for any $\alpha,\beta=0,1$ and $j =2, 3, 4$,
\begin{align}\label{3.26.2}
\begin{split}
|\partial_{\xi_1}^\alpha\partial_{\xi_2}^\beta r_j(\xi_1,\xi_2)| &\lesssim  \langle  \xi_1\rangle\langle  \xi_2\rangle^{-1},\\
|\partial_{\xi_1}^\alpha \partial_{\xi_2}^\beta b_j(\xi_1,\xi_2) | &\lesssim \langle\xi_1\rangle^2.
\end{split}
\end{align}
Therefore, the estimate \eqref{bjestimate} holds for all $j$, $j= 1, 2, 3, 4$.

Note that
\begin{align*}
\langle\partial_x\rangle^{N}\mathcal{O}[f,B_1]U &= \frac{1}{(2 \pi)^2}\langle\partial_x\rangle^{N}\int_{\mathbb{R}^2}  e^{i x (\xi_1 +\xi_2)}
 \widehat{f}(\xi_1)B_1(\xi_1, \xi_2) \widehat{U}(\xi_2) d\xi_1 d\xi_2\\
 &=\frac{1}{(2 \pi)^2}\int_{\mathbb{R}^2}  e^{i x (\xi_1 +\xi_2)}
 \widehat{\langle\partial_x\rangle^\rho f}(\xi_1) M(\xi_1, \xi_2) \widehat{\langle\partial_x\rangle^N U}(\xi_2) d\xi_1 d\xi_2
\end{align*}
with
$$
M(\xi_1,\xi_2)=\frac{\langle\xi_1+\xi_2\rangle^N}{\langle\xi_1\rangle^\rho\langle\xi_2\rangle^N}B_1(\xi_1,\xi_2), \  \ |\xi_1|\ll|\xi_2|.
$$
In view of \eqref{bjestimate}, we have for $\rho\geq 3$,
\begin{align*}
\|M(\xi_1,\xi_2-\xi_1)\|_{L^\infty_{\xi_2}H^1_{\xi_1}}
=\big\|\frac{\langle\xi_2\rangle^N}{\langle\xi_1\rangle^\rho\langle\xi_2-\xi_1\rangle^N}B_1(\xi_1,\xi_2-\xi_1)
\big\|_{L^\infty_{\xi_2}H^1_{\xi_1}}\lesssim \|\langle\xi_1\rangle^{2-\rho}\|_{L^2}\lesssim 1.
\end{align*}
So applying  Lemma \ref{l5.1} yields
$$
\|\langle\partial_x\rangle^{N}\mathcal{O}[f,B_1]U\|_{L^2}\lesssim \|f\|_{C^\rho}\|U\|_{H^N}.
$$
Similarly, we can  prove
$$
\|\langle\partial_x\rangle^{N}\mathcal{O}[f,B_2]U\|_{L^2}\lesssim \|f\|_{C^\rho}\|U\|_{H^N}.
$$
Hence,  \eqref{4.5} follows from these two estimates. This ends the proof of  Proposition \ref{p3.1}.
 $\hfill \Box$

\subsection{Energy normal form transformation}

For the equation \eqref{main},
\begin{align}
 U_t + D U =  \mathcal{O}[r, Q_1]U + \mathcal{O}[u, Q_2]U +  \mathcal{O}[r, S_1]U + \mathcal{O}[u, S_2]U, \label{4-1}
\end{align}
from Proposition \ref{p3.1}, we notice that the low-high term $\mathcal{O}[r,B_1]U+  \mathcal{O}[u, B_2]U $, which is a part of $\mathcal{O}[r, Q_1]U + \mathcal{O}[u, Q_2]U$,  will not lead to loss of derivatives. Now we can use Shatah's  normal form method to eliminate  this quadratic  term.

\begin{prop}\label{p4}
There exist two matrices $A_1$ and $A_2$ defined by
\begin{align}\label{A12}
A_1=\begin{pmatrix}
0 & a_1(\xi_1, \xi_2)  \\
a_4(\xi_1, \xi_2) &  0 \\
\end{pmatrix}, \quad\quad
A_2=\begin{pmatrix}
 a_2(\xi_1, \xi_2) &  0 \\
 0 &  a_3(\xi_1, \xi_2) \\
 \end{pmatrix}
\end{align}
such that
\begin{align}\label{2}
\begin{split}
&D\mathcal{O}[r,A_2]U - \mathcal{O}[\langle\partial_x\rangle r,A_1]U  - \mathcal{O}[r,A_2]DU  = -\mathcal{O}[r, B_1]U,\\
& D\mathcal{O}[u,A_1]U + \mathcal{O}[\langle\partial_x\rangle u,A_2]U  - \mathcal{O}[u,A_1]DU = -\mathcal{O}[u, B_2]U.
\end{split}
\end{align}
Moreover, for any $ \alpha,\beta=0,1$, we have
\begin{align}\label{ajestimate}
|\partial_{\xi_1}^\alpha\partial_{\xi_2}^\beta a_j(\xi_1,\xi_2) | \lesssim \langle\xi_1\rangle^3, \quad j=1,2,3,4.
\end{align}
\end{prop}

{\noindent \emph{Proof.}}
Inserting \eqref{A12} into \eqref{2}, we have
\begin{align}
\int_{\mathbb{R}^2}  &e^{i x (\xi_1 +\xi_2)} \widehat{r}(\xi_1) \Big[ \begin{pmatrix}
                                                                        0 &  -\langle\xi_1+\xi_2\rangle {a}_3(\xi_1, \xi_2) \\
                                                                         \langle\xi_1+\xi_2\rangle  {a}_2(\xi_1, \xi_2)& 0 \\
                                                                      \end{pmatrix}
                                                                      \nonumber\\
& \quad\quad\quad\quad\quad \quad  - \begin{pmatrix}
                                                                        0 &  \langle\xi_1\rangle  {a}_1(\xi_1, \xi_2) \\
                                                                         \langle\xi_1\rangle {a}_4(\xi_1, \xi_2)& 0 \\
                                                                      \end{pmatrix}
                                                                      \nonumber\\
& \quad\quad\quad\quad\quad \quad -\begin{pmatrix}
                                                                        0 & - \langle\xi_2\rangle{a}_2(\xi_1, \xi_2) \\
                                                                         \langle\xi_2\rangle  {a}_3(\xi_1, \xi_2)& 0 \\
                                                                      \end{pmatrix} \Big]\widehat{U}(\xi_2)d\xi_1 d\xi_2\nonumber\\
                                                                    =& \int_{\mathbb{R}^2}  e^{i x (\xi_1 +\xi_2)} \widehat{r}(\xi_1) \begin{pmatrix}
                                                                        0 & -{b_1}(\xi_1, \xi_2) \\
                                                                       -{b_4}(\xi_1, \xi_2) & 0 \\
                                                                      \end{pmatrix}
   \widehat{U}(\xi_2)
d\xi_1 d\xi_2.\nonumber
\end{align}
Similarly,
\begin{align}
\int_{\mathbb{R}^2}  &e^{i x (\xi_1 +\xi_2)} \widehat{u}(\xi_1) \Big[ \begin{pmatrix}
                                                                         -\langle\xi_1+\xi_2\rangle {a}_4(\xi_1, \xi_2) & 0 \\
                                                                         0  & \langle\xi_1+\xi_2\rangle {a}_1(\xi_1, \xi_2)\\
                                                                      \end{pmatrix}
                                                                      \nonumber\\
& \quad\quad\quad\quad\quad \quad  + \begin{pmatrix}
                                                                        \langle\xi_1\rangle {a}_2(\xi_1, \xi_2) & 0 \\
                                                                        0 & \langle\xi_1\rangle {a}_3(\xi_1, \xi_2) \\
                                                                      \end{pmatrix}
                                                                      \nonumber\\
& \quad\quad\quad\quad\quad \quad -\begin{pmatrix}
                                                                        \langle\xi_2\rangle{a}_1(\xi_1, \xi_2) & 0 \\
                                                                         0  & -\langle\xi_2\rangle  {a}_4(\xi_1, \xi_2) \\
                                                                      \end{pmatrix} \Big]\widehat{U}(\xi_2)d\xi_1 d\xi_2\nonumber\\
                                                                     =& \int_{\mathbb{R}^2}  e^{i x (\xi_1 +\xi_2)} \widehat{u}(\xi_1) \begin{pmatrix}
                                                                      - b_2(\xi_1, \xi_2) & 0 \\
                                                                     0  & - b_3(\xi_1, \xi_2) \\
                                                                      \end{pmatrix}
   \widehat{U}(\xi_2)
d\xi_1 d\xi_2.\nonumber
\end{align}
Thus we obtain linear equations
\begin{align}\label{linearequation}
\begin{pmatrix}
    -\langle\xi_1\rangle &\langle\xi_2\rangle  & - \langle\xi_1+\xi_2\rangle &  0\\
     0  & \langle\xi_1+\xi_2\rangle  & - \langle\xi_2\rangle  &  -\langle\xi_1\rangle\\
 - \langle\xi_2\rangle   &\langle\xi_1\rangle & 0 & - \langle\xi_1+\xi_2\rangle \\
\langle\xi_1+\xi_2\rangle &  0 & \langle\xi_1\rangle  &  \langle\xi_2\rangle  \\
\end{pmatrix} \begin{pmatrix}
                {a}_1 \\
               {a}_2 \\
               {a}_3 \\
                {a}_4 \\
              \end{pmatrix}
              =\begin{pmatrix}
               - b_1(\xi_1, \xi_2) \\
                - b_4(\xi_1, \xi_2) \\
                -  b_2(\xi_1, \xi_2) \\
                - b_3(\xi_1, \xi_2)
              \end{pmatrix}.
\end{align}
The solution of the above system is
\begin{align}\label{coofa}
\begin{split}
a_1=&\frac{1}{G}\big[(-\langle\xi_1\rangle^2+\langle\xi_2\rangle^2+\langle\xi_1+\xi_2\rangle^2)
\cdot(\langle\xi_1\rangle b_1-\langle\xi_2\rangle b_2
+\langle\xi_1+\xi_2\rangle b_3)\\
&\quad\quad\quad-2\langle\xi_2\rangle\langle\xi_1+\xi_2\rangle\cdot(\langle\xi_1\rangle b_4
+\langle\xi_2\rangle b_3-\langle\xi_1+\xi_2\rangle b_2
)\big],\\
a_2=&\frac{1}{G}\big[(\langle\xi_1\rangle^2-\langle\xi_2\rangle^2-\langle\xi_1+\xi_2\rangle^2)\cdot(\langle\xi_1\rangle b_2
-\langle\xi_2\rangle b_1
-\langle\xi_1+\xi_2\rangle b_4)\\
&\quad\quad\quad-2\langle\xi_2\rangle\langle\xi_1+\xi_2\rangle\cdot(\langle\xi_1\rangle b_3+\langle\xi_2\rangle b_4+\langle\xi_1+\xi_2\rangle b_1
)\big],\\
a_3=&\frac{1}{G}\big[(\langle\xi_1\rangle^2-\langle\xi_2\rangle^2-\langle\xi_1+\xi_2\rangle^2)\cdot(\langle\xi_1\rangle b_3
+\langle\xi_2\rangle b_4+\langle\xi_1+\xi_2\rangle b_1)\\
&\quad\quad\quad-2\langle\xi_2\rangle\langle\xi_1+\xi_2\rangle\cdot(\langle\xi_1\rangle b_2-\langle\xi_2\rangle b_1-\langle\xi_1+\xi_2\rangle b_4)\big],\\
a_4=&\frac{1}{G}\big[(-\langle\xi_1\rangle^2+\langle\xi_2\rangle^2+\langle\xi_1+\xi_2\rangle^2)
\cdot(\langle\xi_1\rangle b_4+\langle\xi_2\rangle b_3-
\langle\xi_1+\xi_2\rangle b_2)\\
&\quad\quad\quad-2\langle\xi_2\rangle\langle\xi_1+\xi_2\rangle\cdot(\langle\xi_1\rangle b_1- \langle\xi_2\rangle b_2
+\langle\xi_1+\xi_2\rangle b_3)\big],
\end{split}
\end{align}
where $G = 2\xi_1^2 + 2\xi_2^2 +2(\xi_1+\xi_2)^2 + 3   > 0$.

Now we prove \eqref{ajestimate}. To this end, we first claim that, for any $\alpha,\beta=0,1$,
\begin{align}
&|\partial^\alpha_{\xi_1}\partial^\beta_{\xi_2} (\langle\xi_1\rangle b_j(\xi_1,\xi_2))|\lesssim \langle\xi_1\rangle^3,\ j=1,2,3,4,\label{basic0}\\
&| \partial_{\xi_1}^\alpha \partial_{\xi_2}^\beta(\langle\xi_1+\xi_2\rangle b_3 (\xi_1 , \xi_2)-\langle\xi_2\rangle b_2 (\xi_1 , \xi_2))|\lesssim \langle\xi_1\rangle^3,\label{basic1}\\
&| \partial_{\xi_1}^\alpha\partial_{\xi_2}^\beta (\langle\xi_2\rangle b_3 (\xi_1 , \xi_2)-\langle\xi_1+\xi_2\rangle b_2 (\xi_1 , \xi_2))|\lesssim \langle\xi_1\rangle^3,\label{basic2}\\
&| \partial_{\xi_1}^\alpha \partial_{\xi_2}^\beta(\langle\xi_1+\xi_2\rangle b_1 (\xi_1 , \xi_2)+\langle\xi_2\rangle b_4 (\xi_1 , \xi_2))|\lesssim \langle\xi_1\rangle^3,\label{basic3}\\
&| \partial_{\xi_1}^\alpha\partial_{\xi_2}^\beta (\langle\xi_2\rangle b_1 (\xi_1 , \xi_2)+\langle\xi_1+\xi_2\rangle b_4 (\xi_1 , \xi_2))|\lesssim \langle\xi_1\rangle^3.\label{basic4}
\end{align}
Indeed, the bound \eqref{basic0} is a direct consequence of \eqref{bjestimate}. The proofs for \eqref{basic1}--\eqref{basic4} are similar, so we only show \eqref{basic1}.  From  \eqref{reductionofb423}, we see
\begin{align}\label{basic5}
\begin{split}
\langle\xi_1+\xi_2\rangle b_3 (\xi_1 , \xi_2)-\langle\xi_2\rangle b_2 (\xi_1 , \xi_2) &= \frac{2i\langle\xi_1 \rangle((\xi_1+\xi_2)\xi_2-\langle\xi_2\rangle^2)}{(\langle\xi_1 + \xi_2\rangle^{2N} + \langle\xi_2\rangle^{2N})\langle\xi_2\rangle}\chi(\xi_1,\xi_2) \\
&\quad+\langle\xi_1+\xi_2\rangle r_3(\xi_1 , \xi_2)- \langle\xi_2\rangle r_2(\xi_1 , \xi_2).
\end{split}
\end{align}
Remember that $|\xi_1|\ll |\xi_2|$. The bounds \eqref{3.24} and \eqref{3.25} imply
\begin{align}\label{basic6}
|\partial_{\xi_1}^\alpha \partial_{\xi_2}^\beta \chi(\xi_1,\xi_2)| \lesssim  \langle\xi_1\rangle \langle\xi_2\rangle^{2N},\quad \alpha,\beta=0,1.
\end{align}
Also, using  \eqref{3.26.2}, we have
 \begin{align}\label{basic7}
 |\partial_{\xi_1}^\alpha \partial_{\xi_2}^\beta( \langle\xi_1+\xi_2\rangle r_3(\xi_1 , \xi_2))|+  |\partial_{\xi_1}^\alpha \partial_{\xi_2}^\beta( \langle\xi_2\rangle r_2(\xi_1 , \xi_2))|\lesssim  \langle\xi_1\rangle ,\quad \alpha,\beta=0,1.
 \end{align}
 Inserting the bounds \eqref{basic6}--\eqref{basic7} into \eqref{basic5}, we can obtain  \eqref{basic1} as desired. Then from \eqref{coofa} and the bounds \eqref{basic0}--\eqref{basic4}, it is easy to see \eqref{ajestimate} holds. This completes the proof of Proposition \ref{p4}.  $\hfill \Box$

Similarly, we also use normal form method to cancel the high-high quadratic term $ \mathcal{O}[r, S_1]U + \mathcal{O}[u, S_2]U$. More precisely, we have
\begin{prop}\label{ppp4}
There exist two matrices $C_1$ and $C_2$ defined by
\begin{align}\label{C12}
C_1=\begin{pmatrix}
0 & c_1(\xi_1, \xi_2)  \\
c_4(\xi_1, \xi_2) &  0 \\
\end{pmatrix}, \quad\quad
C_2=\begin{pmatrix}
c_2(\xi_1, \xi_2) &  0 \\
 0 &  c_3(\xi_1, \xi_2) \\
 \end{pmatrix}
\end{align}
such that
\begin{align}\label{22}
\begin{split}
&D\mathcal{O}[r,C_2]U - \mathcal{O}[\langle\partial_x\rangle r,C_1]U  - \mathcal{O}[r,C_2]DU  = -\mathcal{O}[r, S_1]U,\\
& D\mathcal{O}[u,C_1]U + \mathcal{O}[\langle\partial_x\rangle u,C_2]U  - \mathcal{O}[u,C_1]DU = -\mathcal{O}[u, S_2]U.
\end{split}
\end{align}
Moreover, for any $ \alpha,\beta=0,1$, we have
\begin{align}\label{cjestimate}
|\partial_{\xi_1}^\alpha\partial_{\xi_2}^\beta c_j | \lesssim \langle\xi_1\rangle^3, \quad j=1,2,3,4.
\end{align}
\end{prop}

{\noindent \emph{Proof.}}  From \eqref{C12}--\eqref{22},  we obtain linear equations for $c_1$, $c_2$, $c_3$ and $c_4$
\begin{align*}
\begin{pmatrix}
    -\langle\xi_1\rangle &\langle\xi_2\rangle  & - \langle\xi_1+\xi_2\rangle &  0\\
     0  & \langle\xi_1+\xi_2\rangle  & - \langle\xi_2\rangle  &  -\langle\xi_1\rangle\\
 - \langle\xi_2\rangle   &\langle\xi_1\rangle & 0 & - \langle\xi_1+\xi_2\rangle \\
\langle\xi_1+\xi_2\rangle &  0 & \langle\xi_1\rangle  &  \langle\xi_2\rangle  \\
\end{pmatrix} \begin{pmatrix}
                {c}_1 \\
               {c}_2 \\
               { c}_3 \\
                {c}_4 \\
              \end{pmatrix}
              = \begin{pmatrix}
                s_1 \\
               s_4 \\
                s_2 \\
               s_3 \\
              \end{pmatrix},
\end{align*}
where the definitions of $s_1$, $s_2$, $s_3$ and $s_4$ are given by \eqref{S1}--\eqref{S2}. Clearly, the solution $(c_1,c_2,c_3,c_4)$ is given by replacing $b_j$ with $s_j$ ($j=1,2,3,4$) in \eqref{coofa}. Moreover, in the support of $1- \theta(\xi_1, \xi_2)-\theta(\xi_2,\xi_1)$, there holds $|\xi_1|\sim |\xi_2|$.  Hence, with similar argument as the proof of  Proposition \ref{p4}, the bound \eqref{cjestimate} can be obtained easily. $\hfill \Box$

Now we define the energy normal form transformation
\begin{align}\label{ph}
\Phi := U + \mathcal{O}[u,A_1]U + \mathcal{O}[r,A_2]U+  \mathcal{O}[u,C_1]U + \mathcal{O}[r,C_2]U ,
\end{align}
so that
 \begin{align}
\Phi_t + D\Phi &= U_t + D U + D\mathcal{O}[u,A_1]U + D\mathcal{O}[r,A_2]U+  D\mathcal{O}[u,C_1]U + D\mathcal{O}[r,C_2]U \nonumber\\
&\quad + \mathcal{O}[u_t,A_1]U +\mathcal{O}[u,A_1]U_t+ \mathcal{O}[r_t,A_2]U+   \mathcal{O}[r,A_2]U_t\nonumber\\
&\quad + \mathcal{O}[u_t,C_1]U +\mathcal{O}[u,C_1]U_t+ \mathcal{O}[r_t,C_2]U+   \mathcal{O}[r,C_2]U_t.\nonumber
\end{align}
Moreover, by  \eqref{main}, \eqref{ep3}, \eqref{2} and \eqref{22},  we have
\begin{align}
\Phi_t + D\Phi&= \mathcal{O}[r, Q_1]U  +D\mathcal{O}[r,A_2]U -\mathcal{O}[r,A_2]DU-\mathcal{O}[\langle\partial_x\rangle r,A_1]U\nonumber\\
&\quad +\mathcal{O}[u, Q_2]U+D\mathcal{O}[u,A_1]U -\mathcal{O}[u,A_1]DU+\mathcal{O}[\langle\partial_x\rangle u,A_2]U\nonumber\\
&\quad +\mathcal{O}[r, S_1]U +D\mathcal{O}[r,C_2]U -\mathcal{O}[r,C_2]DU-\mathcal{O}[\langle\partial_x\rangle r,C_1]U\nonumber\\
&\quad +\mathcal{O}[u, S_2]U+D\mathcal{O}[u,C_1]U -\mathcal{O}[u,C_1]DU+\mathcal{O}[\langle\partial_x\rangle u,C_2]U\nonumber\\
&\quad +(\mathcal{O}[u,A_1]+\mathcal{O}[r,A_2]+\mathcal{O}[u,C_1]+\mathcal{O}[r,C_2])
(\mathcal{O}[r,Q_1]U+\mathcal{O}[u,Q_2]U)\nonumber\\
&\quad +(\mathcal{O}[u,A_1]+\mathcal{O}[r,A_2]+\mathcal{O}[u,C_1]+\mathcal{O}[r,C_2])
(\mathcal{O}[r,S_1]U+\mathcal{O}[u,S_2]U)\nonumber\\
&\quad  + \mathcal{O}[\frac{\partial_x}{\langle\partial_x\rangle}( (\langle\partial_x\rangle u)^2 + (r_x)^2 ),A_1]U+ \mathcal{O}[ 2 \langle\partial_x\rangle u  \ r_x,A_2]U \nonumber\\
&\quad  + \mathcal{O}[\frac{\partial_x}{\langle\partial_x\rangle}((\langle\partial_x\rangle u)^2 + (r_x)^2 ),C_1]U+ \mathcal{O}[ 2 \langle\partial_x\rangle u  \ r_x,C_2]U\nonumber\\
&= \mathcal{O}[r, Q_1-B_1]U  +\mathcal{O}[u, Q_2-B_2]U\nonumber\\
&\quad +(\mathcal{O}[u,A_1]+\mathcal{O}[r,A_2]+\mathcal{O}[u,C_1]+\mathcal{O}[r,C_2])
(\mathcal{O}[r,Q_1]U+\mathcal{O}[u,Q_2]U)\nonumber\\
&\quad +(\mathcal{O}[u,A_1]+\mathcal{O}[r,A_2]+\mathcal{O}[u,C_1]+\mathcal{O}[r,C_2])
(\mathcal{O}[r,S_1]U+\mathcal{O}[u,S_2]U)\nonumber\\
&\quad  +I_1+I_2,\nonumber
\end{align}
where
\begin{align}
I_1&:=\mathcal{O}[\frac{\partial_x}{\langle\partial_x\rangle}( (\langle\partial_x\rangle u)^2 + (r_x)^2 ),A_1]U+ \mathcal{O}[ 2 \langle\partial_x\rangle u  \ r_x,A_2]U,\label{Term1} \\
I_2&:= \mathcal{O}[\frac{\partial_x}{\langle\partial_x\rangle}((\langle\partial_x\rangle u)^2 + (r_x)^2 ),C_1]U+ \mathcal{O}[ 2 \langle\partial_x\rangle u  \ r_x,C_2]U.\label{Term2}
\end{align}
Now using \eqref{ph}, we conclude that
\begin{align}
&\Phi_t + D \Phi - \mathcal{O}[r, Q_1-B_1]\Phi -\mathcal{O}[u, Q_2-B_2]\Phi= I_1 + I_2 + I_3 + I_4.\label{4.1}
\end{align}
Here, $I_1$, $I_2$ are defined by \eqref{Term1}, \eqref{Term2}, respectively, and
\begin{align}\label{Term3}
I_3&:= -\big[\mathcal{O}[r,Q_1]+\mathcal{O}[u,Q_2], \mathcal{O}[u,A_1]+\mathcal{O}[r,A_2]\big]U,
\end{align}
where the notation $[\cdot,\cdot]$ denotes the commutator, and
\begin{align}
\begin{split}\label{Term4}
I_4&:=  -(\mathcal{O}[r,Q_1]+\mathcal{O}[u,Q_2])(\mathcal{O}[u,C_1]+\mathcal{O}[r,C_2])U\\
&\quad+ (\mathcal{O}[u,C_1]+\mathcal{O}[r,C_2])(\mathcal{O}[r,Q_1]+\mathcal{O}[u,Q_2])U\\
&\quad+(\mathcal{O}[r,B_1] + \mathcal{O}[u,B_2])( \mathcal{O}[u,A_1]U+\mathcal{O}[r,A_2]U+\mathcal{O}[u,C_1]U+\mathcal{O}[r,C_2]U)\\
 &\quad +(\mathcal{O}[u,A_1]+\mathcal{O}[r,A_2]+\mathcal{O}[u,C_1]+\mathcal{O}[r,C_2])
(\mathcal{O}[r,S_1]U+\mathcal{O}[u,S_2]U).
\end{split}
\end{align}

\subsection{Energy estimate}

\begin{prop}\label{prope} Solutions of the  equation \eqref{4.1} satisfy
\begin{align*}
\frac{d}{ dt}\|\Phi\|_{H^N}^2 \lesssim  \|U\|_{C^5}^2\|U\|_{H^N}\|\Phi\|_{H^N}.
\end{align*}
\end{prop}

Proposition \ref{prope} will be proved by  Lemmas \ref{l5.2-1}--\ref{l5.2} and Lemma \ref{l7}.

\begin{lem}\label{l5.2-1} Let $A_1$ and $A_2$ be given by \eqref{A12}, then for any $\rho\geq 4$,
\begin{align}
&\|\langle\partial_x\rangle^{N}\mathcal{O}[f,A_1]U \|_{L^2}+\|\langle\partial_x\rangle^{N}\mathcal{O}[f,A_2]U  \|_{L^2}  \lesssim  \|f\|_{C^\rho}  \|U\|_{H^N}.\label{4.6}
\end{align}
\end{lem}

{\noindent \emph{Proof.}} The proof is similar to \eqref{4.5}.  By the definition of $\mathcal{O}[u,A_1]U $, we have
\begin{align*}
\langle\partial_x\rangle^{N}\mathcal{O}[f,A_1]U & = \frac{1}{(2 \pi)^2} \langle\partial_x\rangle^{N}\int_{\mathbb{R}^2}  e^{i x (\xi_1 +\xi_2)} \widehat{f}(\xi_1) A_1(\xi_1,\xi_2)\widehat{U}(\xi_2)d\xi_1 d\xi_2\\
&= \frac{1}{(2 \pi)^2} \int_{\mathbb{R}^2}  e^{i x (\xi_1 +\xi_2)} \widehat{\langle\partial_x\rangle^\rho f}(\xi_1) M(\xi_1,\xi_2) \widehat{\langle\partial_x\rangle^N U}(\xi_2)d\xi_1 d\xi_2,
\end{align*}
where
$$
M(\xi_1, \xi_2)= \frac{\langle\xi_1 + \xi_2\rangle ^N }{\langle\xi_1\rangle^{\rho}\langle\xi_2\rangle^N}A_1(\xi_1, \xi_2),\ |\xi_1|\ll |\xi_2|.$$
Note that \eqref{A12} and \eqref{ajestimate}  imply
\begin{align*}
\|M(\xi_1, \xi_2- \xi_1)\|_{L_{\xi_2}^\infty H_{\xi_1}^1} \lesssim \|\langle\xi_1\rangle^{3-\rho}\|_{L^2}\lesssim 1,\ \rho\geq 4.
\end{align*}
Using Lemma \ref{l5.1}, we thus obtain
$$
\|\langle\partial_x\rangle^{N}\mathcal{O}[f,A_1]U \|_{L^2} \lesssim  \|f\|_{C^\rho}  \|U\|_{H^N}.
$$
The estimate for $\|\mathcal{O}[f,A_2]U \|_{H^N}$  is the same.  $\hfill\Box$

As a direct consequence of Lemma \ref{l5.2-1}, we have
\begin{lem} For any $\alpha\in \mathbb{N}$,  we have
\begin{align}\label{4.8-2}
&\|\langle\partial_x\rangle^{N}\mathcal{O}[\langle\partial_x\rangle^\alpha f ,A_1]U \|_{L^2} + \|\langle\partial_x\rangle^{N}\mathcal{O}[\langle\partial_x\rangle^\alpha f ,A_2]U \|_{L^2} \lesssim \|f\|_{C^{\alpha +4}} \|U\|_{H^N}.
\end{align}
\end{lem}
Moreover, with the same argument as Lemma \ref{l5.2-1}, we can obtain the following lemma.
\begin{lem}\label{l5.6} Let $Q_1$ and $Q_2$ be defined by \eqref{Q1}--\eqref{Q2}, then for any $\rho \geq 2$, there holds
\begin{align}
&\|\langle\partial_x\rangle^{N}\mathcal{O}[f,Q_1]U \|_{L^2}+\|\langle\partial_x\rangle^{N}\mathcal{O}[f,Q_2]U \|_{L^2} \lesssim \|f\|_{C^{\rho}}  \|U\|_{H^{N+1}}.\label{l5.3}
\end{align}
\end{lem}

According to \eqref{property1},  the support  of $1- \theta(\xi_1, \xi_2)-\theta(\xi_2, \xi_1) $ satisfies $|\xi_1| \sim |\xi_2|$. Therefore,  we have the following results called ``derivative sharing" lemma.

\begin{lem}\label{l5.2} For any $\rho, \mu \in \mathbb{N}$, $\rho + \mu = N +2 $,  then
\begin{align}
&\|\langle\partial_x\rangle^{N}\mathcal{O}[f,C_1]U\|_{L^2} + \|\langle\partial_x\rangle^{N}\mathcal{O}[f,C_2]U\|_{L^2} \lesssim \|f\|_{C^{\rho+2}} \|U\|_{H^\mu},\label{4.7}\\
& \|\langle\partial_x\rangle^{N}\mathcal{O}[f,S_1]U\|_{L^2}+\|\langle\partial_x\rangle^{N}\mathcal{O}[f,S_2]U\|_{L^2}  \lesssim \|f\|_{C^{\rho+2}} \|U\|_{H^\mu}, \label{4.7.2}
\end{align}
where $C_1$, $C_2$ and $S_1$, $S_2$ are defined by \eqref{C12} and \eqref{S1}--\eqref{S2}, respectively.
\end{lem}
{\noindent\emph{ Proof.}} For $j=1,2$, we have
\begin{align}
\langle\partial_x\rangle^{N}\mathcal{O}[f,C_j]U&= \frac{1}{(2 \pi)^2}\langle\partial_x\rangle^{N}\int_{\mathbb{R}^2}  e^{i x (\xi_1 +\xi_2)} \widehat{f}(\xi_1)C_j(\xi_1, \xi_2) \widehat{U}(\xi_2) d\xi_1 d\xi_2\nonumber\\
&=  \frac{1}{(2 \pi)^2} \int_{\mathbb{R}^2}  e^{i x (\xi_1 +\xi_2)} \widehat{\langle\partial_x\rangle^{\rho+2} f}(\xi_1)\frac{\langle\xi_1 + \xi_2\rangle^N C_j(\xi_1, \xi_2)}{\langle\xi_1\rangle^{\rho+2} \langle\xi_2\rangle^{\mu} }\widehat{\langle\partial_x\rangle^\mu U}(\xi_2) d\xi_1 d\xi_2.\nonumber
\end{align}
Let
 $$
 M_j(\xi_1, \xi_2)= \frac{\langle\xi_1 + \xi_2\rangle^N C_j(\xi_1, \xi_2)}{\langle\xi_1\rangle^{\rho+2} \langle\xi_2\rangle^{\mu} },\ |\xi_1|\sim |\xi_2|,\ \ j=1,2,
 $$
from \eqref{C12} and \eqref{cjestimate}, it is easy to see
\begin{align}
\|M_j(\xi_1, \xi_2- \xi_1)\|_{L_{\xi_2}^\infty H_{\xi_1}^1}= \big\|\frac{\langle\xi_1 \rangle^N C_j(\xi_1, \xi_2 - \xi_1)}{\langle\xi_1\rangle^{\rho+2} \langle\xi_2-\xi_1\rangle^{\mu} }\big\|_{L_{\xi_2}^\infty H_{\xi_1}^1} \lesssim 1,\quad j=1,2.\label{5.8.1}
\end{align}
Hence, the desired estimate \eqref{4.7} follows from \eqref{5.8.1} and Lemma \ref{l5.1}. The proof for \eqref{4.7.2} is  similar as above.
$\hfill\Box$

\begin{lem} \label{lem3.8-2} Let
\begin{align*}
\mathcal{O}[f_1,f_2,M]V(x) :=  \frac{1}{(2 \pi)^3} \int_{\mathbb{R}^3} e^{i x(\xi + \eta +\sigma)}\widehat{f}_1(\xi) \widehat{f}_2(\eta)  M(\xi, \eta, \sigma) \widehat{V}(\sigma) d\xi  d \eta d\sigma,\nonumber
\end{align*}
 then we have
\begin{align*}
&\|\mathcal{O}[f_1,f_2,M]V\|_{L^2(\mathbb{R})}\lesssim \|   M(\xi, \eta- \xi, \sigma-\eta)\|_{L_\sigma^\infty H_{\xi}^1H^1_{\eta}}\|f_1\|_{L^\infty }\|f_2\|_{L^\infty }\|V\|_{L^2}.
\end{align*}
\end{lem}

This lemma can be proved by applying similar argument as the proof of  Lemma \ref{l5.1}.
The following two lemmas are  crucial in proving Lemma \ref{l7} below.

\begin{lem}\label{lem99}
Assume $|\xi_1|, |\eta|\ll |\xi_2|$, then for any $\alpha,\beta,\gamma=0,1$, we have
\begin{align}
|\partial_{\xi_1}^{\alpha}\partial_{\xi_2}^\beta(q_1(\xi_1,\xi_2)+q_4(\xi_1,\xi_2))|&\lesssim \langle\xi_1\rangle,\label{q1+4}\\
|\partial_{\xi_1}^{\alpha}\partial_{\xi_2}^\beta(q_2(\xi_1,\xi_2)-q_3(\xi_1,\xi_2))|&\lesssim \langle\xi_1\rangle,\label{q2-3}\\
|\partial^{\alpha}_{\eta}\partial_{\xi_1}^{\beta}\partial_{\xi_2}^\gamma(q_j(\eta,\xi_1+\xi_2)-q_j(\eta,\xi_2))|&\lesssim \langle\eta\rangle\langle\xi_1\rangle,\ j=1,2,3,4.\label{qj-j}
\end{align}
 \end{lem}

From the definitions \eqref{Q1}--\eqref{Q2}, one can easily  obtain the  bounds \eqref{q1+4}--\eqref{qj-j}.

\begin{lem}\label{lem999}
Assume $|\xi_1|, |\eta|\ll |\xi_2|$, then for any $\alpha,\beta,\gamma=0,1$, we have
\begin{align}
|\partial_{\xi_1}^{\alpha}\partial_{\xi_2}^\beta(a_1(\xi_1,\xi_2)+a_4(\xi_1,\xi_2))|&\lesssim \langle\xi_1\rangle^4\langle\xi_2\rangle^{-1},\label{a1+4}\\
|\partial_{\xi_1}^{\alpha}\partial_{\xi_2}^\beta(a_2(\xi_1,\xi_2)-a_3(\xi_1,\xi_2))|&\lesssim \langle\xi_1\rangle^4\langle\xi_2\rangle^{-1},\label{a2-3}\\
|\partial^{\alpha}_{\eta}\partial_{\xi_1}^{\beta}\partial_{\xi_2}^\gamma(a_j(\xi_1,\eta+\xi_2)-a_j(\xi_1,\xi_2))|&\lesssim \langle\eta\rangle\langle\xi_1\rangle^3\langle\xi_2\rangle^{-1},\ j=1,2,3,4.\label{aj-j}
\end{align}
\end{lem}

{\noindent \emph{Proof.}} It follows from \eqref{coofa} that
\begin{align*}
&a_1(\xi_1,\xi_2)+a_4(\xi_1,\xi_2)\\
&\quad=\frac{1}{G}(-\langle\xi_1\rangle^2+\langle\xi_2\rangle^2+\langle\xi_1
+\xi_2\rangle^2-2\langle\xi_2\rangle\langle\xi_1+\xi_2\rangle)\cdot(\langle\xi_1\rangle b_1-\langle\xi_2\rangle b_2
+\langle\xi_1+\xi_2\rangle b_3)\\
&\quad\quad +\frac{1}{G}(-\langle\xi_1\rangle^2+\langle\xi_2\rangle^2+\langle\xi_1
+\xi_2\rangle^2-2\langle\xi_2\rangle\langle\xi_1+\xi_2\rangle)\cdot(\langle\xi_1\rangle b_4+\langle\xi_2\rangle b_3
-\langle\xi_1+\xi_2\rangle b_2).
\end{align*}
Recall that $G=2\xi_1^2+2\xi_2^2+2(\xi_1+\xi_2)^2+3$. Now, using \eqref{basic0}--\eqref{basic2}, we can obtain \eqref{a1+4} as desired. The proof of \eqref{a2-3} is similar, so we skip it. By observing the structure of the expressions for $a_i$, we see that in order to prove  \eqref{aj-j}, it suffices to show
\begin{align*}
&|\partial_{\xi_1}^{\alpha}\partial_{\xi_2}^\beta(b_j(\xi_1,\eta+\xi_2)-b_j(\xi_1,\xi_2))|\lesssim \langle\eta\rangle\langle\xi_1\rangle^2\langle\xi_2\rangle^{-1},\ \ j=1,2,3,4,\\
&|\partial_{\eta}^\alpha\partial_{\xi_1}^\beta\partial_{\xi_2}^\gamma[\langle\xi_1+\xi_2+\eta\rangle
b_3(\xi_1,\xi_2+\eta)-\langle\xi_2+\eta\rangle b_2(\xi_1,\xi_2+\eta)\\
&\qquad\qquad\qquad-\langle\xi_1+\xi_2\rangle b_3(\xi_1,\xi_2)+\langle\xi_2\rangle b_2(\xi_1,\xi_2)]|\lesssim \langle\eta\rangle\langle\xi_1\rangle^3\langle\xi_2\rangle^{-1},\\
& |\partial_{\eta}^\alpha\partial_{\xi_1}^\beta\partial_{\xi_2}^\gamma[\langle\xi_2+\eta\rangle b_3(\xi_1,\xi_2+\eta)-\langle\xi_1+\xi_2+\eta\rangle b_2(\xi_1,\xi_2+\eta)\\
&\qquad\qquad\qquad-\langle\xi_2\rangle b_3(\xi_1,\xi_2)+\langle\xi_1+\xi_2\rangle b_2(\xi_1,\xi_2)]| \lesssim \langle\eta\rangle\langle\xi_1\rangle^3\langle\xi_2\rangle^{-1},
\end{align*}
and
\begin{align*}
&|\partial_{\eta}^\alpha\partial_{\xi_1}^\beta\partial_{\xi_2}^\gamma[\langle\xi_1+\xi_2+\eta\rangle
b_4(\xi_1,\xi_2+\eta)+\langle\xi_2+\eta\rangle b_1(\xi_1,\xi_2+\eta)\\
&\qquad\qquad\qquad-\langle\xi_1+\xi_2\rangle b_4(\xi_1,\xi_2)-\langle\xi_2\rangle b_1(\xi_1,\xi_2)]|\lesssim \langle\eta\rangle\langle\xi_1\rangle^3\langle\xi_2\rangle^{-1},\\
& |\partial_{\eta}^\alpha\partial_{\xi_1}^\beta\partial_{\xi_2}^\gamma[\langle\xi_2+\eta\rangle b_4(\xi_1,\xi_2+\eta)+\langle\xi_1+\xi_2+\eta\rangle b_1(\xi_1,\xi_2+\eta)\\
&\qquad\qquad\qquad-\langle\xi_2\rangle b_4(\xi_1,\xi_2)-\langle\xi_1+\xi_2\rangle b_1(\xi_1,\xi_2)]| \lesssim \langle\eta\rangle\langle\xi_1\rangle^3\langle\xi_2\rangle^{-1}.
\end{align*}
These estimates follow by  \eqref{2.5}, \eqref{reductionofb423} and an elementary  but tedious computation. We omit the details for simplicity.
$\hfill\Box$

\begin{lem}\label{l7} The following four commutator estimates hold:
\begin{align}
\|\big[\mathcal{O}[r,Q_1],\mathcal{O}[u,A_1]\big]U\|_{H^N} &\lesssim  \|r\|_{C^5}  \|u\|_{C^5} \|U\|_{H^N},\label{commutator1}\\
\|\big[\mathcal{O}[r,Q_1],\mathcal{O}[r,A_2]\big]U\|_{H^N} &\lesssim  \|r\|_{C^5}^2
 \|U\|_{H^N},\label{commutator2}\\
\|\big[\mathcal{O}[u,Q_2],\mathcal{O}[u,A_1]\big]U\|_{H^N}& \lesssim  \|u\|_{C^5}^2 \|U\|_{H^N},\label{commutator3}\\
\|\big[\mathcal{O}[u,Q_2],\mathcal{O}[r,A_2]\big]U\|_{H^N} &\lesssim  \|r\|_{C^5}  \|u\|_{C^5} \|U\|_{H^N}.\label{commutator4}
\end{align}
\end{lem}

\noindent{\emph{Proof.}}  We first show \eqref{commutator1}.  Note that
\begin{align*}
\mathcal{O}[r,Q_1](\mathcal{O}[u,A_1]U)&=\frac{1}{(2\pi)^3}\int_{\mathbb{R}^3}e^{ix(\eta+\xi_2)}\widehat{r}(\eta)Q_1(\eta,\xi_2)
\widehat{u}(\xi_1)A_1(\xi_1,\xi_2-\xi_1)\widehat{U}(\xi_2-\xi_1)d\eta d\xi_1d\xi_2\\
&=\frac{1}{(2\pi)^3}\int_{\mathbb{R}^3}e^{ix(\eta+\xi_1+\xi_2)}\widehat{r}(\eta)\widehat{u}(\xi_1)Q_1(\eta,\xi_1+\xi_2)
A_1(\xi_1,\xi_2)\widehat{U}(\xi_2)d\eta d\xi_1d\xi_2,
\end{align*}
and
\begin{align*}
\mathcal{O}[u,A_1](\mathcal{O}[r,Q_1]U)
&=\frac{1}{(2\pi)^3}\int_{\mathbb{R}^3}e^{ix(\xi_1+\xi_2)}\widehat{u}(\xi_1)A_1(\xi_1,\xi_2)
\widehat{r}(\eta)Q_1(\eta,\xi_2-\eta)\widehat{U}(\xi_2-\eta)d\eta d\xi_1d\xi_2\\
&=\frac{1}{(2\pi)^3}\int_{\mathbb{R}^3}e^{ix(\eta+\xi_1+\xi_2)}\widehat{r}(\eta)\widehat{u}(\xi_1)A_1(\xi_1,\xi_2+\eta)
Q_1(\eta,\xi_2)\widehat{U}(\xi_2)d\eta d\xi_1d\xi_2,
\end{align*}
Hence, we obtain
\begin{align*}
\big[\mathcal{O}[r,Q_1],\mathcal{O}[u,A_1]\big]U =\frac{1}{(2\pi)^3}\int_{\mathbb{R}^3}e^{ix(\eta+\xi_1+\xi_2)}\widehat{r}(\eta)\widehat{u}(\xi_1)M_1(\eta,\xi_1,\xi_2)
\widehat{U}(\xi_2)d\eta d\xi_1d\xi_2,
\end{align*}
where
\begin{align}\label{juanji1}
M_1(\eta,\xi_1,\xi_2):=Q_1(\eta,\xi_1+\xi_2)A_1(\xi_1,\xi_2)-A_1(\xi_1,\xi_2+\eta)Q_1(\eta,\xi_2).
\end{align}
From the support property of $Q_{1}$ and $A_1$,  we know that the support of $M(\eta, \xi_1, \xi_2)$ satisfies  $|\eta|, |\xi_1| \ll |\xi_2|$. Using this fact, in order to prove \eqref{commutator1}, it suffices to show
\begin{align}\label{temphaha2}
|\partial_{\eta}^\alpha \partial_{\xi_1}^\beta\partial_{\xi_2}^\gamma M_1(\eta,\xi_1,\xi_2)|&\lesssim \langle\xi_1\rangle^4\langle\eta\rangle^4, \quad \alpha,\beta,\gamma=0,\ 1.
\end{align}
Indeed, if \eqref{temphaha2} holds, we have
\begin{align*}
|\partial_{\eta}^\alpha \partial_{\xi_1}^\beta M_1(\eta,\xi_1-\eta,\xi_2-\xi_1)|&\lesssim \langle\xi_1-\eta\rangle^4\langle\eta\rangle^4,\ \alpha, \beta=0,1,
\end{align*}
then according to Lemma \ref{lem3.8-2}, the estimate  \eqref{commutator1} thus follows.

To prove \eqref{temphaha2}, we decompose the symbol $M_1$ into $M_{11}+M_{12}+M_{13}$ with
\begin{align*}
M_{11}(\eta,\xi_1,\xi_2)&:=(Q_1(\eta,\xi_1+\xi_2)-Q_1(\eta,\xi_2))A_1(\xi_1,\xi_2),\\
M_{12}(\eta,\xi_1,\xi_2)&:=(A_1(\xi_1,\xi_2)-A_1(\xi_1,\xi_2+\eta))Q_1(\eta,\xi_2),\\
M_{13}(\eta,\xi_1,\xi_2)&:=Q_1(\eta,\xi_2)A_1(\xi_1,\xi_2)-A_1(\xi_1,\xi_2)Q_1(\eta,\xi_2)\\
&=\begin{pmatrix}
1&  0 \\
 0 &  -1 \\
 \end{pmatrix}\cdot (q_1(\eta,\xi_2)a_4(\xi_1,\xi_2)-q_4(\eta,\xi_2)a_1(\xi_1,\xi_2)).
\end{align*}
For the symbol $M_{11}$, we use \eqref{ajestimate} and \eqref{qj-j}, then
$$
|\partial_{\eta}^\alpha \partial_{\xi_1}^\beta\partial_{\xi_2}^\gamma M_{11}(\eta,\xi_1,\xi_2)|\lesssim \langle\eta\rangle\langle\xi_1\rangle^4.
$$
For $M_{12}$, by \eqref{Q1} and \eqref{aj-j}, it is easy to see
$$
|\partial_{\eta}^\alpha \partial_{\xi_1}^\beta\partial_{\xi_2}^\gamma M_{12}(\eta,\xi_1,\xi_2)|\lesssim \langle\eta\rangle^2\langle\xi_1\rangle^3.
$$
For the last symbol $M_{13}$, we note that
\begin{align*}
q_1(\eta,\xi_2)a_4(\xi_1,\xi_2)-q_4(\eta,\xi_2)a_1(\xi_1,\xi_2)&=
(q_1(\eta,\xi_2)+q_4(\eta,\xi_2))a_4(\xi_1,\xi_2)\\
&\quad -q_4(\eta,\xi_2)(a_1(\xi_1,\xi_2)+a_4(\xi_1,\xi_2)),
\end{align*}
then it can be inferred from \eqref{q1+4}, \eqref{ajestimate}, \eqref{a1+4} and \eqref{Q1} that
$$
|\partial_{\eta}^\alpha \partial_{\xi_1}^\beta\partial_{\xi_2}^\gamma M_{13}(\eta,\xi_1,\xi_2)|\lesssim \langle\eta\rangle\langle\xi_1\rangle^4.
$$
Combing the above three bounds yield \eqref{temphaha2}. This finishes the proof of \eqref{commutator1}.

We then turn to show \eqref{commutator2}--\eqref{commutator4}. Notice that
\begin{align*}
\big[\mathcal{O}[r,Q_1],\mathcal{O}[r,A_2]\big]U =\frac{1}{(2\pi)^3}\int_{\mathbb{R}^3}e^{ix(\eta+\xi_1+\xi_2)}\widehat{r}(\eta)\widehat{r}(\xi_1)M_2(\eta,\xi_1,\xi_2)
\widehat{U}(\xi_2)d\eta d\xi_1d\xi_2,\\
\big[\mathcal{O}[u,Q_2],\mathcal{O}[u,A_1]\big]U =\frac{1}{(2\pi)^3}\int_{\mathbb{R}^3}e^{ix(\eta+\xi_1+\xi_2)}\widehat{u}(\eta)\widehat{u}(\xi_1)M_3(\eta,\xi_1,\xi_2)
\widehat{U}(\xi_2)d\eta d\xi_1d\xi_2,\\
\big[\mathcal{O}[u,Q_2],\mathcal{O}[r,A_2]\big]U :=\frac{1}{(2\pi)^3}\int_{\mathbb{R}^3}e^{ix(\eta+\xi_1+\xi_2)}\widehat{u}(\eta)\widehat{r}(\xi_1)M_4(\eta,\xi_1,\xi_2)
\widehat{U}(\xi_2)d\eta d\xi_1d\xi_2,
\end{align*}
where
\begin{align*}
M_2(\eta,\xi_1,\xi_2):=Q_1(\eta,\xi_1+\xi_2)A_2(\xi_1,\xi_2)-A_2(\xi_1,\xi_2+\eta)Q_1(\eta,\xi_2),\\
M_3(\eta,\xi_1,\xi_2):=Q_2(\eta,\xi_1+\xi_2)A_1(\xi_1,\xi_2)-A_1(\xi_1,\xi_2+\eta)Q_2(\eta,\xi_2),\\
M_4(\eta,\xi_1,\xi_2)=Q_2(\eta,\xi_1+\xi_2)A_2(\xi_1,\xi_2)-A_2(\xi_1,\xi_2+\eta)Q_2(\eta,\xi_2).
\end{align*}
Applying Lemma \ref{lem3.8-2} and repeating similar argument as proof of \eqref{temphaha2}, \eqref{commutator2}--\eqref{commutator4} can be proved as desired. Since the proof is very similar to the symbol \eqref{juanji1}, we omit further details.\hfill $\Box$\newline

 \noindent{ \emph {Proof of Proposition \ref{prope}.}}  Performing energy estimate at $H^N$ level for  \eqref{4.1}, we have
\begin{align}
&\mathrm{Re}\langle\langle\partial_x\rangle^N \Phi_t + \langle\partial_x\rangle^N  D \Phi ,\langle\partial_x\rangle^N\Phi\rangle - \mathrm{Re}\langle\langle\partial_x\rangle^N  \mathcal{O}[r, Q_1-B_1]\Phi, \langle\partial_x\rangle^N\Phi\rangle\nonumber\\
&\quad-\mathrm{Re}\langle\langle\partial_x\rangle^N  \mathcal{O}[u, Q_2-B_2]\Phi, \langle\partial_x\rangle^N\Phi\rangle=  \mathrm{Re} \langle  \langle\partial_x\rangle^N (I_1 + I_2 + I_3 + I_4),  \langle\partial_x\rangle^N \Phi\rangle, \label{ip}
\end{align}
where
\begin{align}\label{4.2}
\mathrm{Re}\langle\langle\partial_x\rangle^N \Phi_t, \langle\partial_x\rangle^N\Phi\rangle= \frac{1}{2}\frac{d}{ dt}\|\Phi\|_{H^N}^2, \qquad
\mathrm{Re} \langle \langle\partial_x\rangle^N D \Phi, \langle\partial_x\rangle^N\Phi\rangle = 0,
\end{align}
and from  \eqref{2.1}--\eqref{2.2},
\begin{align}
\mathrm{Re}\langle\langle\partial_x\rangle^N  \mathcal{O}[r, Q_1-B_1]\Phi, \langle\partial_x\rangle^N\Phi\rangle=0,
\quad \mathrm{Re}\langle\langle\partial_x\rangle^N  \mathcal{O}[u, Q_2-B_2]\Phi, \langle\partial_x\rangle^N\Phi\rangle=0. \label{4.3}
\end{align}
It remains to estimate the nonlinear terms in the right hand side of \eqref{ip}.

First, we consider $I_1$. From \eqref{Term1} and \eqref{4.8-2},
\begin{align}
|\mathrm{Re} \langle  \langle\partial_x\rangle^N I_1, \langle\partial_x\rangle^N\Phi\rangle|&\lesssim |\langle \langle\partial_x\rangle^N\mathcal{O}[\frac{\partial_x}{\langle\partial_x\rangle}( (\langle\partial_x\rangle u)^2 + (r_x)^2 ),A_1]U,  \langle\partial_x\rangle^N\Phi\rangle|\nonumber\\
&\quad +|\langle \langle\partial_x\rangle^N \mathcal{O}[ 2 \langle\partial_x\rangle u \ r_x,A_2]U,\langle\partial_x\rangle^N\Phi\rangle|\nonumber\\
&\lesssim (\|r^2\|_{C^5} +  \|u^2\|_{C^5}+\|ru\|_{C^5})\|U\|_{H^N} \|\Phi\|_{H^N}  \nonumber\\
& \lesssim (\|r\|_{C^5} +  \|u\|_{C^5})^2 \|U\|_{H^N}\|\Phi\|_{H^N}.\label{5.34}
\end{align}
Similarly, from \eqref{Term2} and \eqref{4.7},  there holds
\begin{align}
&\quad |\mathrm{Re} \langle  \langle\partial_x\rangle^N I_2 ,  \langle\partial_x\rangle^N\Phi\rangle|
 \lesssim (\|r\|_{C^5} +  \|u\|_{C^5})^2 \|U\|_{H^N}\|\Phi\|_{H^N}.\label{5.35}
\end{align}

Next,  we consider $I_4 $. Decompose $I_4$ (see \eqref{Term4}) into $I_{41}+I_{42}+I_{43}+I_{44}+I_{45}$ with
\begin{align*}
I_{41}&:= -(\mathcal{O}[r,Q_1]+\mathcal{O}[u,Q_2])(\mathcal{O}[u,C_1]+\mathcal{O}[r,C_2])U,\\
I_{42}&:= (\mathcal{O}[u,C_1]+\mathcal{O}[r,C_2])(\mathcal{O}[r,Q_1]+\mathcal{O}[u,Q_2])U,\\
I_{43}&:= (\mathcal{O}[r,B_1] + \mathcal{O}[u,B_2])( \mathcal{O}[u,A_1]U+\mathcal{O}[r,A_2]U),\\
I_{44}&:= (\mathcal{O}[r,B_1] + \mathcal{O}[u,B_2])(\mathcal{O}[u,C_1]U+\mathcal{O}[r,C_2]U),\\
I_{45}&:=(\mathcal{O}[u,A_1]+\mathcal{O}[r,A_2]+\mathcal{O}[u,C_1]+\mathcal{O}[r,C_2])
(\mathcal{O}[r,S_1]U+\mathcal{O}[u,S_2]U).
\end{align*}
From Lemma \ref{l5.6} and Lemma \ref{l5.2},
\begin{align}
|\mathrm{Re} \langle  \langle\partial_x\rangle^N I_{41} ,  \langle\partial_x\rangle^N\Phi\rangle|
 &\lesssim (\|r\|_{C^5} +  \|u\|_{C^5})\|(\mathcal{O}[u,C_1]+\mathcal{O}[r,C_2])U\|_{H^{N+1}} \|\Phi\|_{H^N},\nonumber\\
 &\lesssim  (\|r\|_{C^5} +  \|u\|_{C^5})^2\|U\|_{H^{N}} \|\Phi\|_{H^N},\\
|\mathrm{Re} \langle  \langle\partial_x\rangle^N I_{42} ,  \langle\partial_x\rangle^N\Phi\rangle|
 &\lesssim (\|r\|_{C^5} +  \|u\|_{C^5})\|(\mathcal{O}[r,Q_1]+\mathcal{O}[u,Q_2])U\|_{H^{N-1}} \|\Phi\|_{H^N},\nonumber \\
 &\lesssim  (\|r\|_{C^5} +  \|u\|_{C^5})^2\|U\|_{H^{N}} \|\Phi\|_{H^N}.
\end{align}
From  \eqref{4.5} and Lemma \ref{l5.2-1}, we have
\begin{align}
|\mathrm{Re} \langle  \langle\partial_x\rangle^N I_{43},  \langle\partial_x\rangle^N\Phi\rangle|&\lesssim  \|(\mathcal{O}[r,B_1] + \mathcal{O}[u,B_2])( \mathcal{O}[u,A_1]U+\mathcal{O}[r,A_2]U)\|_{H^N} \|\Phi\|_{H^N}\nonumber\\
&\lesssim (\|r\|_{C^5} +  \|u\|_{C^5}) \|\mathcal{O}[u,A_1]U+\mathcal{O}[r,A_2]U\|_{H^N} \|\Phi\|_{H^5}\nonumber\\
& \lesssim (\|r\|_{C^5} +  \|u\|_{C^5}\|)^2 \|U\|_{H^N}\|\Phi\|_{H^N}.
\end{align}
For the terms $I_{44}$ and $I_{45}$,  we use \eqref{4.5}, Lemma \ref{l5.2-1} and  Lemma \ref{l5.2} to obtain
\begin{align}
|\mathrm{Re} \langle  \langle\partial_x\rangle^N I_{44},  \langle\partial_x\rangle^N\Phi\rangle|&\lesssim
(\|r\|_{C^5} +  \|u\|_{C^5}) \|\mathcal{O}[u,C_1]U+\mathcal{O}[r,C_2]U\|_{H^N} \|\Phi\|_{H^N}\nonumber\\
& \lesssim (\|r\|_{C^5} +  \|u\|_{C^5})^2 \|U\|_{H^N}\|\Phi\|_{H^N},\\
|\mathrm{Re} \langle  \langle\partial_x\rangle^N I_{45},  \langle\partial_x\rangle^N\Phi\rangle|&\lesssim
(\|r\|_{C^5} +  \|u\|_{C^5} ) \|\mathcal{O}[r,S_1]U+\mathcal{O}[u,S_2]U\|_{H^N} \|\Phi\|_{H^N}\nonumber\\
& \lesssim (\|r\|_{C^5} +  \|u\|_{C^5} )^2 \|U\|_{H^N}\|\Phi\|_{H^N}.
\end{align}

At last, we consider the term $I_3$ (see \eqref{Term3}), which is a commutator operator. Indeed, applying Lemma \ref{l7}, we see
\begin{align}
& |\mathrm{Re} \langle \langle\partial_x\rangle^N I_3, \langle\partial_x\rangle^N \Phi \rangle | \lesssim (\|r\|_{C^5} +  \|u\|_{C^5})^2\|U\|_{H^N} \|\Phi\|_{H^N}. \label{5.47}
\end{align}

Now, combing the estimates \eqref{ip}--\eqref{5.47}, we obtain
\begin{align*}
\frac{d}{ dt}\|\Phi\|_{H^N}^2 \lesssim  \|U\|_{C^5}^2 \|U\|_{H^N}\|\Phi\|_{H^N}.
\end{align*}
This ends the proof of Proposition \ref{prope}.
$\hfill\Box$

Finally, we prove the energy estimate stated at the beginning of this section.
\newline

\noindent{\emph {Proof of Proposition \ref{eprop}.}}  It follows  from \eqref{ph}, Lemma \ref{l5.2-1} and Lemma \ref{l5.2} that
\begin{align*}
\|\Phi(t)\|_{H^N} \lesssim \|U(t)\|_{H^N} + \|U(t)\|_{C^5}\|U(t)\|_{H^N},
\end{align*}
and equally
\begin{align*}
\|U(t)\|_{H^N} \lesssim \|\Phi(t)\|_{H^N} + \|U(t)\|_{C^5}\|U(t)\|_{H^N}.
\end{align*}
Using \eqref{main a prior bound1}, we notice that if $\epsilon_1$ is sufficiently small, then
\begin{align*}
\|U(t)\|_{H^N} \lesssim \|\Phi(t)\|_{H^N} \lesssim \|U(t)\|_{H^N}.
\end{align*}
Hence, Proposition \ref{prope} and the a-priori bound \eqref{main a prior bound1} yield
$$
\frac{d}{dt}\|\Phi\|_{H^{N}}^2\lesssim \|U\|_{C^5}^2\|U\|_{H^N}^2\lesssim \epsilon_1^4(1+t)^{2p_0-1}.
$$
Integrating this estimate and using  \eqref{state1}, we deduce the desired bound \eqref{main energy bound1}.
$\hfill\Box$

\section{Low energy estimate for $\Gamma U$ and $xU$} \label{6}
In this section, we will prove the following two propositions, which lead to the low energy estimate ($H^{N_1}$ norm) of $\Gamma U=(x\partial_t+t\partial_x)U$ and $xU$, where $N_1\ll N$.

\begin{prop} \label{weprop} Let $U(t)\in C([0,T];H^N)$ be the solution of system \eqref{ep3}. Assume \eqref{state1} holds and \begin{align}\label{main a prior bound2}
\sup_{t\in [0,T]}\big[(1+t)^{-p_0}\|U(t)\|_{H^N}&+(1+t)^{-p_0}\|\Gamma U(t)\|_{H^{N_1}}\nonumber \\
&+(1+t)^{1/2}\|U(t)\|_{W^{N_1+10,\infty}}\big]\lesssim  \epsilon_1,
\end{align}
where $N_1=15$, $N=300$, $0<p_0<10^{-3}$ and $0<\epsilon_0\ll\epsilon_1\ll 1$, then
\begin{align}\label{main wenergy bound1}
\sup_{t\in [0,T]}\big[(1+t)^{-p_0}\|\Gamma U(t)\|_{H^{N_1}}\big]\lesssim \epsilon_0+ \epsilon_1^2.
\end{align}
\end{prop}

\begin{prop}\label{energyxuprop} Under the same assumptions as Proposition \ref{weprop}, we have
\begin{align}
\sup_{t\in [0,T]}\|xU(t)\|_{H^{N_1}}& \lesssim   \epsilon_1(1+ t)^{1+ p_0}.   \label{energyxu}
\end{align}
\end{prop}

In the following, Sections 3.1-3.3 are devoted to proving Proposition \ref{weprop}, and  Proposition \ref{energyxuprop} is proved in Section 3.4.

\subsection{Shatah's normal form for  quadratic terms without loss of derivatives
 }

To prove Proposition \ref{weprop}, we have to derive the equation for $\Gamma U$. Recall the Euler-Poisson system
\begin{align}\label{ff}
\begin{pmatrix}
           r \\
           u  \\
         \end{pmatrix}_t+  \begin{pmatrix}
           0&-\langle\partial_x\rangle \\
           \langle\partial_x\rangle&0  \\
         \end{pmatrix}
         \begin{pmatrix}
           r \\
           u  \\
         \end{pmatrix}
         =\begin{pmatrix}
           f_1 \\
           f_2  \\
         \end{pmatrix}
         :=\begin{pmatrix}
           2 \langle\partial_x\rangle u  \ r_x \\
           \frac{\partial_x}{\langle\partial_x\rangle}[ (\langle\partial_x\rangle u)^2 + (r_x)^2 ]  \\
         \end{pmatrix}.
\end{align}
For simplicity, we write the above system as
\begin{align}\label{ep32}
U_t + DU =(f_1,f_2)^T.
\end{align}
Operating $\Gamma$ on both sides of the system \eqref{ff}, then using  the relations
\begin{align*}
[\Gamma,\partial_x]=-\partial_t,\ \ \
[\Gamma,\partial_t]=-\partial_x,\ \ \ [\Gamma,\langle\partial_x\rangle]=\frac{\partial_x}{\langle\partial_x\rangle}\partial_t,\ \ \ [\Gamma,\frac{\partial_x}{\langle\partial_x\rangle}]=-\frac{1}{\langle\partial_x\rangle^3}\partial_t,
\end{align*}
 we obtain equations for $\Gamma r, \Gamma u$
\begin{align*}
\begin{pmatrix}
           \Gamma r \\
           \Gamma u  \\
         \end{pmatrix}_t+  \begin{pmatrix}
           0&-\langle\partial_x\rangle \\
           \langle\partial_x\rangle&0  \\
         \end{pmatrix}
         \begin{pmatrix}
           \Gamma r \\
           \Gamma u  \\
         \end{pmatrix}=\begin{pmatrix}
           2 \langle\partial_x\rangle (\Gamma u)  \ r_x+2 \langle\partial_x\rangle u  \ (\Gamma r)_x \\
           \frac{2\partial_x}{\langle\partial_x\rangle}[ \langle\partial_x\rangle (\Gamma u)\langle\partial_x\rangle u + (\Gamma r)_xr_x ]  \\
         \end{pmatrix}+\begin{pmatrix}
           g'_1 \\
           g'_2  \\
         \end{pmatrix}+\begin{pmatrix}
           g''_1 \\
           g''_2  \\
         \end{pmatrix},
\end{align*}
where $g'_1=g'_1( r, u)$, $g'_2=g'_2( r, u) $ are quadratic terms  without containing $\Gamma r$ and $\Gamma u$,
\begin{align*}
g'_1( r, u):=&  -2(r_x)^2-2(\langle\partial_x\rangle u)^2+\frac{\partial_x^2}{\langle\partial_x\rangle^2}[(r_x)^2+(\langle\partial_x\rangle u)^2],\\
g'_2( r, u):=& -\frac{6\partial_x}{\langle\partial_x\rangle}(r_x\langle\partial_x\rangle u)
+\frac{2}{\langle\partial_x\rangle^3}(\langle\partial_x\rangle^2r\langle\partial_x\rangle u)-\frac{2}{\langle\partial_x\rangle^3}(r_x\langle\partial_x\rangle u_x ),
\end{align*}
and $g''_1=g''_1( r, u)$, $g''_2=g''_2( r, u) $ are cubic terms,
\begin{align*}
g''_1( r, u):=&  -4 r_x( \langle\partial_x\rangle u)^2 +2 r_x \frac{\partial_x^2}{\langle\partial_x\rangle^2}[(r_x)^2+(\langle\partial_x\rangle u)^2],\\
g''_2( r, u)=&-\frac{4\partial_x}{\langle\partial_x\rangle}[(r_x)^2\langle\partial_x\rangle u]+\frac{2\partial_x}{\langle\partial_x\rangle}\left[\frac{\partial_x^2}{\langle\partial_x\rangle^2}[(r_x)^2
+(\langle\partial_x\rangle u)^2]\langle\partial_x\rangle u\right]\\
&-\frac{4}{\langle\partial_x\rangle^3}[(r_x \langle\partial_x\rangle u)_xr_x]
-\frac{4}{\langle\partial_x\rangle^3}[(r_xr_{xx}+\langle\partial_x\rangle u\langle\partial_x\rangle u_x)\langle\partial_x\rangle u].
\end{align*}
Since our aim is to estimate $\|\Gamma U\|_{H^{N_1}}$ (recall that $N_1\ll N$),  we simply decompose the quadratic terms  $g'_1$, $g'_2 $  into
\begin{align*}
(g'_1,g'_2)^T=\mathcal{O}\big[u,\begin{pmatrix}
0&  \widetilde{q}_1(\xi_1, \xi_2)\\
\widetilde{q}_4(\xi_1, \xi_2)&0
\end{pmatrix}\big]U+\mathcal{O}\big[r,\begin{pmatrix}
\widetilde{q}_2(\xi_1, \xi_2)  & 0\\
0 & \widetilde{q}_3(\xi_1, \xi_2)
\end{pmatrix}\big]U
\end{align*}
with
\begin{align*}
\widetilde{q}_1(\xi_1, \xi_2)&:=-2\langle\xi_1\rangle\langle\xi_2\rangle-\frac{(\xi_1+\xi_2)^2}
{\langle\xi_1+\xi_2\rangle^2}\langle\xi_1\rangle\langle\xi_2\rangle,\\
\widetilde{q}_4(\xi_1,\xi_2)&:=\frac{3(\xi_1+\xi_2)\langle\xi_1\rangle\xi_2}{\langle\xi_1+\xi_2\rangle}
+\frac{\langle\xi_1\rangle\langle\xi_2\rangle^2}{\langle\xi_1+\xi_2\rangle^3}
+\frac{\xi_1\langle\xi_1\rangle\xi_2}{\langle\xi_1+\xi_2\rangle^3},\\
\widetilde{q}_2(\xi_1, \xi_2)&:=2\xi_1\xi_2+\frac{(\xi_1+\xi_2)^2}
{\langle\xi_1+\xi_2\rangle^2}\xi_1\xi_2,\\
\widetilde{q}_3(\xi_1, \xi_2)&:=\frac{3(\xi_1+\xi_2)\xi_1\langle\xi_2\rangle}{\langle\xi_1+\xi_2\rangle}
+\frac{\langle\xi_1\rangle^2\langle\xi_2\rangle}{\langle\xi_1+\xi_2\rangle^3}
+\frac{\xi_1\xi_2\langle\xi_2\rangle}{\langle\xi_1+\xi_2\rangle^3}.
\end{align*}
For the terms including $\Gamma r$ and $\Gamma u$, we use similar decomposition as in Section 2.1,
\begin{align*}
\begin{pmatrix}
           2 \langle\partial_x\rangle (\Gamma u)  \ r_x+2 \langle\partial_x\rangle u  \ (\Gamma r)_x \\
           \frac{2\partial_x}{\langle\partial_x\rangle}[ \langle\partial_x\rangle (\Gamma u)\langle\partial_x\rangle u + (\Gamma r)_xr_x ]  \\
         \end{pmatrix}=\begin{pmatrix}
           F'_1 \\
           F'_2  \\
         \end{pmatrix}+\begin{pmatrix}
           F''_1 \\
           F''_2  \\
         \end{pmatrix}
\end{align*}
with $F'_j=F_j(\Gamma r, \Gamma u, r, u)$ and $F''_j= F'_j(\Gamma r, \Gamma u, r, u)$ ($j=1,2$)  defined by
\begin{align}
(F'_1, F'_2)^T :&= \mathcal{O}(r, Q_1) \Gamma U + \mathcal{O}(u, Q_2)\Gamma U +  \mathcal{O}(r, S_1)\Gamma U  + \mathcal{O} (u, S_2)\Gamma U,\label{5.3}\\
(F''_1, F''_2)^T:&= \mathcal{O}(\Gamma r, Q_1)U + \mathcal{O} (\Gamma u, Q_2)U +  \mathcal{O}(\Gamma r, S_1)U  + \mathcal{O} (\Gamma u, S_2)U,\label{5.3.1}
\end{align}
where the matrices $Q_1$, $Q_2$, $S_1$ and $S_2$ are given in \eqref{Q1}--\eqref{S2}. In conclusion, we obtain
\begin{align}\label{zrzu}
(\Gamma U)_t + D\Gamma U = (F'_1, F'_2)^T + (F''_1, F''_2)^T + (g'_1,g'_2)^T+(g''_1,g''_2)^T.
\end{align}

From \eqref{5.3.1}, we see the quadratic terms $F''_1(\Gamma r, \Gamma u, r, u)$ and  $F''_2(\Gamma r, \Gamma u, r, u)$ will not lead to loss of derivatives when taking the $H^{N_1}$ norm energy estimate, since $\Gamma r$ and $\Gamma u$ have lower frequencies  compared to $U$. Notice also that  the quadratic terms $g'_1( r, u)$ and $g'_2( r, u)$ do not contain $\Gamma r$ and $\Gamma u$.  For these reasons,  we only need to take Shatah's normal form transformation
\begin{align} \label{5.7}
  \widetilde{\Omega }  := &\Gamma U+\mathcal{O}[r, G_1]U + \mathcal{O}[u, G_2]U+\mathcal{O}[\Gamma u,  H_1]U+\mathcal{O}[\Gamma r, H_2]U
\end{align}
 for the system \eqref{zrzu} to cancel $g'_1( r, u)$, $g'_2( r, u)$,  $F''_1(\Gamma r, \Gamma u, r, u)$ and  $F''_2(\Gamma r, \Gamma u, r, u)$ .
Similar to \eqref{2}, the matrices $G_1$, $G_2$ could be obtained from  the equations
\begin{align}\label{normal11}
\begin{split}
-(g_{1}',g_{2}')^T&= D\mathcal{O}[r, G_1]U - \mathcal{O}[\langle\partial_x\rangle r,G_2]U  - \mathcal{O}[r, G_1]DU \\
         &\quad+ D\mathcal{O}[u, G_2]U + \mathcal{O}[\langle\partial_x\rangle u,G_1]U  - \mathcal{O}[u, G_2]DU ,
\end{split}
\end{align}
and $H_1$, $H_2$ can be determined by
\begin{align}\label{normal12}
\begin{split}
-(F''_1, F''_2)^T&=D\mathcal{O}[\Gamma u, H_1]U +\mathcal{O}[\langle\partial_x\rangle \Gamma u,H_2]U  - \mathcal{O}[\Gamma u, H_1]DU \\
        &\quad+ D\mathcal{O}[\Gamma r, H_2]U -\mathcal{O}[\langle\partial_x\rangle \Gamma r,H_1]U  - \mathcal{O}[\Gamma r, H_2]DU.
\end{split}
\end{align}
Indeed, the elements of $G_1$, $G_2$ (or $H_1$, $H_2$) satisfy similar linear equations as \eqref{linearequation}, which can be uniquely solved as \eqref{coofa}. Now using \eqref{ff}, \eqref{ep32} and \eqref{zrzu}--\eqref{normal12},  we reduce \eqref{zrzu}  to
\begin{align}
          \widetilde{\Omega}_t +D\widetilde{\Omega }
           = (F'_1, F'_2)^T  +(g_1,g_2)^T,\label{zrzu1}
\end{align}
where $(F'_1, F'_2)^T$ is given as \eqref{5.3}, and  $g_1( r, u, \Gamma r, \Gamma u)$, $g_2( r, u, \Gamma r, \Gamma u)$ are  cubic terms taking the following form
\begin{align}\label{gg}
\begin{split}
(g_1,g_2)^T := &(\mathcal{O}[f_1, G_1] + \mathcal{O}[f_2, G_2])U+ (\mathcal{O}[r, G_1]+ \mathcal{O}[u, G_2])(f_1,f_2)^T\\
&+\mathcal{O}[F'_2+F''_2 + g'_2+g''_2, H_1]U +\mathcal{O}[F'_1 + F''_1 + g'_1+g''_1, H_2]U\\
&+(\mathcal{O}[\Gamma u, H_1)]+ \mathcal{O}[\Gamma r, H_2])(f_1,f_2)^T+(g''_1,g''_2)^T.
\end{split}
\end{align}
Moreover, according to the properties of the symbols $G_1$, $G_2$, $H_1$ and $H_2$, we clearly have
\begin{align}\label{GHnorm}
\begin{split}
\|\mathcal{O}[f,G_1]V\|_{H^{N_1}}+\|\mathcal{O}[f,G_2]V\|_{H^{N_1}}&\lesssim \|f\|_{C^5}\|V\|_{H^N}+\|f\|_{H^{N}}\|V\|_{C^5},\\
\|\mathcal{O}[f,H_1]V\|_{H^{N_1}}+\|\mathcal{O}[f,H_2]V\|_{H^{N_1}}&\lesssim \|f\|_{H^{5}}\|V\|_{C^{N_1+5}}.
\end{split}
\end{align}

\subsection{Energy normal form for  quadratic terms with $\Gamma r$ and $\Gamma u$
 }
 Note that  the quadratic terms $F'_1(\Gamma r, \Gamma u, r, u)$ and  $F'_2(\Gamma r, \Gamma u, r, u)$ in \eqref{zrzu1} will  lead to loss of derivatives in energy estimate for $\Gamma U$, as in these terms $\Gamma U$  has higher frequencies  compared to $U$. So in this subsection, we apply similar modified normal form process as in Section 2  to eliminate the derivative quadratic terms ${F}'_1$ and ${F}'_2$ in \eqref{zrzu1}. Taking the energy normal form transformation
\begin{align}\label{6.9}
\Omega:= & \widetilde{\Omega} + \mathcal{O}[u,A_1]\Gamma U + \mathcal{O}[r,A_2]\Gamma U + \mathcal{O}[u,C_1]\Gamma U+ \mathcal{O}[r,C_2]\Gamma U,
\end{align}
where $A_1,  A_2, C_1, C_2$ are completely the same as those defined in \eqref{A12} and \eqref{C12}. By  repeating  similar process as \eqref{4.1}, we  can obtain the equation for $\Omega$
\begin{align}\label{equationomega}
\begin{split}
\Omega_t + D\Omega&= \mathcal{O}[r,Q_1-B_1]\Gamma U + \mathcal{O}[u,Q_2-B_2]\Gamma U+ (g_1,g_2)^T\\
         &\quad+ (\mathcal{O}[f_2,A_1]+ \mathcal{O}[f_1,A_2] )\Gamma U+(\mathcal{O}[f_2,C_1]+ \mathcal{O}[f_1,C_2] )\Gamma U\\
         &\quad+ (\mathcal{O}[u,A_1]+ \mathcal{O}[r,A_2]) (F'_1+ F''_1 +g'_1+g''_1, F'_2 + F''_2 +g'_2+g''_2)^T \\
         &\quad+ (\mathcal{O}[u,C_1]+ \mathcal{O}[r,C_2]) (F'_1 + F''_1 +g'_1+g''_1, F'_2 + F''_2 +g'_2+g''_2)^T,
       \end{split}
       \end{align}
where $f_1, f_2$, $F'_1, F'_2$, $F''_1, F''_2$ and $g_1, g_2$ are defined by \eqref{ff}, \eqref{5.3}, \eqref{5.3.1} and \eqref{gg}, respectively. From \eqref{5.7} and \eqref{6.9}, we have
\begin{align}\label{6.11}
\begin{split}
            \Omega = \Gamma U&+ \mathcal{O}[r, G_1]U + \mathcal{O}[u, G_2]U+\mathcal{O}[\Gamma u, H_1]U+\mathcal{O}[\Gamma r, H_2]U\\
            &+ \mathcal{O}[u,A_1]\Gamma U + \mathcal{O}[r,A_2]\Gamma U + \mathcal{O}[u,C_1]\Gamma U+ \mathcal{O}[r,C_2]\Gamma U.
\end{split}
\end{align}
Using \eqref{6.11},  we rewrite the  equation \eqref{equationomega}  as
 \begin{align}
        \Omega_t&+ D\Omega-\mathcal{O}[r,Q_1-B_1]\Omega - \mathcal{O}[u,Q_2-B_2]\Omega= J_1 + J_2 + J_3 + J_4  + (g_1,g_2)^T , \label{5.13}
\end{align}
where
\begin{align*}
J_1:=& \big[\mathcal{O}[u,A_1]+ \mathcal{O}[r,A_2], \mathcal{O}[r,Q_1]+ \mathcal{O}[u,Q_2]\big] \Gamma U,\\
J_2:=& -(\mathcal{O}[r,Q_1-B_1]+\mathcal{O}[u,Q_2-B_2])(\mathcal{O}[r, G_1]U + \mathcal{O}[u, G_2]U)\\
&+(\mathcal{O}[u,A_1]+ \mathcal{O}[r,A_2]+\mathcal{O}[u,C_1]+ \mathcal{O}[r,C_2])  (  g'_1+g''_1,  g'_2+g''_2)^T,\\
J_3:=&  (\mathcal{O}[r,B_1]+ \mathcal{O}[u,B_2])(\mathcal{O}[u,A_1] + \mathcal{O}[r,A_2]+ \mathcal{O}[u,C_1]+ \mathcal{O}[r,C_2])\Gamma U \\
&-(\mathcal{O}[r,Q_1]+ \mathcal{O}[u,Q_2])( \mathcal{O}[u,C_1]+ \mathcal{O}[r,C_2])\Gamma U+ ( \mathcal{O}[u,C_1]+ \mathcal{O}[r,C_2])(F'_1, F'_2)^T\nonumber\\
&+(\mathcal{O}[u,A_1]+ \mathcal{O}[r,A_2])( \mathcal{O}[u,S_1]+ \mathcal{O}[r,S_2])\Gamma U\\
&+ (\mathcal{O}[f_2,A_1]+ \mathcal{O}[f_1,A_2] )\Gamma U+(\mathcal{O}[f_2,C_1]+ \mathcal{O}[f_1,C_2] )\Gamma U,\\
J_4: =&   -(\mathcal{O}[r,Q_1-B_1]+\mathcal{O}[u,Q_2-B_2])(\mathcal{O}[\Gamma u, H_1]U + \mathcal{O}[\Gamma r, H_2]U)\\
&+(\mathcal{O}[u,A_1]+ \mathcal{O}[r,A_2]+\mathcal{O}[u,C_1]+ \mathcal{O}[r,C_2])  ( F''_1,  F''_2)^T.
\end{align*}

\subsection{Low energy estimate for $\Gamma U$}

{\noindent \emph{Proof of Proposition \ref{weprop}.}}  The energy estimate for  \eqref{5.13} is similar to \eqref{4.1}. Taking energy estimate at $H^{N_1}$ level  for \eqref{5.13}, we have
\begin{align*}
&\mathrm{Re} \langle\langle\partial_x\rangle^{N_1} \Omega_t + \langle\partial_x\rangle^{N_1} D \Omega  - \langle\partial_x\rangle^{N_1} (\mathcal{O}[r, Q_1-B_1]-\mathcal{O}[u, Q_2-B_2])\Omega, \langle\partial_x\rangle^{N_1}\Omega\rangle\\
&\quad\quad=  \mathrm{Re} \langle \langle\partial_x\rangle^{N_1} (J_1 + J_2 + J_3+J_4 )+\mathrm{Re} \langle \langle\partial_x\rangle^{N_1}(g_1,g_2)^T,  \langle\partial_x\rangle^{N_1}\Omega\rangle.
\end{align*}
Clearly, there holds
\begin{align*}
\mathrm{Re} \langle \langle\partial_x\rangle^{N_1} \Omega_t,  \langle\partial_x\rangle^{N_1}\Omega\rangle = \frac{d}{ 2dt}\|\Omega\|_{H^{N_1}}^2, \ \ \ \
\mathrm{Re} \langle  \langle\partial_x\rangle^{N_1} D \Omega,  \langle\partial_x\rangle^{N_1}\Omega\rangle =  0.
\end{align*}
Moreover, from \eqref{2.1}--\eqref{2.2},
\begin{align*}
\mathrm{Re} \langle  \langle\partial_x\rangle^{N_1} (\mathcal{O}[r, Q_1-B_1]+\mathcal{O}[u, Q_2-B_2]\Omega ), \langle\partial_x\rangle^{N_1} \Omega\rangle = 0.
\end{align*}
For $J_1$, we  use Lemma \ref{l7} to obtain
\begin{align}
|\langle  \langle\partial_x\rangle^{N_1} J_1,  \langle\partial_x\rangle^{N_1}\Omega\rangle| \lesssim \|U\|_{C^5}^2 \|\Gamma U\|_{H^{N_1}}\|\Omega\|_{H^{N_1}}.\nonumber
\end{align}
Note that all the terms in $J_2$ are cubic terms containing only $r, u$, from \eqref{4.5}, Lemmas \ref{l5.2-1}--\ref{l5.2} and \eqref{GHnorm}, we have
\begin{align}
| \langle  \langle\partial_x\rangle^{N_1} J_2,  \langle\partial_x\rangle^{N_1}\Omega\rangle| \lesssim \|U\|_{C^5}^2\|U\|_{H^{N}}\|\Omega\|_{H^{N_1}}.\nonumber
\end{align}
 Similar to the estimates for $I_1$, $I_2$ and $I_4$ in \eqref{4.1},  the term $J_3$ can be bounded by
\begin{align}
|\langle\partial_x\rangle^{N_1} J_3,  \langle\partial_x\rangle^{N_1}\Omega\rangle| \lesssim \|U\|_{C^{5}}^2  \|\Gamma U\|_{H^{N_1}}  \|\Omega\|_{H^{N_1}}.\nonumber
\end{align}
Note that all the terms containing $\Gamma U$ in $J_4$ has lower frequencies compared to $U$, so
\begin{align}
|\langle\partial_x\rangle^{N_1} J_4,  \langle\partial_x\rangle^{N_1}\Omega\rangle| \lesssim \|U\|_{C^{N_1+9}}^2  \|\Gamma U\|_{H^{N_1}}  \|\Omega\|_{H^{N_1}}.\nonumber
\end{align}
Similarly, using \eqref{GHnorm}, we obtain
\begin{align}
|\langle\partial_x\rangle^{N_1}(g_1,g_2)^T ,  \langle\partial_x\rangle^{N_1}\Omega\rangle| \lesssim \|U\|_{C^{N_1+9}}^2  (\|\Gamma U\|_{H^{N_1}}+\|U\|_{H^N})  \|\Omega\|_{H^{N_1}}.\nonumber
\end{align}

Therefore, we conclude that
\begin{align}\label{6.522}
\frac{d}{ dt}\|\Omega\|_{H^{N_1}}  \lesssim  \|U\|_{C^{N_1+9}}^2 ( \|\Gamma U\|_{H^{N_1}}  + \|U\|_{H^{N}})\lesssim \epsilon_1^3(1+t)^{p_0-1}.
\end{align}
where we have used \eqref{main a prior bound2} in the last step.
Note that from \eqref{GHnorm}, \eqref{main a prior bound2} and \eqref{4.6}--\eqref{4.7}, we have
\begin{align*}
\|\mathcal{O}[r, G_1]U\|_{H^{N_1}}+\|\mathcal{O}[u, G_2]U\|_{H^{N_1}}&\lesssim \|U\|_{C^5}\|U\|_{H^N}\lesssim \epsilon_1^2(1+t)^{p_0-1/2}\lesssim \epsilon_1^2,\\
\|\mathcal{O}[\Gamma u, H_1]U\|_{H^{N_1}}+\|\mathcal{O}[\Gamma r, H_2]U\|_{H^{N_1}}&\lesssim \|\Gamma U\|_{H^{N_1}}\|U\|_{C^{N_1+9}}\lesssim \epsilon_1^2(1+t)^{p_0-1/2}\lesssim \epsilon_1^2,\\
\|\mathcal{O}[u,A_1]\Gamma U\|_{H^{N_1}}+\|\mathcal{O}[r,A_2]\Gamma U \|_{H^{N_1}}&\lesssim \|U\|_{C^5}\|\Gamma U\|_{H^{N_1}}\lesssim \epsilon_1^2\\
\|\mathcal{O}[u,C_1]\Gamma U\|_{H^{N_1}}+\|\mathcal{O}[r,C_2]\Gamma U\|_{H^{N_1}}&\lesssim \|U\|_{C^5}\|\Gamma U\|_{H^{N_1}}\lesssim \epsilon_1^2.
\end{align*}
Hence, we deduce from \eqref{6.11} that
\begin{align*}
\|\Omega\|_{H^{N_1}}\lesssim \|\Gamma U\|_{H^{N_1}}+\epsilon_1^2,\quad\quad\quad \|\Gamma U\|_{H^{N_1}}\lesssim \|\Omega\|_{H^{N_1}}+\epsilon_1^2.
\end{align*}
Integrating \eqref{6.522} and using \eqref{state1}, we obtain
\begin{align*}
\|\Gamma U\|_{H^{N_1}}\lesssim \epsilon_0+ (1+t)^{p_0}\epsilon_1^2.
\end{align*}
Proposition \ref{weprop} thus follows.
\hfill $\Box$

\subsection{Low energy estimate for $xU$}

 In this subsection, we aim  to prove Proposition \ref{energyxuprop}.  Using the identities
 \begin{align*}
[x,\partial_x]=-I,\ \ \ [x,\partial_t]=0,\ \ \ [x,\langle\partial_x\rangle]=\frac{\partial_x}{\langle\partial_x\rangle},\ \ \ [x,\frac{\partial_x}{\langle\partial_x\rangle}]=-\frac{1}{\langle\partial_x\rangle^3},
\end{align*}
we see  $xr$ and $xu$ satisfy
\begin{align}
\begin{pmatrix}
           x r \\
           x u \\
         \end{pmatrix}_t + \begin{pmatrix}
                0  & -\langle\partial_x\rangle \\
               \langle\partial_x\rangle  & 0 \\
               \end{pmatrix}
                \begin{pmatrix}
           x r \\
           x u \\
         \end{pmatrix} = \begin{pmatrix}
           (x r)_x\langle\partial_x\rangle u+ r_x\langle\partial_x\rangle(x u) \\
            \frac{\partial_x}{\langle\partial_x\rangle}[(x r)_x r_x]
+\frac{\partial_x}{\langle\partial_x\rangle}[\langle\partial_x\rangle(x u)\langle\partial_x\rangle u] \\
         \end{pmatrix} + \begin{pmatrix}
            N'_1\\
            N'_2 \\
         \end{pmatrix},\label{x}
\end{align}
where $N'_j=N'_j(r,u)$ ($j=1,2$) are linear and quadratic terms not including $xr$ and $xu$,
\begin{align*}
N'_1(r, u)&=  \frac{\partial_x}{\langle\partial_x\rangle} u - r \langle\partial_x\rangle u +  r_x \frac{\partial_x}{\langle\partial_x\rangle} u,\\
N'_2(r, u)&= -\frac{\partial_x}{\langle\partial_x\rangle} r
- \frac{\partial_x }{\langle\partial_x\rangle}(rr_x) +  \frac{\partial_x}{\langle\partial_x\rangle} [\frac{\partial_x}{\langle\partial_x\rangle} u (\langle\partial_x\rangle u)]-\frac{1}{\langle\partial_x\rangle^3}(r_x)^2-\frac{1}{\langle\partial_x\rangle^3}(\langle\partial_x\rangle u)^2.
\end{align*}
As in \eqref{5.3}, the first nonlinear term in \eqref{x} can be decomposed into
\begin{align}
 & \begin{pmatrix}
           (x r)_x\langle\partial_x\rangle u+ r_x\langle\partial_x\rangle(x u) \\
            \frac{\partial_x}{\langle\partial_x\rangle}[(x r)_x r_x]
+\frac{\partial_x}{\langle\partial_x\rangle}[\langle\partial_x\rangle(x u)\langle\partial_x\rangle u] \\
         \end{pmatrix}
         &= \frac{1}{2}\begin{pmatrix}
           F_1(x r, x u, r, u) \\
           F_2(x r, x u, r, u)  \\
         \end{pmatrix} + \frac{1}{2}\begin{pmatrix}
           F'_1(x r, x u, r, u) \\
           F'_2(x r, x u, r, u)  \\
         \end{pmatrix}\nonumber
 \end{align}
with
 \begin{align}
(F_1,  F_2)^T &=
         \mathcal{O}[r, Q_1]xU + \mathcal{O}[u, Q_2]xU  +  \mathcal{O}[r, S_1]xU  + \mathcal{O}[u, S_2]xU ,\nonumber\\
(F'_1,  F'_2)^T &=\mathcal{O}[xr, Q_1]U + \mathcal{O}[xu, Q_2]U  +  \mathcal{O}[xr, S_1]U  + \mathcal{O}[xu, S_2]U.\nonumber
\end{align}
Note that $F_1(x r, x u, r, u)$ and $F_2(x r, x u, r, u)$ will  lead to loss of derivatives for the energy estimate, as $x U$ has higher frequencies compared to $r$ and $u$.
In order to treat this case,  we take energy normal form transformation similar to \eqref{6.9}. Let
\begin{align}\label{definition of Theta}
\begin{split}
\Theta:=xU&+\frac{1}{2}(\mathcal{O}[u,A_1] xU+ \mathcal{O}[r,A_2]xU+ \mathcal{O}[u,C_1]xU+ \mathcal{O}[r,C_2]xU)\\
&+\frac{1}{2}(\mathcal{O}[xu, H_1]U+\mathcal{O}[xr, H_2]U),
\end{split}
\end{align}
where $H_1$ and $H_2$ are the same as \eqref{5.7},  then  \eqref{x} is changed  into
\begin{align*}
 \Theta_t + D \Theta
         &= \frac{1}{2}\mathcal{O}[r,Q_1-B_1]xU + \frac{1}{2}\mathcal{O}[u,Q_2-B_2]xU + (N'_1,N'_2)^T\\
         &+ \frac{1}{2}(\mathcal{O}[f_2,A_1] + \mathcal{O}[f_1,A_2]+ \mathcal{O}[f_2,C_1]+ \mathcal{O}[f_1,C_2] )xU\\
         &+ \frac{1}{2}(\mathcal{O}[u,A_1] + \mathcal{O}[r,A_2])(\frac{1}{2}(F_1+F'_1)+N'_1,\frac{1}{2}(F_2+F'_2)+N'_2)^T \\
         &+ \frac{1}{2}( \mathcal{O}[u,C_1]+ \mathcal{O}[r,C_2])(\frac{1}{2}(F_1+F'_1)+N'_1,\frac{1}{2}(F_2+F'_2)+N'_2)^T \\
         &+\frac{1}{2}(\mathcal{O}[xu, H_1]+\mathcal{O}[xr, H_2])(f_1,f_2)^T\\
         &+\frac{1}{2}(\mathcal{O}[\frac{1}{2}(F_2+F'_2)+N'_2, H_1]+\mathcal{O}[\frac{1}{2}(F_1+F'_1)+N'_1, H_2])U.
\end{align*}
 Similar to \eqref{5.13}, we can obtain
\begin{align}\label{6.21}
\Theta_t + D \Theta
         -\frac{1}{2}\mathcal{O}[r,Q_1-B_1]\Theta -\frac{1}{2} [u,Q_2-B_2]\Theta
         = L_1 + L_2 + L_3,
\end{align}
where
\begin{align*}
L_1:=&  \frac{1}{4}\big[\mathcal{O}[u,A_1] + \mathcal{O}[r,A_2], \mathcal{O}[r,Q_1]+ \mathcal{O}[u,Q_2]\big] xU,\\
L_2:=&(N'_1,N'_2)^T+ \frac{1}{2}(\mathcal{O}[u,A_1]+ \mathcal{O}[r,A_2]+ \mathcal{O}[u,C_1]+ \mathcal{O}[r,C_2])(N'_1,N'_2)^T\\
&+\frac{1}{2}(\mathcal{O}[N'_2, H_1]+\mathcal{O}[N'_1, H_2])U,\\
L_3:=&-\frac{1}{4}(\mathcal{O}[r,Q_1-B_1]+  \mathcal{O}[u,Q_2-B_2])(\mathcal{O}[u,C_1]+ \mathcal{O}[r,C_2])xU\\
&-\frac{1}{4}(\mathcal{O}[r,Q_1-B_1]+  \mathcal{O}[u,Q_2-B_2])(\mathcal{O}[xu, H_1]+\mathcal{O}[xr, H_2])U\\
&+\frac{1}{4}(\mathcal{O}[r,B_1]+\mathcal{O}[u,B_2])(\mathcal{O}[u,A_1] + \mathcal{O}[r,A_2])xU\\
&+\frac{1}{4}(\mathcal{O}[u,A_1]+ \mathcal{O}[r,A_2])(\mathcal{O}[r,S_1]+ \mathcal{O}[u,S_2])xU
+\frac{1}{4}(\mathcal{O}[u,A_1]+ \mathcal{O}[r,A_2])(F'_1,F'_2)^T\\
&+\frac{1}{4}(\mathcal{O}[u,C_1]+ \mathcal{O}[r,C_2])(F_1+F'_1,F_2+F'_2)^T\\
&+\frac{1}{2}(\mathcal{O}[f_2,A_1] + \mathcal{O}[f_1,A_2]+ \mathcal{O}[f_2,C_1]+ \mathcal{O}[f_1,C_2] )xU\\
&+\frac{1}{2}(\mathcal{O}[xu, H_1]+\mathcal{O}[xr, H_2])(f_1,f_2)^T+\frac{1}{4}(\mathcal{O}[F_2+F'_2, H_1]+\mathcal{O}[F_1+F'_1, H_2])U.
\end{align*}
Remember that, there are  linear terms  in $N'_1$ and $N'_2$. Now, applying  similar argument as Section 3.3, we can obtain
\begin{align*}
| \langle  \langle\partial_x\rangle^{N_1} L_1,  \langle\partial_x\rangle^{N_1}\Theta\rangle| &\lesssim \|U\|_{C^5}^2\|xU\|_{H^{N_1}}\|\Theta\|_{H^{N_1}},\\
| \langle  \langle\partial_x\rangle^{N_1} L_2,  \langle\partial_x\rangle^{N_1}\Theta\rangle| &\lesssim (\|U\|_{H^{N}}+\|U\|_{C^5}\|U\|_{H^{N}}+\|U\|_{C^5}^2\|U\|_{H^{N}})\|\Theta\|_{H^{N_1}},\\
| \langle  \langle\partial_x\rangle^{N_1} L_3,  \langle\partial_x\rangle^{N_1}\Theta\rangle| &\lesssim \|U\|_{C^{N_1+9}}^2\|xU\|_{H^{N_1}}\|\Theta\|_{H^{N_1}}.
\end{align*}
Therefore, the $H^{N_1}$ energy estimate for \eqref{6.21} is
\begin{align*}
\frac{d}{ dt}\|\Theta\|_{H^{N_1}} &\lesssim   \|U\|_{C^{N_1+9}}^2   \|xU\|_{H^{N_1}} +\|U\|_{C^5}^2\|U\|_{H^{N}}+\|U\|_{H^{N}}\\
& \lesssim \epsilon_1^2(1+t)^{-1} \|xU\|_{H^{N_1}}+\epsilon_1(1+t)^{p_0}.
\end{align*}
Note that \eqref{definition of Theta} implies $\|xU\|_{H^{N_1}} \sim \|\Theta\|_{H^{N_1}}$ if $\epsilon_1$ is  small enough, so we have
\begin{align*}
\frac{d}{ dt}\|\Theta\|_{H^{N_1}} \lesssim \epsilon_1^2(1+t)^{-1} \|\Theta\|_{H^{N_1}}+\epsilon_1(1+t)^{p_0}.
\end{align*}
Using Gronwall's inequality, we obtain
$$
\|xU\|_{H^{N_1}}\sim \|\Theta\|_{H^{N_1}} \lesssim \epsilon_0(1+t)^{C\epsilon_1^2}+\epsilon_1(1+t)^{1+p_0}\lesssim\epsilon_1(1+t)^{1+p_0}
$$
provided that  $\epsilon_1$ is sufficiently small. This ends the proof of Proposition \ref{energyxuprop}.

\section{Modified scattering and decay estimate} \label{ddecay}

In this section,  we will prove Proposition \ref{decaye} below.  Recalling \eqref{defofh} and \eqref{ep002}, we have
\begin{align}\label{EP36}
h_t+i\langle\partial_x\rangle h &=\frac{1}{2i}(h+\overline{h})_x\langle\partial_x\rangle(h-\overline{h})
+\frac{\partial_x}{4i\langle\partial_x\rangle}[\langle\partial_x\rangle(h-\overline{h})]^2
-\frac{\partial_x}{4i\langle\partial_x\rangle}[(h+\overline{h})_x]^2\nonumber\\
&=\mathcal{O}[h, q^{++}]h + \mathcal{O}[h, q^{+-}]\overline{h} + \mathcal{O}[\overline{h}, q^{--}]\overline{h},
\end{align}
where the expressions for the symbols $q^{++}$, $q^{+-}$ and  $q^{--}$ are
\begin{align}\label{ps1}
\begin{split}
q^{++}(\xi,\eta)&:=\frac{1}{2}\xi\langle\eta\rangle
+\frac{\xi+\eta}{4\langle\xi+\eta\rangle}\langle\xi\rangle\langle\eta\rangle+\frac{\xi+\eta}{4\langle\xi+\eta\rangle}\xi\eta,\\
q^{+-}(\xi,\eta)&:=-\frac{1}{2}\xi\langle\eta\rangle+\frac{1}{2}\langle\xi\rangle\eta
-\frac{\xi+\eta}{2\langle\xi+\eta\rangle}\langle\xi\rangle\langle\eta\rangle+\frac{\xi+\eta}{2\langle\xi+\eta\rangle}\xi\eta,\\
q^{--}(\xi,\eta)&:=-\frac{1}{2}\xi\langle\eta\rangle
+\frac{\xi+\eta}{4\langle\xi+\eta\rangle}\langle\xi\rangle\langle\eta\rangle+\frac{\xi+\eta}{4\langle\xi+\eta\rangle}\xi\eta.
\end{split}
\end{align}
We first apply Shatah's normal form transformation to eliminate the quadratic terms in the equation \eqref{EP36}. Let
\begin{align}\label{ps-a1}
g:=h+ \mathcal{O}[h, b^{++}]h + \mathcal{O}[h, b^{+-}]\overline{h}  + \mathcal{O}[\overline{h}, b^{--}]\overline{h}
=h+\sum_{\iota_1\iota_2\in \Lambda}\mathcal{O}[h^{\iota_1}, b^{\iota_1\iota_2}]h^{\iota_2} ,
\end{align}
where $\Lambda:=\{++,+-,--\}$, $h^+:=h$, $h^-:=\overline{h}$ and
\begin{align}\label{ps2}
b^{\iota_1\iota_2}(\xi,\eta):=\frac{iq^{\iota_1\iota_2}(\xi,\eta)}
{\langle\xi+\eta\rangle-\iota_1\langle\xi\rangle-\iota_2\langle\eta\rangle},\ \ \ \iota_1\iota_2\in\Lambda.
\end{align}
We remark that for any $\xi,\eta\in \mathbb{\mathbb{R}}$,
$$|\langle\xi+\eta\rangle\pm\langle\xi\rangle\pm\langle\eta\rangle|\geq (\langle\xi+\eta\rangle+\langle\xi\rangle+\langle\eta\rangle)^{-1}>0.$$
By \eqref{ps-a1} and \eqref{ps2},  equation \eqref{EP36}  is changed  into
\begin{align}\label{ps3}
g_t+i\langle\partial_x\rangle g=\mathcal{N}(h),
\end{align}
where $\mathcal{N}(h)$ denotes cubic nonlinear term,
\begin{align}\label{ps4}
\begin{split}
\mathcal{N}(h):=&\sum_{\iota_1\iota_2\in \Lambda}\mathcal{O}[\mathcal{O}[h^{\iota_1},q^{{\iota_1\iota_2}}]h^{\iota_2},b^{++}]h
+\sum_{\iota_1\iota_2\in\Lambda}\mathcal{O}[h,b^{++}]\mathcal{O}[h^{\iota_1},q^{{\iota_1\iota_2}}]h^{\iota_2}\\
&+\sum_{\iota_1\iota_2\in\Lambda}\mathcal{O}[\mathcal{O}[h^{\iota_1},q^{{\iota_1\iota_2}}]h^{\iota_2},b^{+-}]\overline{h}
+\sum_{\iota_1\iota_2\in\Lambda}\mathcal{O}[h,b^{+-}]\overline{\mathcal{O}[h^{\iota_1},q^{{\iota_1\iota_2}}]h^{\iota_2}}\\
&+\sum_{\iota_1\iota_2\in \Lambda}\mathcal{O}[\overline{\mathcal{O}[h^{\iota_1},q^{{\iota_1\iota_2}}]h^{\iota_2}},b^{--}]\overline{h}
+\sum_{\iota_1\iota_2\in\Lambda}\mathcal{O}[\overline{h},b^{--}]\overline{\mathcal{O}[h^{\iota_1},q^{{\iota_1\iota_2}}]h^{\iota_2}}.
\end{split}
\end{align}
Let $w$ be the linear profile of $g$, that is
\begin{align}\label{ps5}
w(t):=e^{it\langle\partial_x\rangle}g(t),
\end{align}
then from \eqref{ps3}, $w$ satisfies
\begin{align}\label{ps6}
w_t=e^{it\langle\partial_x\rangle}(\partial_t+i\langle\partial_x\rangle)g=e^{it\langle\partial_x\rangle}\mathcal{N}(h).
\end{align}

Now we state the main result of this section.

\begin{prop}\label{decaye} Let $h\in C([0,T];H^N)$ be the solution of \eqref{EP36}, and $g,w$ be given by \eqref{ps-a1}, \eqref{ps5}, respectively. Assume that
\begin{align}\label{intial bound}
&\|h(0)\|_{H^{N}}+\|xh(0)\|_{H^{N_1+1}}+\|\langle\xi\rangle^{N_1+10}\widehat{h(0)}\|_{L^\infty}\leq \epsilon_0,
\end{align}
and
\begin{align}
&\sup_{t\in[0,T]}[(1+t)^{-p_0}\|h(t)\|_{H^{N}}+(1+t)^{-p_0}\|\Gamma h(t)\|_{H^{N_1}}
+\|\langle\xi\rangle^{N_1+10}\widehat{w}(t)\|_{L^\infty}\nonumber \\
&\qquad\qquad\qquad\qquad\qquad\qquad\qquad\qquad+(1+t)^{1/2}\|h(t)\|_{W^{N_1+10,\infty}}]\leq \epsilon_1,\label{a-priori wwbound}
\end{align}
where $N=300$, $N_1=15$, $0<p_0<10^{-3}$ and $0<\epsilon_0\ll\epsilon_1\ll1$. Then we have
\begin{align}
&\sup_{t\in[0,T]}[(1+t)^{-p_0}\|xw(t)\|_{H^{N_1-4}}]\lesssim \epsilon_0+\epsilon_1^2,\label{xwbound}\\
&\sup_{t\in[0,T]}\|\langle\xi\rangle^{N_1+10}\widehat{w}(t)\|_{L^\infty}\lesssim \epsilon_0+\epsilon_1^3,\label{desired wwbound3}\\
&\sup_{t\in[0,T]}[(1+t)^{1/2}\|h(t)\|_{W^{N_1+10,\infty}}]\lesssim \epsilon_0+\epsilon_1^2.\label{desired wwbound1}
\end{align}
\end{prop}

To prove Proposition \ref{decaye}, we need to construct a new linear dispersive estimate for the solution of  \eqref{EP36}.

\begin{lem}\label{linear dis est} For all $t\geq 0$, there holds that
\begin{align}\label{linear estimate}
\|e^{\pm it\langle\partial_x\rangle}f\|_{L^\infty}\lesssim (1+t)^{-1/2}\|\widehat{f} \|_{L^\infty}
+(1+t)^{-5/8}(\|f\|_{H^2}+\|xf\|_{H^1}).
\end{align}
\end{lem}

The proof for this estimate is given in Lemma A.1 of the appendix. Let $f=w$ (or $\overline{w}$) in \eqref{linear estimate}, Lemma \ref{linear dis est} shows  that the $L^\infty$ norm of the solution $g$ is controlled by the $L^\infty$ norm of $\widehat{w}$ and the Sobolev norms of $w$ and $xw$. The estimates for these norms are presented in the following subsections.

\subsection{Proof of \eqref{xwbound}}
We need the following isotropic multiplier estimate for  $\mathcal{O} [h^{\iota_1}, q^{\iota_1\iota_2}]h^{\iota_2}$   and $\mathcal{N}(h)$.

\begin{lem}\label{multiplierlemma2d}
Let $m(\xi,\eta)$ be a Fourier multiplier satisfying
\begin{align}\label{multiplier bound3}
\|m\|_{L^2(\mathbb{R}^2)}+\|\partial_\xi^2 m\|_{L^2(\mathbb{R}^2)}+\|\partial_\eta^2 m\|_{L^2(\mathbb{R}^2)}\lesssim 1,
\end{align}
then for any $p_0,\ p_1,\ p_2\in [1,+\infty]$ with $p_0^{-1}=p_1^{-1}+p_2^{-1}$, we have
\begin{align*}
\|\mathcal{O}[f_1, m]f_2\|_{L^{p_0}(\mathbb{R})}
\lesssim
\|f_1\|_{L^{p_1}(\mathbb{R})}\|f_2\|_{L^{p_2}(\mathbb{R})}.
\end{align*}
\end{lem}

For the proof of this multiplier lemma,  see Lemma B.2  in the appendix.

\begin{lem}\label{wwlemma1}
Under the same assumptions as Proposition \ref{decaye}, there hold
\begin{align}
&\|\mathcal{O} [h^{\iota_1}, b^{\iota_1\iota_2}]h^{\iota_2}\|_{H^{N-5}}\lesssim \epsilon_1^2(1+t)^{p_0-1/2},\label{ps10}\\
&\|\mathcal{O} [h^{\iota_1}, b^{\iota_1\iota_2}]h^{\iota_2}\|_{W^{N_1+10,\infty}}\lesssim \epsilon_1(1+t)^{p_0/2-1/4}\|h\|_{W^{N_1+10,\infty}}, \label{ps12}\\
&\|\Gamma\mathcal{O} [h^{\iota_1}, b^{\iota_1\iota_2}]h^{\iota_2}\|_{H^{N_1-5}}\lesssim \epsilon_1^2(1+t)^{p_0-1/2},\label{ps11}
\end{align}
where  $\iota_1\iota_2\in\Lambda=\{++,+-,--\}$. Moreover, we have
\begin{align}\label{ps13}
\|x\mathcal{N}(h)\|_{H^{N_1-5}}\lesssim \epsilon_1^3(1+t)^{p_0}.
\end{align}
\end{lem}

{\noindent \emph{Proof.}}  It follows from \eqref{a-priori wwbound} that
\begin{align}\label{ps14-1}
\|h\|_{H^{N}}\lesssim \epsilon_1(1+t)^{p_0},\ \ \ \|\Gamma h\|_{H^{N_1}}\lesssim \epsilon_1(1+t)^{p_0},\ \ \
\|h\|_{W^{N_1+10,\infty}}\lesssim \epsilon_1(1+t)^{-1/2}.
\end{align}
By the definition \eqref{defofO},
$$
\mathscr{F}(\langle\partial_x\rangle^{N-5}\mathcal{O} [h^{\iota_1}, b^{\iota_1\iota_2}]h^{\iota_2})(\xi)
=\frac{1}{2\pi}\int_{\mathbb{R}}m^{\iota_1\iota_2}(\xi-\eta,\eta)(\langle\xi-\eta\rangle^N
+\langle\eta\rangle^N)\widehat{h^{\iota_1}}(\xi-\eta)\widehat{h^{\iota_2}}(\eta)d\eta
$$
with
$$
m^{\iota_1\iota_2}(\xi-\eta,\eta):=\frac{\langle\xi\rangle^{N-5}b^{\iota_1\iota_2}(\xi-\eta,\eta)}
{\langle\xi-\eta\rangle^N+\langle\eta\rangle^N}.
$$
Note that
\begin{align}\label{ps14}
\left|\partial_\xi^{a_1}\partial_\eta^{a_2} \left[\frac{1}{\langle\xi+\eta\rangle\pm\langle\xi\rangle\pm\langle\eta\rangle}\right]\right|\lesssim \max(\langle\xi+\eta\rangle,\langle\xi\rangle,\langle\eta\rangle),\ \ a_1,\ a_2\geq 0,
\end{align}
then we deduce from \eqref{ps1} and \eqref{ps2} that
\begin{align}\label{ps-a2}
|\partial_\xi^{a_1}\partial_\eta^{a_2}  b^{\iota_1\iota_2}(\xi,\eta)|\lesssim  (\max(\langle\xi+\eta\rangle,\langle\xi\rangle,\langle\eta\rangle))^3.
\end{align}
In view of \eqref{ps-a2}, it is easy to check that $m^{\iota_1\iota_2}(\xi,\eta)$ satisfies \eqref{multiplier bound3}, then Lemma \ref{multiplierlemma2d} shows
 \begin{align}
 \|\mathcal{O} [h^{\iota_1}, b^{\iota_1\iota_2}]h^{\iota_2}\|_{H^{N-5}}\lesssim \|h^{\iota_1}\|_{H^{N}}\|h^{\iota_2}\|_{L^\infty} + \|h^{\iota_1}\|_{L^\infty}\|h^{\iota_2}\|_{H^{N}}\lesssim \epsilon_1^2(1+t)^{p_0-1/2}, \label{4.23}
 \end{align}
where we have used \eqref{ps14-1} in the last step. Hence, the bound \eqref{ps10} follows.

For  \eqref{ps12},
using the interpolation inequality,
$$
\|h\|_{W^{N_1+15,\infty}}\lesssim  \|h\|_{W^{N_1+10,\infty}}^{1/2}\|h\|_{W^{N_1+20,\infty}}^{1/2}\lesssim\epsilon_1(1+t)^{p_0/2-1/4},
$$
then by Lemma \ref{multiplierlemma2d}, we have
\begin{align*}
 \|\mathcal{O} [h^{\iota_1}, b^{\iota_1\iota_2}]h^{\iota_2}\|_{W^{N_1+10,\infty}}
 &\lesssim \|h\|_{L^{\infty}}\|h\|_{W^{N_1+15,\infty}}\lesssim \epsilon_1(1+t)^{p_0/2-1/4}\|h\|_{L^{\infty}}.
 \end{align*}

To prove \eqref{ps11}, we first consider the case $\iota_1\iota_2=++$. A direct computation gives
\begin{align}\label{ps-a3}
\begin{split}
\mathscr{F}(\Gamma\mathcal{O} [h, b^{++}]h)(\xi)=&\mathscr{F}(\mathcal{O} [\Gamma h, b^{++}]h)(\xi)+\mathscr{F}(\mathcal{O} [h, b^{++}]\Gamma h)(\xi)\\
&+\frac{i}{2\pi}\int_\mathbb{R}\partial_\xi b^{++}(\xi-\eta,\eta)\widehat{h}_t(\xi-\eta)\widehat{h}(\eta)d\eta\\
&+\frac{i}{2\pi}\int_\mathbb{R}(\partial_\xi b^{++} (\xi-\eta,\eta)+\partial_\eta b^{++} (\xi-\eta,\eta))\widehat{h}(\xi-\eta)\widehat{h}_t(\eta)d\eta.
\end{split}
\end{align}
From the equation \eqref{EP36} and the bound \eqref{ps14-1}, it is easy to see
$$
\|h_t\|_{H^{N-1}}\lesssim  \epsilon_1(1+t)^{p_0},\ \ \ \ \|h_t\|_{L^\infty}\lesssim\epsilon_1(1+t)^{-1/2}.
$$
Then  using  \eqref{ps-a2}--\eqref{ps-a3} and Lemma \ref{multiplierlemma2d}, we obtain
\begin{align*}
\|\Gamma\mathcal{O} [h, b^{++}]h\|_{H^{N_1-5}}&\lesssim \|\Gamma h\|_{H^{N_1}}\|h\|_{W^{N_1,\infty}}+\|h_t\|_{H^{N_1}}\|h\|_{L^{\infty}}+\|h_t\|_{L^{\infty}}\|h\|_{H^{N_1}}\\
& \lesssim\epsilon_1^2(1+t)^{p_0-1/2},
\end{align*}
which proves \eqref{ps11} for $\iota_1\iota_2=++$. The proof for  $\iota_1\iota_2=+-,--$ is the same as above.

In order to prove  \eqref{ps13}, it suffices to show  that each term in \eqref{ps4} satisfies the bound \eqref{ps13}. Here, we only consider the term $\mathcal{O}[\mathcal{O}[h,q^{++}]h,b^{++}]h$ in detail. Note that
\begin{align*}
\mathscr{F}(\mathcal{O}[\mathcal{O}[h,q^{++}]h,b^{++}]h)(\xi)&=\frac{1}{2\pi}\int_{\mathbb{R}} b^{++}(\eta,\xi-\eta)
\mathscr{F}(\mathcal{O}[h,q^{++}]h)(\eta)(\mathscr{F}h)(\xi-\eta)d\eta.
\end{align*}
Applying $\partial_\xi$ to this identity yields
$$
x\mathcal{O}[\mathcal{O}[h,q^{++}]h,b^{++}]h=A_1+A_2,
$$
where
\begin{align*}
\widehat{A}_1(\xi)&:=\frac{i}{2\pi}\int_{\mathbb{R}}\partial_\xi b^{++}(\eta,\xi-\eta)\mathscr{F}(\mathcal{O}[h,q^{++}]h)(\eta)(\mathscr{F}h)(\xi-\eta)d\eta,\\
\widehat{A}_2(\xi)&:=\frac{i}{2\pi}\int_{\mathbb{R}} b^{++}(\eta,\xi-\eta)\mathscr{F}(\mathcal{O}[h,q^{++}]h)(\eta)\partial_\xi(\mathscr{F}h)(\xi-\eta)d\eta.
\end{align*}
Then using Lemma \ref{multiplierlemma2d}, \eqref{ps14-1}, \eqref{ps-a2} and \eqref{ps1}, we have
\begin{align*}
\|A_1\|_{H^{N_1-5}}&\lesssim \|\mathcal{O}[h,q^{++}]h\|_{H^{N_1}}\|h\|_{W^{N_1,\infty}}\lesssim \|h\|_{W^{N_1+10,\infty}}^2\|h\|_{H^N}
\lesssim \epsilon_1^3(1+t)^{p_0-1}.
\end{align*}
For the term $A_2$,  Proposition \ref{energyxuprop} and \eqref{defofh} yield
$$\|xh\|_{H^{N_1}}\sim \|xU\|_{H^{N_1}}\lesssim \epsilon_1(1+t)^{1+p_0},$$
hence, we obtain
\begin{align*}
\|A_2\|_{H^{N_1-5}}&\lesssim \|\mathcal{O}[h,q^{++}]h\|_{W^{N_1,\infty}}\|xh\|_{H^{N_1}}\lesssim
\|h\|_{W^{N_1+10,\infty}}^2\|xh\|_{H^{N_1}}
\lesssim \epsilon_1^3(1+t)^{p_0}.
\end{align*}
Therefore, we conclude that
\begin{align*}
\|x\mathcal{O}[\mathcal{O}[h,q^{++}]h,b^{++}]h\|_{H^{N_1-5}}\lesssim \epsilon_1^3(1+t)^{p_0},
\end{align*}
and the desired bound \eqref{ps13} thus follows.
\hfill $\Box$\newline

{\noindent \emph{Proof of \eqref{xwbound}.}}  To estimate $xw$, an important tool  is introducing  the vector filed
\begin{align}\label{ps7}
\widetilde{\Gamma}:=t\partial_x-i\langle\partial_x\rangle x,
\end{align}
which satisfies
\begin{align}\label{ps9}
\langle\xi\rangle\partial_\xi(e^{it\langle\xi\rangle}\mathscr{F} g)=e^{it\langle\xi\rangle}\mathscr{F}(\widetilde{\Gamma} g).
\end{align}
Thus, $\|xw\|_{H^{s+1}} \sim \| \widetilde{\Gamma} g\|_{H^s}$. Moreover, the relationship between $\widetilde{\Gamma}$ and the homogeneous vector field operator $\Gamma=x\partial_t+t\partial_x$ is
\begin{align*}
\widetilde{\Gamma} g=\Gamma g-x(\partial_t+i\langle\partial_x\rangle)g+\frac{i\partial_x}{\langle\partial_x\rangle}g
=\Gamma g-x\mathcal{N}(h)+\frac{i\partial_x}{\langle\partial_x\rangle}g.
\end{align*}
Using the bounds  \eqref{main energy bound1}, \eqref{main wenergy bound1}, \eqref{ps10} and  \eqref{ps11}, we deduce
\begin{align}
\|g\|_{H^{N-5}}&\lesssim \|h\|_{H^{N-5}}+\sum_{\iota_1\iota_2\in \Lambda}\|\mathcal{O} [h^{\iota_1}, b^{\iota_1\iota_2}]h^{\iota_2}\|_{H^{N-5}}\lesssim (\epsilon_0+\epsilon_1^2)(1+t)^{p_0},\label{wn-5}\\
\|\Gamma g\|_{H^{N_1-5}}&\lesssim \|\Gamma h\|_{H^{N_1-5}}+\sum_{\iota_1\iota_2\in \Lambda}\|\Gamma \mathcal{O} [h^{\iota_1}, b^{\iota_1\iota_2}]h^{\iota_2}\|_{H^{N_1-5}}\lesssim (\epsilon_0+\epsilon_1^2)(1+t)^{p_0}.\nonumber
\end{align}
Hence, combining \eqref{ps13} and the above estimates, we obtain
$$
\|\widetilde{\Gamma} g\|_{H^{N_1-5}}\lesssim \|\Gamma g\|_{H^{N_1-5}}+\|x\mathcal{N}(h)\|_{H^{N_1-5}}+\|g\|_{H^{N_1-5}}\lesssim (\epsilon_0+\epsilon_1^2)(1+t)^{p_0}.
$$
Thanks to the identity \eqref{ps9}, there holds
\begin{align*}
\|xw\|_{H^{N_1-4}}=\|\widetilde{\Gamma} g\|_{H^{N_1-5}}\lesssim (\epsilon_0+\epsilon_1^2)(1+t)^{p_0}.
\end{align*}
The proof of \eqref{xwbound} is completed.\hfill $\Box$

\subsection{Proof of \eqref{desired wwbound3}}

Now we consider the  $L^\infty$ bound for $\widehat{w}$ in low order norm and  present the proof of \eqref{desired wwbound3}. Indeed, we will be devoted in proving  a more  stronger result in this subsection.

\begin{prop}\label{scproposition}
Under the same assumption as Proposition \ref{decaye}, there exists $\delta>0$ such that
\begin{align}\label{sc1}
\sup_{0\leq t_1\leq t_2\leq T}(1+t_1)^{\delta}\|\langle\xi\rangle^{N_1+10} e^{i\vartheta(t_1,\xi)}\widehat{w}(t_1,\xi)-\langle\xi\rangle^{N_1+10} e^{i\vartheta(t_2,\xi)}\widehat{w}(t_2,\xi)\|_{L^\infty}\lesssim \epsilon_1^3,
\end{align}
where $w$ is defined by \eqref{ps5}, and  $\vartheta$ is a real-valued function given by  \eqref{definition of theta}.
\end{prop}

Once  Theorem \ref{mainthm}  is proved, the above proposition implies that the function
$$\langle\xi\rangle^{N_1+10} e^{i\vartheta(t,\xi)}\widehat{w}(t,\xi)$$
 forms a Cauchy family as $t\rightarrow \infty$ in $L^\infty$, so there exists a unique $w_\infty(\xi)\in L^\infty$ such that
\begin{align*}
\sup_{t\geq 0}[(1+t)^\delta\|\langle\xi\rangle^{N_1+10} e^{i\vartheta(t,\xi)}\widehat{w}(t,\xi)-w_\infty(\xi)\|_{L^\infty}]\lesssim \epsilon_0.
\end{align*}
This result says the solution of the equation \eqref{ps3} tends to a nonlinear asymptotic state as $t\rightarrow \infty$, thus  such equation possesses a modified scattering behavior with corrected phase $\vartheta(t,\xi)$. Assuming Proposition \ref{scproposition} holds, we now show the proof of \eqref{desired wwbound3}.
\newline

{\noindent\emph{ Proof of \eqref{desired wwbound3}.}}  By setting $t_1=0$ and  $t_2=t$ in the estimate \eqref{sc1}, we see
\begin{align*}
\sup_{t\in [0,T]}\|\langle\xi\rangle^{N_1+10} \widehat{w}(t,\xi)\|_{L^\infty}
\lesssim \epsilon_1^3+\|\langle\xi\rangle^{N_1+10} \widehat{w(0)}\|_{L^\infty}=\epsilon_1^3
+\|\langle\xi\rangle^{N_1+10} \widehat{g(0)}\|_{L^\infty},
\end{align*}
then the bound \eqref{desired wwbound3} follows immediately, provided that we can show
\begin{align}\label{sc2}
\|\langle\xi\rangle^{N_1+10} \widehat{g(0)}\|_{L^\infty}\lesssim \epsilon_0.
\end{align}
Indeed, note that from \eqref{ps-a1},
\begin{align*}
g(0)=h_0+\mathcal{O} [h_0, b^{++}]h_0+\mathcal{O} [h_0, b^{+-}]\overline{h_0}+\mathcal{O} [\overline{h_0}, b^{--}]\overline{h_0},\ \ \ h_0:=h(0).
\end{align*}
Using the initial bound \eqref{intial bound}, we deduce  that,  for all $\iota_1\iota_2\in\Lambda$,
\begin{align*}
\|\langle\xi\rangle^{N_1+10} \mathscr{F}(\mathcal{O} [h_0^{\iota_1}, b^{\iota_1\iota_2}]h_0^{\iota_2})(\xi)\|_{L^\infty}\lesssim \|\mathcal{O} [h_0^{\iota_1}, b^{\iota_1\iota_2}]h_0^{\iota_2}\|_{W^{N_1+10,1}}\lesssim
\|h_0\|_{H^{N}}^2\lesssim \epsilon_0^2.
\end{align*}
Therefore, the bound \eqref{sc2} follows from the above estimate and \eqref{intial bound}.
\hfill $\Box$

From now on, we concentrate  on the proof of Proposition \ref{scproposition}.
Rewrite the nonlinear term of  the equation \eqref{ps3} as
\begin{align}\label{sc-a1}
\mathcal{N}(h)=\mathcal{N}(g)+\mathcal{N}_R,\ \ \ \ \ \ \ \mathcal{N}_R:=\mathcal{N}(h)-\mathcal{N}(g),
\end{align}
then the profile  $w$ satisfies
\begin{align}\label{sc3}
w_t=e^{it\langle\partial_x\rangle}\mathcal{N}(g)+e^{it\langle\partial_x\rangle}\mathcal{N}_R,
\end{align}
where $\mathcal{N}(g)$ denotes cubic term and $\mathcal{N}_R$ is quartic term. From the definition \eqref{ps4}, the first nonlinear term in the RHS of \eqref{sc3} can be expanded as
\begin{align*}
e^{it\langle\xi\rangle}\widehat{\mathcal{N}(g)}(\xi):=i(2\pi)^{-2}[I^{++-}(t,\xi)+I^{+--}(t,\xi)+ I^{+++}(t,\xi)+I^{---}(t,\xi)],
\end{align*}
where
\begin{align}\label{sc5}
I^{\iota_1\iota_2\iota_3}(t,\xi)&:=\int_{\mathbb{R}^2}c^{\iota_1\iota_2\iota_3}(\xi,\eta,\sigma)
e^{it\Psi^{\iota_1\iota_2\iota_3}(\xi,\eta,\sigma)}\widehat{w^{\iota_1}}(t,\xi-\eta)
\widehat{w^{\iota_2}}(t,\eta-\sigma)\widehat{w^{\iota_3}}(t,\sigma)d\eta d\sigma
\end{align}
with $\iota_1\iota_2\iota_3\in \mathscr{T}:=\{++-,+--,+++,---\}$ and $w^+:=w$, $w^-:=\overline{w}$. If there is no confusion occurs, we also simply write \eqref{sc5} as
\begin{align*}
I^{\iota_1\iota_2\iota_3}
&=\int_{\mathbb{R}^2}c^{\iota_1\iota_2\iota_3}e^{it\Psi^{\iota_1\iota_2\iota_3}}\widehat{w^{\iota_1}}(\xi-\eta)
\widehat{w^{\iota_2}}(\eta-\sigma)\widehat{w^{\iota_3}}(\sigma)d\eta d\sigma.
\end{align*}
The phase $\Psi^{\iota_1\iota_2\iota_3}$ is defined by
\begin{align}\label{phase function}
\Psi^{\iota_1\iota_2\iota_3}(\xi,\eta,\sigma):=\langle\xi\rangle-\iota_1\langle\xi-\eta\rangle-\iota_2\langle\eta-\sigma\rangle
-\iota_3\langle\sigma\rangle,
\end{align}
and the symbols $c^{\iota_1\iota_2\iota_3}$ are
\begin{align*}
ic^{++-}(\xi,\eta,\sigma):=&b^{++}(\eta,\xi-\eta)q^{+-}(\eta-\sigma,\sigma)+b^{++}(\xi-\eta,\eta)q^{+-}(\eta-\sigma,\sigma)\\
&+b^{+-}(\xi-\sigma,\sigma)q^{++}(\xi-\eta,\eta-\sigma)+b^{+-}(\xi-\eta,\eta)q^{+-}(-\sigma,\sigma-\eta)\\
&+b^{--}(\xi-\sigma,\sigma)q^{--}(\eta-\xi,\sigma-\eta)
+b^{--}(\sigma,\xi-\sigma)q^{--}(\eta-\xi,\sigma-\eta),\\
ic^{+--}(\xi,\eta,\sigma):=&b^{++}(\eta,\xi-\eta)q^{--}(\eta-\sigma,\sigma)+b^{++}(\xi-\eta,\eta)q^{--}(\eta-\sigma,\sigma)\\
&+b^{+-}(\xi-\sigma,\sigma)q^{+-}(\xi-\eta,\eta-\sigma)+b^{+-}(\xi-\eta,\eta)q^{++}(\sigma-\eta,-\sigma)\\
&+b^{--}(\xi-\sigma,\sigma)q^{+-}(\sigma-\eta,\eta-\xi)
+b^{--}(\sigma,\xi-\sigma)q^{+-}(\sigma-\eta,\eta-\xi),\\
ic^{+++}(\xi,\eta,\sigma):=&b^{++}(\eta,\xi-\eta)q^{++}(\eta-\sigma,\sigma)+b^{++}(\xi-\eta,\eta)q^{++}(\eta-\sigma,\sigma)\\
&+b^{+-}(\xi-\eta,\eta)q^{--}(\sigma-\eta,-\sigma),\\
ic^{---}(\xi,\eta,\sigma):=&b^{+-}(\xi-\sigma,\sigma)q^{--}(\xi-\eta,\eta-\sigma)+b^{--}(\xi-\sigma,\sigma)q^{++}(\eta-\xi,\sigma-\eta)\\
&+b^{--}(\sigma,\xi-\sigma)q^{++}(\eta-\xi,\sigma-\eta),
\end{align*}
where $q^{\iota_1\iota_2}$ and $b^{\iota_1\iota_2}$  are given by \eqref{ps1} and \eqref{ps2}. Therefore, we conclude that
\begin{align}\label{sc6}
\widehat{w}_t(t,\xi)=\sum_{\iota_1\iota_2\iota_3\in \mathscr{T}}i(2\pi)^{-2} I^{\iota_1\iota_2\iota_3}(t,\xi)
+e^{it\langle\xi\rangle}\widehat{\mathcal{N}_R}(t,\xi).
\end{align}
 For the phase $\Psi^{\iota_1\iota_2\iota_3}$, we can compute the space-time resonance set (\cite{GMS})
 $$
 \{(\xi,\eta,\sigma);\ \Psi^{\iota_1\iota_2\iota_3}(\xi,\eta,\sigma)=\Psi_\eta^{\iota_1\iota_2\iota_3}(\xi,\eta,\sigma)
 =\Psi_\sigma^{\iota_1\iota_2\iota_3}(\xi,\eta,\sigma)=0\}.
 $$
Indeed, it is easy to check that the only space-time resonance is in the case $\iota_1\iota_2\iota_3=++-$, and the resonant set is $(\xi,\eta,\sigma)=(\xi,0,-\xi)$. A direct computation gives
\begin{align}\label{resonant1}
\begin{split}
c^*(\xi):=&c^{++-}(\xi,0,-\xi)\\
=&\xi^2\Big [2\langle\xi\rangle-\frac{(2\langle\xi\rangle+\langle2\xi\rangle)
\cdot(\langle\xi\rangle\langle2\xi\rangle+\xi^2+\langle\xi\rangle^2)^2}
{6\langle\xi\rangle\langle2\xi\rangle}+\frac{(\langle\xi\rangle\langle2\xi\rangle-\langle\xi\rangle^2-\xi^2)^2}{2(2\langle\xi\rangle+\langle2\xi\rangle)
\langle\xi\rangle\langle2\xi\rangle}\Big],
\end{split}
\end{align}
thus,
\begin{align}\label{resonant2}
c^*(0)=c_\xi^*(0)=0,\ \ \ \ |c^*(\xi)|\lesssim \xi^2\langle\xi\rangle^3, \ \ \ \ | c^*_\xi(\xi)|\lesssim |\xi|\langle\xi\rangle^3.
\end{align}

Define
\begin{align}\label{definition of theta}
\vartheta(t,\xi):=-\frac{c^*(\xi)\langle\xi\rangle^3}{2\pi}\int_0^t\frac{|\widehat{w}(s,\xi)|^2}{s+1}ds,
\end{align}
then it follows from  \eqref{sc6} and \eqref{definition of theta} that
\begin{align}\label{corretcedprofile}
\begin{split}
\partial_t[e^{i\vartheta(t,\xi)}\widehat{w}(t,\xi)]=&e^{i\vartheta(t,\xi)}i\vartheta_t(t,\xi)\widehat{w}(t,\xi)
+e^{i\vartheta(t,\xi)}\partial_t\widehat{w}(t,\xi)\\
=&\frac{i}{4\pi^2} e^{i\vartheta(t,\xi)}\Big[I^{++-}(t,\xi)
-2\pi\frac{c^*(\xi)\langle\xi\rangle^3|\widehat{w}(t,\xi)|^2\widehat{w}(t,\xi)}{1+t}\Big]\\
&+\frac{i}{4\pi^2}e^{i\vartheta(t,\xi)}\big[I^{+--}(t,\xi)+I^{+++}(t,\xi)+I^{---}(t,\xi)\big]\\
&+e^{i\vartheta(t,\xi)}e^{it\langle\xi\rangle}\widehat{\mathcal{N}_R}(t,\xi).
\end{split}
\end{align}
Now we make frequency decomposition. Let
\begin{align}\label{definition of I decomposition}
I_{k_1k_2k_3}^{\iota_1\iota_2\iota_3}(s,\xi):=&\int_{\mathbb{R}^2}c^{\iota_1\iota_2\iota_3}_{k_1 k_2 k_3}(\xi,\eta,\sigma)
e^{is\Psi^{\iota_1\iota_2\iota_3}(\xi,\eta,\sigma)}\nonumber\\
&\qquad \cdot\widehat{P_{k_1}w^{\iota_1}}(s,\xi-\eta)\widehat{P_{k_2}w^{\iota_2}}(s,\eta-\sigma)
\widehat{P_{k_3}w^{\iota_3}}(s,\sigma)d\eta d\sigma,
\end{align}
where
\begin{align*}
c^{\iota_1\iota_2\iota_3}_{k_1 k_2 k_3}(\xi,\eta,\sigma):=c^{\iota_1\iota_2\iota_3}(\xi,\eta,\sigma)
\varphi_{k_1}(\xi-\eta)\varphi_{k_2}(\eta-\sigma)\varphi_{k_3}(\sigma).
\end{align*}
For our proof, it is sufficient to use the following bound for this symbol
\begin{align}\label{sc8}
|\partial_\xi^{a_1}\partial_\eta^{a_2}\partial_\sigma^{a_3}c^{\iota_1\iota_2\iota_3}_{k_1 k_2 k_3}(\xi,\eta,\sigma)|\lesssim 2^{5\max(k_1,k_2,k_3)_+},\ \ a_1,a_2,a_3\geq 0,
\end{align}
where $a_+:=\max\{a,0\}$. \eqref{sc8} can be obtained from the definitions of $c^{\iota_1\iota_2\iota_3}$ and a direct computation.
The detailed expressions of $c^{\iota_1\iota_2\iota_3}$  won't play an important role  in our succeeding arguments.

In virtue of \eqref{corretcedprofile}--\eqref{definition of I decomposition}, in order to prove \eqref{sc1}, it suffices to prove there exists $\delta>0$ such that
\begin{align}\label{sc5.1}
&\sum_{k_1,k_2,k_3\in \mathbb{Z}}\Big|\int_{t_1}^{t_2}e^{i\vartheta(s,\xi)}\Big[I_{k_1k_2k_3}^{++-}(s,\xi)-2\pi
\frac{c^*(\xi)\langle\xi\rangle^3\widehat{P_{k_1}w}(s,\xi)\widehat{P_{k_2}w}(s,\xi)
\widehat{P_{k_3}\overline{w}}(s,-\xi)}{s+1}\Big]ds\Big|\nonumber\\
&\qquad\qquad\qquad\qquad\qquad\qquad\qquad\qquad\qquad\qquad\qquad
\lesssim \epsilon_1^32^{-\delta m}2^{-(N_1+10)k_+},
\end{align}
and for $\iota_1\iota_2\iota_3\in\{+--,+++,---\}$,
\begin{align}\label{sc5.2}
\sum_{k_1,k_2,k_3\in \mathbb{Z}}\Big|\int_{t_1}^{t_2}e^{i\vartheta(s,\xi)}I_{k_1k_2k_3}^{\iota_1\iota_2\iota_3}(s,\xi)ds\Big|\lesssim \epsilon_1^32^{-\delta m}2^{-(N_1+10)k_+},
\end{align}
where $|\xi|\sim 2^k$, $k\in\mathbb{Z}$, and $t_1,t_2\in [2^{m}-2,2^{m+1}]\cap[0,T]$, $m=1,2,3,\cdots$. Moreover, we shall also prove
\begin{align*}
\Big|\langle\xi\rangle^{N_1+10}\int_{t_1}^{t_2}
e^{i\vartheta(s,\xi)}e^{is\langle\xi\rangle}\widehat{\mathcal{N}_R}(s,\xi)ds\Big|
\lesssim \epsilon_1^4(1+t_1)^{-\delta}.
\end{align*}

To prove these bounds,  we need some basic estimates for the localized function $P_kw$, which are given in the following lemma.

\begin{lem}  With the same assumption as Proposition \ref{decaye}, we have
\begin{align}
\|\widehat{P_kw}\|_{L^\infty}&\lesssim \epsilon_1 2^{-(N_1+10)k_+},\label{Linfinity}\\
\|\partial_\xi\widehat{P_kw}\|_{L^2}&\lesssim \epsilon_1 2^{p_0 m}2^{-(N_1-4)k_+},\label{L^2 of paritial derivative estimate}\\
\|P_kw\|_{L^2}&\lesssim \epsilon_1  2^{p_0 m}2^{-(N-5)k_+},\label{L^2estimate}\\
\|e^{\pm is\langle\partial_x\rangle}P_kw\|_{L^\infty}&\lesssim \epsilon_1  2^{-m/2},\label{Linfinity of profile}\\
\|e^{\pm is\langle\partial_x\rangle}P_kw\|_{L^\infty}&\lesssim  \epsilon_1 2^k2^{-(N_1+10)k_+}, \label{Linfinity of low}\\
\|e^{\pm is\langle\partial_x\rangle}P_kw\|_{L^2}&\lesssim \epsilon_1  2^{k/2}2^{-(N_1+10)k_+},\label{L^2 of low}\\
\|\partial_s\widehat{P_kw}\|_{L^2}&\lesssim \epsilon_1 ^3 2^{p_0 m}2^{-m}2^{-(N-7)k_+},\label{L^2 of partial time}
\end{align}
where  $s\in[2^m-2,2^{m+1}]$, $m\in \mathbb{N}$ and $k_+=\max\{k,0\}$.
\end{lem}

{\noindent \emph{Proof.}}
The bounds \eqref{Linfinity}, \eqref{L^2 of paritial derivative estimate} follow from \eqref{a-priori wwbound}, \eqref{xwbound}, respectively. Using \eqref{ps-a1}, \eqref{ps10} and  \eqref{a-priori wwbound}, we can obtain
\begin{align}\label{sc-a2}
\|w\|_{H^{N-5}}\lesssim \epsilon_1(1+t)^{p_0},
\end{align}
which gives \eqref{L^2estimate}. The bound \eqref{Linfinity of profile} follows from \eqref{linear estimate}, \eqref{a-priori wwbound}, \eqref{xwbound} and \eqref{sc-a2}. Note that
\begin{align*}
|e^{\pm is\langle\partial_x\rangle}P_kw|
=\frac{1}{2\pi}\Big|\int_{\mathbb{R}}e^{ix\xi}e^{\pm is\langle\xi\rangle}\widehat{P_kw}(s,\xi)d\xi\Big|
\lesssim \|\widehat{P_kw}\|_{L^\infty}2^k,
\end{align*}
then \eqref{Linfinity of low} follows from \eqref{Linfinity}. The estimate \eqref{L^2 of low} is proved by the  Plancherel's identity, Cauchy-Schwarz inequality  and \eqref{Linfinity}. For \eqref{L^2 of partial time},  we can obtain  from  \eqref{ps4} and \eqref{ps6} that
\begin{align*}
\|\partial_sw\|_{H^{N-7}}=\|\mathcal{N}(h)\|_{H^{N-7}}\lesssim \|h\|_{H^N}\|h\|_{W^{N_1+10,\infty}}^2\lesssim \epsilon_1^3(1+s)^{-1+p_0},
\end{align*}
so the desired bound \eqref{L^2 of partial time} follows easily.
\hfill $\Box$

We first show  the bounds \eqref{sc5.1} and \eqref{sc5.2}  in two simpler cases.

\begin{lem}\label{simpler lemma1}
 The bounds \eqref{sc5.1} and \eqref{sc5.2} hold if we take the sum over those $(k_1,k_2,k_3)$ satisfying
\begin{align}\label{sc5.14}
\min(k_1,k_2,k_3)\leq -4m\ \ \mathrm{or}\ \ \max(k_1,k_2,k_3)\geq m/200-100.
\end{align}
\end{lem}

{\noindent \emph{Proof.}} Using \eqref{sc8}, \eqref{L^2estimate} and Cauchy-Schwarz inequality, we see that
\begin{align*}
|I_{k_1k_2k_3}^{\iota_1\iota_2\iota_3}(s,\xi)|\lesssim&
\epsilon_1^3 2^{3p_0 m}2^{5\max(k_1,k_2,k_3)_+}\cdot2^{\min(k_1,k_2,k_3)/2}\cdot \\
&2^{-(N-5)\mathrm{min}(k_1,k_2,k_3)_+}\cdot2^{-(N-5)\mathrm{med}(k_1,k_2,k_3)_+}\cdot2^{-(N-5)\max(k_1,k_2,k_3)_+}
\end{align*}
for any $\iota_1\iota_2\iota_3\in\mathscr{T}$. Using \eqref{resonant2} and the $L^\infty$ bound \eqref{Linfinity}, there holds
\begin{align*}
&\Big|\frac{c^*(\xi)\langle\xi\rangle^3\widehat{P_{k_1}w}(s,\xi)\widehat{P_{k_2}w}(s,\xi)
\widehat{P_{k_3}\overline{w}}(s,-\xi)}{s+1}\Big|\\
&\qquad\qquad\lesssim 2^{-m}\epsilon_1^32^{2k}2^{6k_+}2^{-3(N_1+10)k_+}
\textbf{1}_{[0,20]}\max(|k_1-k|,|k_2-k|,|k_3-k|).
\end{align*}
In virtue of \eqref{sc5.14}, Lemma \ref{simpler lemma1} thus follows.
\hfill $\Box$

If one takes the sum for $(k_1,k_2,k_3)$ which satisfies $k_1,k_2,k_3 \in [-4m,m/200-100]\cap \mathbb{Z}$, then there are at most $Cm^3$ terms, which are summable as the desired estimate \eqref{sc5.1} or \eqref{sc5.2} has an exponential factor $2^{-\delta m}$. So  in the following it is sufficient for us to prove that  for fixed $k_1,k_2,k_3$, there exists $\delta>0$ such that
\begin{align}\label{sc9}
&\Big|\int_{t_1}^{t_2}e^{i\vartheta(s,\xi)}\big[I_{k_1k_2k_3}^{++-}(s,\xi)-2\pi
\frac{c^*(\xi)\langle\xi\rangle^3\widehat{P_{k_1}w}(s,\xi)\widehat{P_{k_2}w}(s,\xi)
\widehat{P_{k_3}\overline{w}}(s,-\xi)}{s+1}\big]ds\Big|\nonumber\\
&\qquad\qquad\qquad\qquad\qquad\qquad\qquad\qquad\lesssim \epsilon_1^3 2^{-\delta m}2^{-(N_1+10)k_+},
\end{align}
and  for $\iota_1\iota_2\iota_3\in\{+--,+++,---\}$,
\begin{align}
\Big|\int_{t_1}^{t_2}e^{i\vartheta(s,\xi)}I_{k_1k_2k_3}^{\iota_1\iota_2\iota_3}(s,\xi)ds\Big|\lesssim \epsilon_1^32^{-\delta m}2^{-(N_1+10)k_+}.
\label{sc10}
\end{align}

\begin{lem}\label{simpler lemma2} The estimate \eqref{sc10} holds  if  $k_1,k_2,k_3 \in [-4m,m/200-100]\cap \mathbb{Z}$ and
\begin{align}\label{sc5.15}
\min(k_1,k_2,k_3)+\mathrm{med}(k_1,k_2,k_3)\leq -6m/5.
\end{align}
If, in addition,
\begin{align}\label{sc5.16}
\max(|k_1-k|,|k_2-k|,|k_3-k|)\geq 21,
\end{align}
then  the estimate \eqref{sc9}  also holds.
\end{lem}

{\noindent \emph{Proof.}}
Under the condition \eqref{sc5.15}, we use \eqref{sc8},  the $L^\infty$ bound \eqref{Linfinity} to get
\begin{align*}
|I_{k_1k_2k_3}^{\iota_1\iota_2\iota_3}(s,\xi)|&\lesssim
\epsilon_1^3 2^{5\max(k_1,k_2,k_3)_+}\cdot2^{\min(k_1,k_2,k_3)}2^{\mathrm{med}(k_1,k_2,k_3)}\\
&\quad \cdot 2^{-(N_1+10)k_{1+}}\cdot2^{-(N_1+10)k_{2+}}\cdot2^{-(N_1+10)k_{3+}}\\
&\lesssim \epsilon_1^3 2^{(N_1+15)\max(k_1,k_2,k_3)_+}\cdot2^{-6 m/5}2^{-(N_1+10)k_+}\\
&\lesssim\epsilon_1^3 2^{-31 m/30}2^{-(N_1+10)k_+}
\end{align*}
for any $\iota_1\iota_2\iota_3\in\mathscr{T}$, where in the last step, we have also used
$$
(N_1+15)\max(k_1,k_2,k_3)\leq (N_1+15)m/200 <m/6.
$$
Therefore, the estimate \eqref{sc10} clearly holds. If, in addition,  \eqref{sc5.16} holds, then
$$
\varphi_{k}(\xi)\widehat{P_{k_1}w}(s,\xi)\widehat{P_{k_2}w}(s,\xi)
\widehat{P_{k_3}\overline{w}}(s,-\xi)= 0,
$$
so the estimate \eqref{sc9} follows.
\hfill $\Box$

In view of the above two lemmas, in order to prove  \eqref{sc5.1} and \eqref{sc5.2}, it suffices to show the following proposition.

\begin{prop}\label{scpropostion1}
Let $k\in\mathbb{Z}$, $|\xi|\sim 2^k$ and $t_1,t_2\in [2^{m}-2,2^{m+1}]\cap[0,T]$, $m\geq 25$ be an integer. Assume that $k_1,k_2,k_3$ satisfies
\begin{align}\label{sc5.17}
k_1,k_2,k_3 \in [-4m,m/200-100]\cap \mathbb{Z},
\end{align}
and
\begin{align} \label{sc5.18}
\min(k_1,k_2,k_3)+\mathrm{med}(k_1,k_2,k_3)\geq -6m/5.
\end{align}
Then the estimates \eqref{sc9} and \eqref{sc10} are valid.
\end{prop}

As mentioned before, in order to finish the proof of Proposition \ref{scproposition}, we shall also prove

\begin{prop}\label{scpropostion2}For any $0\leq t_1\leq t_2\leq T$, there exists $\delta>0$ such that
\begin{align}\label{sc7}
\Big|\langle\xi\rangle^{N_1+10}\int_{t_1}^{t_2}e^{i\vartheta(s,\xi)}
e^{is\langle\xi\rangle}\widehat{\mathcal{N}_R}(s,\xi)ds\Big|
\lesssim \epsilon_1^4(1+t_1)^{-\delta}.
\end{align}
\end{prop}

According to the above reductions, we see  Proposition \ref{scproposition} follows easily from  Propositions \ref{scpropostion1}--\ref{scpropostion2}. Hence, the remaining part of this subsection is devoted to the proofs of these two propositions. The bound \eqref{sc9} is proven through Lemmas \ref{lemsc5.5}--\ref{lemsc5.7} below, depending on different cases between the sizes of the input and output frequencies, and the bound \eqref{sc10} is obtained by Lemma \ref{lemsc5.8}. In addition, we will establish  the bound \eqref{sc7} with the help of Lemma \ref{lemmascr1}. In the proofs, we will frequently use the following multiplier lemma.

\begin{lem}\label{multiplierlemma}
There holds
\begin{align*}
\Big|\int_{\mathbb{R}^2 }m(\eta,\sigma)\widehat{f_1}(\eta)\widehat{f_2}(\sigma) \widehat{f_3}(-\eta-\sigma)d\eta d\sigma\Big|
\lesssim \|\mathscr{F}^{-1}m\|_{L^1}\|f_1\|_{L^{p_1}}\|f_2\|_{L^{p_2}}\|f_3\|_{L^{p_3}}
\end{align*}
with $p_1^{-1}+p_2^{-1}+p_3^{-1}=1$ and $p_1,\ p_2,\ p_3\in [1,+\infty]$.
\end{lem}

The proof of Lemma \ref{multiplierlemma} can be found in \cite{IPu1}. To bound the $L^1$ norm of $\mathscr{F}^{-1}m$, we usually use Lemma \ref{L1norm} below.

\begin{lem}\label{L1norm} If $m(\eta,\sigma)$ is a Fourier multiplier with $\eta$ and $\sigma$ localized in the size $2^k$ and $2^l$, respectively, and satisfies
\begin{align}\label{mb1}
|\partial_{\eta}^a\partial_{\sigma}^b m|\lesssim  A2^{-a k}2^{-b l}\  (\mathrm{resp.}\  A)
\end{align}
for any $a,b=0,1,2$, then we have
\begin{align}\label{mb2}
\|\mathscr{F}^{-1}m\|_{L^1(\mathbb{R}^2)}\lesssim A\ (\mathrm{resp.}\  A2^k2^l).
\end{align}
\end{lem}

\begin{lem}\label{integral lemma}  For any $\lambda,\ \mu>0$ and $n\in \mathbb{N}$, there holds that
\begin{align*}
\int_{\mathbb{R}^2}e^{  i\lambda xy}\varphi(\mu^{-1}x)\varphi(\mu^{-1}y)dxdy=2\pi\lambda^{-1}+\lambda^{-1-n}\mu^{-2n} O(1),
\end{align*}
 where $\varphi$ is the smooth radial function used in the Littlewood-Paley decomposition. The implicit constant coming from the term $O(1)$ depends only on $n$ and $\varphi$.
\end{lem}

Lemmas \ref{L1norm}--\ref{integral lemma} are proved in the appendix (see Lemma B.3 and Lemma B.4).

\begin{lem}\label{lemsc5.5} The estimate \eqref{sc9} holds provided that
\begin{align}\label{5.21}
\max(|k_1-k|,|k_2-k|,|k_3-k|)\leq 20.
\end{align}
\end{lem}

{\noindent \emph{Proof.}} It suffices  to prove, for any $s\in [t_1,t_2]$, that
\begin{align}\label{sc5.22}
\Big|I_{k_1k_2k_3}^{++-}(s,\xi)-2\pi
\frac{c^*(\xi)\langle\xi\rangle^3\widehat{P_{k_1}w}(s,\xi)\widehat{P_{k_2}w}(s,\xi)
\widehat{P_{k_3}\overline{w}}(s,-\xi)}{s+1}\Big|\lesssim\epsilon_1^32^{-(1+\delta_1)m}2^{-(N_1+10)k_+}
\end{align}
for some $\delta_1>0$. We split the proof into several steps.

\textbf{Step 1:}  $|\xi|\lesssim 2^{-m}$. In this case, we use \eqref{sc8} and the $L^\infty$ bound \eqref{Linfinity} to obtain
\begin{align*}
\mathrm{\mathrm{LHS\ of\ }} \eqref{sc5.22}&\lesssim  \epsilon_1^32^{2k}2^{5k_+}2^{-3(N_1+10)k_+}+\epsilon_1^32^{-m}2^{2k}2^{6k_+}2^{-3(N_1+10)k_+}\\
&\lesssim \epsilon_1^32^{-2m}2^{-(3N_1+24)k_+},
\end{align*}
which is better than the desired bound.

\textbf{Step 2:}  $|\xi|\gtrsim  2^{-m}$.  For the sake of convenience, we rewrite (by the change of variables $\eta\rightarrow -\eta$, $\sigma\rightarrow-\xi-\sigma-\eta$)
\begin{align*}
I_{k_1k_2k_3}^{++-}=\int_{\mathbb{R}^2}\widetilde{c}^{\ ++-}_{k_1 k_2 k_3}(\xi,\eta,\sigma)e^{is\Psi(\xi,\eta,\sigma)}
\widehat{P_{k_1}w}(\xi+\eta)\widehat{P_{k_2}w}(\xi+\sigma)
\widehat{P_{k_3}\overline{w}}(-\xi-\eta-\sigma)d\eta d\sigma,
\end{align*}
where
\begin{align}
\widetilde{c}^{\ ++-}(\xi,\eta,\sigma)&:=c(\xi,-\eta,-\xi-\sigma-\eta), \label{sc20}\\
\Psi(\xi,\eta,\sigma)&:=\langle\xi\rangle-\langle\xi+\eta\rangle-\langle\xi+\sigma\rangle
+\langle\xi+\eta+\sigma\rangle.\label{sc5.24}
\end{align}
Note that the set of space-time resonance  for  $\Psi$ now reduces to $(\xi,\eta,\sigma)=(\xi,0,0)$.  We see from \eqref{resonant1} and \eqref{sc20} that
$$
\widetilde{c}^{\ ++-}(\xi,0,0)=c(\xi,0,-\xi)=c^*(\xi).
$$
Let $\bar{l}$ be the smallest integer satisfying $2^{\bar{l}}\geq 2^{-9m/20}$. Note that $2^k\gtrsim 2^{-m}$ implies  $\bar{l}\leq k+10$. Now, we decompose
\begin{align}\label{5.26}
I_{k_1k_2k_3}^{++-}(s,\xi)=\sum_{l_1,l_2=\bar{l}}^{k+100}J_{{l_1}{l_2}}(s,\xi),
\end{align}
where
\begin{align*}
J_{{l_1}{l_2}}(s,\xi):=\int_{\mathbb{R}^2}\widetilde{c}^{\ ++-}_{k_1 k_2 k_3}e^{is\Psi}
\widehat{P_{k_1}w}(\xi+\eta)\widehat{P_{k_2}w}(\xi+\sigma)
\widehat{P_{k_3}\overline{w}}(-\xi-\eta-\sigma)\varphi_{l_1}^{(\bar{l})}(\eta)\varphi_{l_2}^{(\bar{l})}(\sigma)d\eta d\sigma,
\end{align*}
and
\begin{equation}\label{bump}
\varphi_l^{(l_0)}(\xi):=\left\{
\begin{array}{ll}
\varphi(|\xi|/2^l)-\varphi(|\xi|/2^{l-1}),&l\geq l_0+1,\\
\varphi(|\xi|/2^{l_0}),&l=l_0.
\end{array}
\right.
\end{equation}
In the following, we consider three different cases.

\textbf{Case 2a:} $\sigma$ is away from the space-time resonance set.  We aim to show that
\begin{align}\label{sc5.27}
|J_{{l_1}{l_2}}(s,\xi)|\lesssim \epsilon_1^32^{-m}2^{-\delta_1 m}2^{-(N_1+10)k_+},\ \ l_2\geq \max(l_1,\bar{l}+1).
\end{align}
From \eqref{sc5.24}, it is easy to see
\begin{align}\label{sc5.28}
|\partial_{\eta}\Psi|=\Big|-\frac{\xi+\eta}{\langle\xi+\eta\rangle}
+\frac{\xi+\eta+\sigma}{\langle\xi+\eta+\sigma\rangle}\Big|\gtrsim 2^{l_2}2^{-3k_+},
\end{align}
 whenever $|\xi+\eta|\sim |\xi+\sigma|\sim|\xi+\eta+\sigma|\sim 2^k$ and $|\sigma|\sim 2^{l_2}$. With integration by parts in $\eta$, we have
$$
|J_{{l_1}{l_2}}(s,\xi)|\leq |J_{{l_1}{l_2},1}(s,\xi)|+|J_{{l_1}{l_2},2}(s,\xi)|+|J_{{l_1}{l_2},3}(s,\xi)|,
$$
where
\begin{align*}
J_{{l_1}{l_2},1}&=\int_{\mathbb{R}^2}m_1e^{is\Psi}
\partial_\eta\widehat{P_{k_1}w}(s,\xi+\eta)\widehat{P_{k_2}w}(s,\xi+\sigma)
\widehat{P_{k_3}\overline{w}}(s,-\xi-\eta-\sigma)d\eta d\sigma,\\
J_{{l_1}{l_2},2}&=\int_{\mathbb{R}^2}m_1e^{is\Psi}
\widehat{P_{k_1}w}(s,\xi+\eta)\widehat{P_{k_2}w}(s,\xi+\sigma)
\partial_\eta\widehat{P_{k_3}\overline{w}}(s,-\xi-\eta-\sigma)d\eta d\sigma,\\
J_{{l_1}{l_2},3}&=\int_{\mathbb{R}^2}\partial_\eta m_1e^{is\Psi}
\widehat{P_{k_1}w}(s,\xi+\eta)\widehat{P_{k_2}w}(s,\xi+\sigma)
\widehat{P_{k_3}\overline{w}}(s,-\xi-\eta-\sigma)d\eta d\sigma,
\end{align*}
with
$$
m_1(\eta,\sigma):=\varphi_{l_1}^{(\bar{l})}(\eta)\varphi_{l_2}^{(\bar{l})}(\sigma)\cdot(s\partial_\eta\Psi)^{-1}\cdot \widetilde{c}^{\ ++-}_{k_1 k_2 k_3}.$$
 Using \eqref{sc5.28} and the fact $l_2\geq l_1$, we compute
\begin{align*}
|\partial_{\eta}^{a}\partial_{\sigma}^{b}[\varphi_{l_1}^{(\bar{l})}(\eta)\varphi_{l_2}^{(\bar{l})}(\sigma)
(s\partial_\eta\Psi)^{-1}]|\lesssim 2^{-m}2^{-l_2}2^{3k_+}2^{-a l_1}2^{-b l_2},\  \ \ a,b=0,1,2.
\end{align*}
Then by \eqref{mb1}--\eqref{mb2}, we have
$$
\|\mathscr{F}^{-1}[\varphi_{l_1}^{(\bar{l})}(\eta)\varphi_{l_2}^{(\bar{l})}(\sigma)
(s\partial_\eta\Psi)^{-1}]\|_{L^1(\mathbb{R}^2)} \lesssim 2^{-m}2^{-l_2}2^{3k_+}.
$$
Recalling the bound \eqref{sc8} for $\widetilde{c}^{\ ++-}_{k_1 k_2 k_3}$, we deduce from Lemma \ref{L1norm} that
$$
\|\mathscr{F}^{-1}\widetilde{c}^{\ ++-}_{k_1 k_2 k_3}\|_{L^1(\mathbb{R}^2)} \lesssim 2^{5k_+}2^{2k}.
$$
Combing the above two bounds give
\begin{align}\label{sc5.29}
\|\mathscr{F}^{-1}m_1\|_{L^1(\mathbb{R}^2)}&\lesssim\|\mathscr{F}^{-1}[\varphi_{l_1}^{(\bar{l})}(\eta)\varphi_{l_2}^{(\bar{l})}(\sigma)
(s\partial_\eta\Psi)^{-1}]\|_{L^1(\mathbb{R}^2)}\|\mathscr{F}^{-1}\widetilde{c}^{\ ++-}_{k_1 k_2 k_3}\|_{L^1(\mathbb{R}^2)}\nonumber\\
&\lesssim 2^{-m}2^{-l_2}2^{8k_+}2^{2k}.
\end{align}
Similarly, we can obtain
\begin{align*}
\|\mathscr{F}^{-1}(\partial_\eta m_1)\|_{L^1(\mathbb{R}^2)}&\lesssim 2^{-m}2^{-l_1}2^{-l_2}2^{8k_+}2^{2k}.
\end{align*}
Now, we apply Lemma \ref{multiplierlemma} with
\begin{align*}
\widehat{\alpha}(\eta)&:=e^{-is\langle\xi+\eta\rangle}\partial_\eta\widehat{P_{k_1}w}(s,\xi+\eta),,\\
\widehat{\beta}(\sigma)&:=e^{-is\langle\xi+\sigma\rangle}\widehat{P_{k_2}w}(s,\xi+\sigma),\ \ \ |\sigma|\sim 2^{l_2},\\
\widehat{\gamma}(\zeta)&:=e^{is\langle-\xi+\zeta\rangle}\widehat{P_{k_3}\overline{w}}(s,-\xi+\zeta),
\end{align*}
then
\begin{align*}
|J_{l_1l_2,1}|\lesssim \|\mathscr{F}^{-1}m_1\|_{L^1}\|\alpha\|_{L^2}\|\beta\|_{L^2}\|\gamma\|_{L^\infty}.
\end{align*}
Using the fact $|\sigma|\sim 2^{l_2}$, \eqref{L^2 of paritial derivative estimate}, \eqref{Linfinity} and \eqref{Linfinity of profile}, we see
\begin{align*}
\|\alpha\|_{L^2}\lesssim \epsilon_1 2^{p_0m} 2^{-(N_1-4)k_+},\ \ \|\beta\|_{L^2}\lesssim \epsilon_1 2^{l_2/2}2^{-(N_1+10)k_+},\ \  \|\gamma\|_{L^\infty}\lesssim \epsilon_1 2^{-m/2}.
\end{align*}
Therefore, these estimates and \eqref{sc5.29} lead to
\begin{align*}
|J_{l_1l_2,1}|&\lesssim 2^{-m}2^{-l_2}2^{8k_+}2^{2k}\cdot\epsilon_1 2^{p_0 m}2^{-(N_1-4)k_+}\cdot
\epsilon_1 2^{l_2/2}2^{-(N_1+10)k_+}
\cdot \epsilon_1 2^{-m/2}\\
&=\epsilon_1^3 2^{-3m/2}2^{p_0m}2^{-l_2/2}2^{2k}2^{-(2N_1-2)k_+}.
\end{align*}
Since $2^{-l_2/2}\lesssim 2^{9m/40}$, $J_{l_1l_2,1}$ can be  bounded by $\epsilon_1^3 2^{-51m/40}2^{p_0 m}2^{-(N_1+10)k_+}$.  With similar argument as above, we obtain  the same bound for $|J_{l_1l_2,2}|$. For the term $J_{l_1l_2,3}$, we apply Lemma \ref{multiplierlemma} with
\begin{align*}
\widehat{\tilde{\alpha}}(\eta)&:=e^{-is\langle\xi+\eta\rangle}\widehat{P_{k_1}w}(s,\xi+\eta),\ \ \ |\eta|\lesssim 2^{l_1},\\
\widehat{\tilde{\beta}}(\sigma)&:=e^{-is\langle\xi+\sigma\rangle}\widehat{P_{k_2}w}(s,\xi+\sigma),\ \ \ |\sigma|\sim 2^{l_2},\\
\widehat{\tilde{\gamma}}(\zeta)&:=e^{is\langle-\xi+\zeta\rangle}\widehat{P_{k_3}\overline{w}}(s,-\xi+\zeta)
\end{align*}
to obtain
\begin{align*}
|J_{l_1l_2,3}|&\lesssim \|\mathscr{F}^{-1}(\partial_\eta m_1)\|_{L^1}\|\tilde{\alpha}\|_{L^2}\|\tilde{\beta}\|_{L^2}\|\tilde{\gamma}\|_{L^\infty}\\
&\lesssim2^{-m}2^{-l_1}2^{-l_2}2^{8k_+}2^{2k}\cdot\epsilon_1 2^{l_1/2}2^{-(N_1+10)k_+}\cdot\epsilon_1 2^{l_2/2}2^{-(N_1+10)k_+}\cdot \epsilon_1 2^{-m/2}\\
&\lesssim \epsilon_1^3 2^{-m}2^{-m/2}2^{-l_1/2}2^{-l_2/2}2^{-(N_1+10)k_+}\\
&\lesssim \epsilon_1^3 2^{-m}2^{-m/20}2^{-(N_1+10)k_+}.
\end{align*}
Therefore, the estimate \eqref{sc5.27} is established.

\textbf{Case 2b: } $\eta$ is away from the space-time resonance set.  In this case, applying similar argument as above, we can prove that
\begin{align}\label{5.30}
|J_{{l_1}{l_2}}(s,\xi)|\lesssim\epsilon_1^32^{-m}2^{-\delta_1 m}2^{-(N_1+10)k_+},\ \ l_1\geq \max(l_2,\bar{l}+1).
\end{align}
Further details are omitted here since the proof is almost the same as Case 2a.

\textbf{Case 2c: } $(\eta,\sigma)$ is near the space-time resonance set. In this case, the above strategy is not workable as both $\eta$ and $\sigma$ can be very small, and a phase correction is needed to close the argument. Our aim is to show
\begin{align}\label{sc5.31}
\Big|J_{\bar{l}\ \bar{l}}(s,\xi)-2\pi
\frac{c^*(\xi)\langle\xi\rangle^3\widehat{P_{k_1}w}(s,\xi)\widehat{P_{k_2}w}(s,\xi)
\widehat{P_{k_3}\overline{w}}(s,-\xi)}{s+1}\Big|\lesssim \epsilon_1^32^{-m}2^{-\delta_1 m}2^{-(N_1+10)k_+}
\end{align}
for some $\delta_1>0$.  To prove \eqref{sc5.31}, we use
\begin{align*}
\mathrm{LHS\ of\ }\eqref{sc5.31}\leq & \big|J_{\bar{l}\ \bar{l}}(s,\xi)-\widetilde{J}_{\bar{l}\ \bar{l}}(s,\xi)\big|
+ \big|\widetilde{J}_{\bar{l}\ \bar{l}}(s,\xi)-\overline{J}_{\bar{l}\ \bar{l}}(s,\xi)\big|\\
&+\Big|\overline{J}_{\bar{l}\ \bar{l}}(s,\xi)-2\pi\frac{c^*(\xi)\langle\xi\rangle^3\widehat{P_{k_1}w}(s,\xi)\widehat{P_{k_2}w}(s,\xi)
\widehat{P_{k_3}\overline{w}}(s,-\xi)}{s+1}\Big|,
\end{align*}
where
\begin{align*}
\widetilde{J}_{\bar{l}\ \bar{l}}(s,\xi):=&\int_{\mathbb{R}^2}\widetilde{c}^{\ ++-}_{k_1 k_2 k_3}(\xi,\eta,\sigma)e^{\frac{is\eta\sigma}{\langle\xi\rangle^3}}\\
&\cdot\widehat{P_{k_1}w}(\xi+\eta)\widehat{P_{k_2}w}(\xi+\sigma)
\widehat{P_{k_3}\overline{w}}(-\xi-\eta-\sigma)\varphi(2^{-\bar{l}}\eta)\varphi(2^{-\bar{l}}\sigma)d\eta d\sigma,\\
\overline{J}_{\bar{l}\ \bar{l}}(s,\xi):=&\int_{\mathbb{R}^2}c^*(\xi)e^{\frac{is\eta\sigma}{\langle\xi\rangle^3}}
\widehat{P_{k_1}w}(\xi)\widehat{P_{k_2}w}(\xi)
\widehat{P_{k_3}\overline{w}}(-\xi)\varphi(2^{-\bar{l}}\eta)\varphi(2^{-\bar{l}}\sigma)d\eta d\sigma.
\end{align*}
By using Taylor's expansion, we have
\begin{align*}
\Psi(\xi,\eta,\sigma)=&\Psi(\xi,0,0)+\Psi_\eta(\xi,0,0)\eta+\Psi_\sigma(\xi,0,0)\sigma\\
&+\frac{1}{2}\Psi_{\eta\eta}(\xi,0,0)\eta^2+\frac{1}{2}\Psi_{\sigma\sigma}(\xi,0,0)\sigma^2
+\Psi_{\eta\sigma}(\xi,0,0)\eta\sigma+\mathrm{remainder},
\end{align*}
which implies, by \eqref{sc5.24} and the fact $|\eta|,|\sigma|\sim 2^{\bar{l}}$,
\begin{align*}
|\Psi(\xi,\eta,\sigma)-\langle\xi\rangle^{-3}\eta\sigma|\lesssim 2^{-4k_+}(|\eta|+|\sigma|)^3\lesssim 2^{-4k_+}2^{3\bar{l}}.
\end{align*}
Combining  \eqref{sc8}, \eqref{sc20} and \eqref{Linfinity}  yields ( $2^{\bar{l}}\sim 2^{-9m/20}$)
\begin{align}
\left|J_{\bar{l}\ \bar{l}}(s,\xi)-\widetilde{J}_{\bar{l}\ \bar{l}}(s,\xi)\right|
& \lesssim \int_{\mathbb{R}^2}|\widetilde{c}^{\ ++-}_{k_1 k_2 k_3}(\xi,\eta,\sigma)|\cdot|e^{is \Psi}-e^{\frac{is\eta\sigma}{\langle\xi\rangle^3}}|\nonumber\\
&\quad \cdot|\widehat{P_{k_1}w}(\xi+\eta)\widehat{P_{k_2}w}(\xi+\sigma)
\widehat{P_{k_3}\overline{w}}(-\xi-\eta-\sigma)|d\eta d\sigma\nonumber\\
&\lesssim
2^{5k_+}\cdot 2^m 2^{-4k_+}2^{3\bar{l}}\cdot\epsilon_1^3 2^{-3(N_1+10)k_+}\cdot2^{2\bar{l}}\nonumber\\
&\lesssim \epsilon_1^3 2^{-5m/4}2^{-(N_1+10)k_+}.\label{sc5.32}
\end{align}
In order to estimate the term $\widetilde{J}_{\bar{l}\ \bar{l}}(s,\xi)-\overline{J}_{\bar{l}\ \bar{l}}(s,\xi)$, note that
$$
|\widetilde{c}^{\ ++-}(\xi,\eta,\sigma)-c^*(\xi)|=|\widetilde{c}^{\ ++-}(\xi,\eta,\sigma)-\widetilde{c}^{\ ++-}(\xi,0,0)|
\lesssim 2^{\bar{l}}2^{5k_+},
$$
and by \eqref{L^2 of paritial derivative estimate},
$$
|\widehat{P_{k_1}w}(s,\xi+\zeta)-\widehat{P_{k_1}w}(s,\xi)|
\lesssim \|\partial_\xi\widehat{P_{k_1}w}\|_{L^2}2^{\bar{l}/2} \lesssim \epsilon_12^{p_0m}2^{-(N_1-4)k_+}2^{\bar{l}/2},\ \ \ |\zeta|\lesssim 2^{\bar{l}}.
$$
So it is easy to see
\begin{align*}
&\Big|\widetilde{c}^{\ ++-}(\xi,\eta,\sigma)\widehat{P_{k_1}w}(s,\xi+\eta)\widehat{P_{k_2}w}(s,\xi+\sigma)
\widehat{P_{k_3}\overline{w}}(s,-\xi-\eta-\sigma)\\
&\quad\quad\quad-c^*(\xi)\widehat{P_{k_1}w}(s,\xi)\widehat{P_{k_2}w}(s,\xi)
\widehat{P_{k_3}\overline{w}}(s,-\xi)\Big|\\
&\qquad\lesssim \epsilon_1^32^{\bar{l}}2^{-(3N_1+25)k_+}+ \epsilon_1^32^{\bar{l}/2}2^{p_0m}2^{-(3N_1+11)k_+}
\end{align*}
whenever $|\eta|,|\sigma|\lesssim 2^{\bar{l}}$, where we have used \eqref{Linfinity} in the above estimate. Therefore
\begin{align}\label{sc5.33}
\big|\widetilde{J}_{\bar{l}\ \bar{l}}(s,\xi)-\overline{J}_{\bar{l}\ \bar{l}}(s,\xi)\big|&\lesssim \epsilon_1^32^{3\bar{l}}2^{-(3N_1+25)k_+}+ \epsilon_1^32^{5\bar{l}/2}2^{p_0m}2^{-(3N_1+11)k_+}\nonumber\\
&\lesssim \epsilon_1^32^{-m}2^{-m/8}2^{p_0m}2^{-(N_1+10)k_+}.
\end{align}
Now, using \eqref{Linfinity} and applying Lemma \ref{integral lemma} with $\lambda=s/\langle\xi\rangle^3$, $\mu=2^{\bar{l}}$ and $n=1$, we have
\begin{align}
&\Big|\overline{J}_{\bar{l}\ \bar{l}}(s,\xi)-2\pi\frac{c^*(\xi)\langle\xi\rangle^3\widehat{P_{k_1}w}(s,\xi)\widehat{P_{k_2}w}(s,\xi)
\widehat{P_{k_3}\overline{w}}(s,-\xi)}{s+1}\Big|\nonumber\\
&\quad\quad\lesssim \Big|\overline{J}_{\bar{l}\ \bar{l}}(s,\xi)-2\pi\frac{c^*(\xi)\langle\xi\rangle^3\widehat{P_{k_1}w}(s,\xi)\widehat{P_{k_2}w}(s,\xi)
\widehat{P_{k_3}\overline{w}}(s,-\xi)}{s}\Big|\nonumber\\
&\quad\quad\quad+ |c^*(\xi)\langle\xi\rangle^3\widehat{P_{k_1}w}(s,\xi)\widehat{P_{k_2}w}(s,\xi)
\widehat{P_{k_3}\overline{w}}(s,-\xi)|\cdot(\frac{1}{s}-\frac{1}{s+1})\nonumber\\
&\quad\quad \lesssim  \epsilon_1^32^{2k}2^{-(3N_1+27)k_+}\Big|\int_{\mathbb{R}^2}e^{\frac{is\eta\sigma}{\langle\xi\rangle^3}}
\varphi(2^{-\bar{l}}\eta)\varphi(2^{-\bar{l}}\sigma)d\eta d\sigma-2\pi\langle\xi\rangle^3s^{-1}\Big|\nonumber\\
&\quad\quad\quad +\epsilon_1^32^{2k}2^{-(3N_1+24)k_+}2^{-2m}\nonumber\\
&\quad\quad \lesssim  \epsilon_1^32^{-(3N_1+25)k_+}2^{-2m}2^{6k_+}2^{-2\bar{l}}+\epsilon_1^3 2^{-(3N_1+22)k_+}2^{-2m}\nonumber\\
&\quad\quad \lesssim  \epsilon_1^32^{-m}2^{-m/10}2^{-(N_1+10)k_+}. \label{sc5.34}
\end{align}
Therefore, \eqref{sc5.31} follows from \eqref{sc5.32}--\eqref{sc5.34}. This ends the proof of the lemma.
\hfill $\Box$

\begin{lem}\label{lemsc5.6} The estimate \eqref{sc9} holds  under the conditions \eqref{sc5.17}, \eqref{sc5.18} and
\begin{align}\label{sc5.35}
\max(|k_1-k|,|k_2-k|,|k_3-k|)\geq 21,\ \ \ \ \ \max(|k_1-k_3|,|k_2-k_3|)\geq 6.
\end{align}
\end{lem}

{\noindent \emph{Proof.}}
Recall that
\begin{align*}
I_{k_1k_2k_3}^{++-}=\int_{\mathbb{R}^2}\widetilde{c}^{\ ++-}_{k_1k_2k_3}(\xi,\eta,\sigma)e^{is\Psi(\xi,\eta,\sigma)}
\widehat{P_{k_1}w}(\xi+\eta)\widehat{P_{k_2}w}(\xi+\sigma)
\widehat{P_{k_3}\overline{w}}(-\xi-\eta-\sigma)d\eta d\sigma,
\end{align*}
where
\begin{align*}
\Psi(\xi,\eta,\sigma)=\langle\xi\rangle-\langle\xi+\eta\rangle-\langle\xi+\sigma\rangle
+\langle\xi+\eta+\sigma\rangle,
\end{align*}
and our aim is to show that there exists $\delta_2>0$ such that
\begin{align}\label{sc5.36}
|I_{k_1k_2k_3}^{++-}(s,\xi)|\lesssim \epsilon_1^32^{-m}2^{-\delta_2 m}2^{-(N_1+10)k_+}.
\end{align}

According to \eqref{sc5.35}, we may assume $|k_1-k_3|\geq 6$. Since $-\sigma=(\xi+\eta)+(-\xi-\eta-\sigma)$,  then we have $|\sigma|\sim \max(|\xi+\eta|,|\xi+\eta+\sigma|)=2^{\max({k_1},k_3)}$ and
\begin{align}
|\partial_\eta \Psi|&=\Big|-\frac{\xi+\eta}{\langle\xi+\eta\rangle}+\frac{\xi+\eta+\sigma}{\langle\xi+\eta+\sigma\rangle}\Big|
\gtrsim 2^{-3\max(k_1,k_3)_{+}}2^{\max(k_1,k_3)},\label{sc5.37}\\
|\partial_\eta^2 \Psi|&=\Big|-\frac{1}{\langle\xi+\eta\rangle^3}+\frac{1}{\langle\xi+\eta+\sigma\rangle^3}\Big|
\lesssim 2^{-5\min(k_1,k_3)_{+}}2^{2\max(k_1,k_3)}.\label{sc5.38}
\end{align}
Integration by parts with respect to $\eta$ gives
\begin{align*}
|I_{k_1k_2k_3}^{++-}(s,\xi)|\leq |F_1(s,\xi)|+|F_2(s,\xi)|+|F_3(s,\xi)|,
 \end{align*}
where
\begin{align*}
F_1(s,\xi)&:=\int_{\mathbb{R}^2}e^{is\Psi}m_2
\partial_\eta\widehat{P_{k_1}w}(\xi+\eta)\widehat{P_{k_2}w}(\xi+\sigma)
\widehat{P_{k_3}\overline{w}}(-\xi-\eta-\sigma)d\eta d\sigma,\\
F_2(s,\xi)&:=\int_{\mathbb{R}^2}e^{is\Psi}m_2
\widehat{P_{k_1}w}(\xi+\eta)\widehat{P_{k_2}w}(\xi+\sigma)
\partial_\eta{\widehat{P_{k_3}\overline{w}}}(-\xi-\eta-\sigma)d\eta d\sigma,\\
F_3(s,\xi)&:=\int_{\mathbb{R}^2} e^{is\Psi}\partial_\eta m_2
\widehat{P_{k_1}w}(\xi+\eta)\widehat{P_{k_2}w}(\xi+\sigma)
\widehat{P_{k_3}\overline{w}}(-\xi-\eta-\sigma)d\eta d\sigma
\end{align*}
with
$$m_2=m_2(\eta,\sigma):=( s\partial_\eta \Psi)^{-1}\cdot\widetilde{c}^{\ ++-}_{k_1k_2k_3}. $$
Using the bounds \eqref{sc8}, \eqref{sc5.37}, \eqref{sc5.38} and Lemma \ref{L1norm}, we can obtain
\begin{align}
&\|\mathscr{F}^{-1}m_2\|_{L^1(\mathbb{R}^2)}\lesssim 2^{-m}2^{10\max(k_1,k_3)_{+}}2^{-\max(k_1,k_3)},\label{sc5.39}\\
&\|\mathscr{F}^{-1}(\partial_\eta m_2)\|_{L^1(\mathbb{R}^2)}\lesssim 2^{-m}2^{13\max(k_1,k_3)_{+}}2^{-5\min(k_1,k_3)_{+}}.
\label{sc5.40}
\end{align}
Applying Lemma \ref{multiplierlemma} with
\begin{align*}
\widehat{\alpha}(\eta)&:=e^{-is\langle\xi+\eta\rangle}\partial_\eta\widehat{P_{k_1}w}(s,\xi+\eta),\\
\widehat{\beta}(\sigma)&:=e^{-is\langle\xi+\sigma\rangle}\widehat{P_{k_2}w}(s,\xi+\sigma),\\
\widehat{\gamma}(\zeta)&:=e^{is\langle-\xi+\zeta\rangle}\widehat{P_{k_3}\overline{w}}(s,-\xi+\zeta),
\end{align*}
we use \eqref{sc5.39}, \eqref{L^2 of paritial derivative estimate}, \eqref{L^2 of low} and  \eqref{Linfinity of profile} to get
\begin{align*}
|F_1(s,\xi)|&\lesssim \|\mathscr{F}^{-1}m_2\|_{L^1}\|\alpha\|_{L^2}\|\beta\|_{L^2}\|\gamma\|_{L^\infty}\\
&\lesssim 2^{-m}2^{10\max(k_1,k_3)_{+}}2^{-\max(k_1,k_3)}\cdot \epsilon_12^{p_0m}2^{-(N_1-4)k_{1+}}\\
&\quad \cdot
\epsilon_1 2^{\max(k_1,k_3)/2}2^{-(N_1+10)k_{2+}}\cdot \epsilon_1 2^{-m/2}\\
&\lesssim \epsilon_1^32^{-3m/2}2^{p_0m}2^{-(N_1+10)k_+}2^{(N_1+10)k_+}2^{-\max(k_1,k_3)/2}.
\end{align*}
Notice that  $k\leq \max(k_1,k_2,k_3)+2$, so the assumptions \eqref{sc5.17}--\eqref{sc5.18} yield
\begin{align*}
(N_1+10)\max(k_1,k_2,k_3)\leq m/6,\ \ \ \  \mathrm{med}(k_1,k_2,k_3)\geq -3m/5.
\end{align*}
Since
\begin{align*}
\max(k_1,k_3)&= \mathrm{med}(k_1,k_2,k_3),\  \mathrm{if}\  k_2=\max(k_1,k_2,k_3),\\
\max(k_1,k_3)&\geq  \mathrm{med}(k_1,k_2,k_3),\ \mathrm{if}\ k_2<\max(k_1,k_2,k_3),
\end{align*}
we also get $\max(k_1,k_3)\geq -3m/5$. Therefore, we conclude
\begin{align*}
|F_1(s,\xi)|&\lesssim \epsilon_1^32^{-3m/2}2^{p_0m}2^{-(N_1+10)k_+}2^{(N_1+10)\max(k_1,k_2,k_3)_{+}}2^{-\max(k_1,k_3)/2}\\
&\lesssim \epsilon_1^32^{-m}2^{p_0m}2^{-m/30}2^{-(N_1+10)k_+}.
\end{align*}
With the same treatment, we can get the same bound for $|F_2(s,\xi)|$. Finally, using Lemma \ref{multiplierlemma}, \eqref{sc5.40}, \eqref{L^2estimate} and \eqref{Linfinity of profile}, we can obtain
\begin{align*}
|F_3(s,\xi)|&\lesssim 2^{-m}2^{13\max(k_1,k_3)_{+}}2^{-5\min(k_1,k_3)_+}\cdot \epsilon_1^2 2^{2p_0m}
2^{-(N-5)\max(k_1,k_2,k_3)_+}\cdot \epsilon_1 2^{-m/2}\\
&\lesssim \epsilon_1^32^{-m}2^{2p_0m}2^{-m/2}2^{-(N_1+10)k_+}.
\end{align*}
By combining the estimates for $F_1$, $F_2$ and $F_3$, we deduce the desired bound \eqref{sc5.36}.
\hfill $\Box$

\begin{lem}\label{lemsc5.7} The estimate \eqref{sc9} holds under the hypotheses \eqref{sc5.17}, \eqref{sc5.18} and
\begin{align}\label{sc5.41}
\max(|k_1-k|,|k_2-k|,|k_3-k|)\geq 21, \ \ \ \
\max(|k_1-k_3|,|k_2-k_3|)\leq 5.
\end{align}
\end{lem}

{\noindent \emph{Proof.}}
Recall  that we want to show
\begin{align}\label{sc5.44}
|I_{k_1k_2k_3}^{++-}(s,\xi)|\lesssim \epsilon_1^32^{-m}2^{-\delta_3 m}2^{-(N_1+10)k_+}
\end{align}
for some $\delta_3>0$, where the definition of $I_{k_1k_2k_3}^{++-}$ is the same as in Lemma \ref{lemsc5.6}. According to \eqref{sc5.41}, we may assume $k_1,k_2,k_3\geq k+11$, then it follows from \eqref{sc5.18} that
\begin{align}\label{sc5.45}
2^{k_1}\sim 2^{k_2}\sim 2^{k_3}\gtrsim 2^{-3m/5}.
\end{align}
Since $\eta=(\xi+\eta)-\xi$ and $\sigma=(\xi+\sigma)-\xi$, we also have $|\eta|\sim |\sigma|\sim 2^{k_1}$.  Therefore,
\begin{align*}
|\partial_\eta \Psi|&=\Big|-\frac{\xi+\eta}{\langle\xi+\eta\rangle}+\frac{\xi+\eta+\sigma}{\langle\xi+\eta+\sigma\rangle}\Big|
\sim  2^{-3k_{1+}}2^{k_1},\\
|\partial_\eta^2 \Psi|&=\Big|-\frac{1}{\langle\xi+\eta\rangle^3}+\frac{1}{\langle\xi+\eta+\sigma\rangle^3}\Big|
\sim 2^{-5k_{1+}}2^{2k_1}.
\end{align*}
Now, with integration by parts in $\eta$, we see
\begin{align*}
|I_{k_1k_2k_3}^{++-}(s,\xi)|\leq |G_1(s,\xi)|+|G_2(s,\xi)|+|G_3(s,\xi)|,
 \end{align*}
where
\begin{align*}
G_1(s,\xi)&:=\int_{\mathbb{R}^2}e^{is\Psi}m_3
\partial_\eta\widehat{P_{k_1}w}(\xi+\eta)\widehat{P_{k_2}w}(\xi+\sigma)
\widehat{P_{k_3}\overline{w}}(-\xi-\eta-\sigma)d\eta d\sigma,\\
G_2(s,\xi)&:=\int_{\mathbb{R}^2}e^{is\Psi}m_3
\widehat{P_{k_1}w}(\xi+\eta)\widehat{P_{k_2}w}(\xi+\sigma)
\partial_\eta{P_{k_3}\overline{w}}(-\xi-\eta-\sigma)d\eta d\sigma,\\
G_3(s,\xi)&:=\int_{\mathbb{R}^2} e^{is\Psi}\partial_\eta m_3
\widehat{P_{k_1}w}(\xi+\eta)\widehat{P_{k_2}w}(\xi+\sigma)
\widehat{P_{k_3}\overline{w}}(-\xi-\eta-\sigma)d\eta d\sigma
\end{align*}
with
$$m_3=m_3(\eta,\sigma):=( s\partial_\eta \Psi)^{-1}\cdot \widetilde{c}^{\ ++-}_{k_1k_2k_3}.$$
From \eqref{sc8}, Lemma \ref{L1norm} and the bounds for $\partial_\eta \Psi$ and $\partial^2_\eta \Psi$, it is easy to see
\begin{align*}
\|\mathscr{F}^{-1}m_3\|_{L^1(\mathbb{R}^2)}\lesssim 2^{-m}2^{10k_{1+}}2^{-k_1},\ \ \ \ \ \|\mathscr{F}^{-1}(\partial_\eta m_2)\|_{L^1(\mathbb{R}^2)}\lesssim  2^{-m}2^{8k_{1+}}.
\end{align*}
Applying Lemma \ref{multiplierlemma} with
\begin{align*}
\widehat{\alpha}(\eta)&:=e^{-is\langle\xi+\eta\rangle}\partial_\eta\widehat{P_{k_1}w}(s,\xi+\eta),\\
\widehat{\beta}(\sigma)&:=e^{-is\langle\xi+\sigma\rangle}\widehat{P_{k_2}w}(s,\xi+\sigma),\\
\widehat{\gamma}(\zeta)&:=e^{is\langle-\xi+\zeta\rangle}\widehat{P_{k_3}\overline{w}}(s,-\xi+\zeta),
\end{align*}
and using \eqref{sc5.45}, we deduce
\begin{align*}
|G_1(s,\xi)|&\lesssim \|\mathscr{F}^{-1}m_3\|_{L^1(\mathbb{R}^2)}\|\alpha\|_{L^2}\|\beta\|_{L^2}\|\gamma\|_{L^\infty}\\
&\lesssim  2^{-m}2^{10k_{1+}}2^{-k_1}\cdot \epsilon_12^{p_0m}2^{-(N_1-4)k_{1+}}\cdot \epsilon_12^{k_2/2}2^{-(N_1+10)k_{2+}}\cdot \epsilon_1 2^{-m/2}\\
&\lesssim \epsilon_1^32^{-m}2^{p_0m}2^{-m/30}2^{-(N_1+10)k_+}.
\end{align*}
Similarly, we can obtain the same bound for the term $|G_2(s,\xi)|$.  To estimate $|G_3(s,\xi)|$, we again use Lemma \ref{multiplierlemma} to get
\begin{align*}
|G_3(s,\xi)|&\lesssim  2^{-m}2^{8k_{1+}}\cdot  \epsilon_12^{p_0m}2^{-(N-5)k_{1+}}\cdot  \epsilon_12^{p_0m}2^{-(N-5)k_{2+}}\cdot \epsilon_1 2^{-m/2} \\
&\lesssim  \epsilon_1^3 2^{-3m/2}2^{2p_0m}2^{-(N_1+10)k_+}.
\end{align*}
The proof of Lemma \ref{lemsc5.7} is completed.
\hfill $\Box$

\begin{lem}\label{lemsc5.8}
The estimate \eqref{sc10} holds under the assumptions of Proposition \ref{scpropostion1}.
\end{lem}

{\noindent \emph{Proof.}}
By a simple change of variables, we rewrite the LHS of \eqref{sc10} as
\begin{align*}
&\int_{t_1}^{t_2}e^{i\vartheta(s,\xi)}I_{k_1k_2k_3}^{\iota_1\iota_2\iota_3}(s,\xi)ds\\
& =\int_{t_1}^{t_2}e^{i\vartheta}\Big[\int_{\mathbb{R}^2}\widetilde{c}^{\ \iota_1\iota_2\iota_3}_{k_1 k_2 k_3}e^{is\Psi^{\iota_1\iota_2\iota_3}}
\widehat{P_{k_1}w^{\iota_1}}(\xi+\eta)\widehat{P_{k_2}w^{\iota_2}}(\xi+\sigma)
\widehat{P_{k_3}w^{\iota_3}}(-\xi-\eta-\sigma)d\eta d\sigma\Big]ds,
\end{align*}
where
\begin{align*}
\widetilde{c}^{\ \iota_1\iota_2\iota_3}_{k_1 k_2 k_3}(\xi,\eta,\sigma)&
:=c^{\iota_1\iota_2\iota_3}_{k_1 k_2 k_3}(\xi+\eta,\xi+\sigma,-\xi-\eta-\sigma),\\
\Psi^{\iota_1\iota_2\iota_3}(\xi,\eta,\sigma)
&=\langle\xi\rangle-\iota_1\langle\xi+\eta\rangle-\iota_2\langle\xi+\sigma\rangle
-\iota_3\langle\xi+\eta+\sigma\rangle.
\end{align*}
Note that the phase $\Psi^{\iota_1\iota_2\iota_3}(\xi,\eta,\sigma)$ never vanishes when $\iota_1\iota_2\iota_3\in\{+--,+++,---\}$. So we use integration by parts in $s$ to obtain
\begin{align*}
\int_{t_1}^{t_2}e^{iH(s,\xi)}I_{k_1k_2k_3}^{\iota_1\iota_2\iota_3}(s,\xi)ds:=K_1(t_1,\xi)+K_2(t_2,\xi)+L_1(\xi)+L_2(\xi),
\end{align*}
where
\begin{align*}
&K_1(t_1,\xi):=\\
&-
e^{i\vartheta(t_1,\xi)}\int_{\mathbb{R}^2}\frac{\widetilde{c}^{\ \iota_1\iota_2\iota_3}_{k_1 k_2 k_3}e^{it_1\Psi^{\iota_1\iota_2\iota_3}}}{i\Psi^{\iota_1\iota_2\iota_3}}
\widehat{P_{k_1}w^{\iota_1}}(t_1,\xi+\eta)\widehat{P_{k_2}w^{\iota_2}}(t_1,\xi+\sigma)
\widehat{P_{k_3}w^{\iota_3}}(t_1,-\xi-\eta-\sigma)d\eta d\sigma,\\
&K_2(t_2,\xi):=\\
&\quad
e^{i\vartheta(t_2,\xi)}\int_{\mathbb{R}^2}\frac{\widetilde{c}^{\ \iota_1\iota_2\iota_3}_{k_1 k_2 k_3}e^{it_2\Psi^{\iota_1\iota_2\iota_3}}}{i\Psi^{\iota_1\iota_2\iota_3}}
\widehat{P_{k_1}w^{\iota_1}}(t_2,\xi+\eta)\widehat{P_{k_2}w^{\iota_2}}(t_2,\xi+\sigma)
\widehat{P_{k_3}w^{\iota_3}}(t_2,-\xi-\eta-\sigma)d\eta d\sigma,
\end{align*}
and
\begin{align*}
&L_1(\xi):=\\
&-\int_{t_1}^{t_2}e^{i\vartheta}\vartheta_s\Big[\int_{\mathbb{R}^2}\frac{\widetilde{c}^{\ \iota_1\iota_2\iota_3}_{k_1 k_2 k_3}e^{is\Psi^{\iota_1\iota_2\iota_3}}}{\Psi^{\iota_1\iota_2\iota_3}}
\widehat{P_{k_1}w^{\iota_1}}(\xi+\eta)\widehat{P_{k_2}w^{\iota_2}}(\xi+\sigma)
\widehat{P_{k_3}w^{\iota_3}}(-\xi-\eta-\sigma)d\eta d\sigma\Big]ds,\\
&L_2(\xi):=\\
&-\int_{t_1}^{t_2}e^{i\vartheta}\Big[\int_{\mathbb{R}^2}\frac{\widetilde{c}^{\ \iota_1\iota_2\iota_3}_{k_1 k_2 k_3}e^{is\Psi^{\iota_1\iota_2\iota_3}}}{i\Psi^{\iota_1\iota_2\iota_3}}
\partial_s\widehat{P_{k_1}w^{\iota_1}}(\xi+\eta)\widehat{P_{k_2}w^{\iota_2}}(\xi+\sigma)
\widehat{P_{k_3}w^{\iota_3}}(-\xi-\eta-\sigma)d\eta d\sigma\Big]ds\\
&-\int_{t_1}^{t_2}e^{i\vartheta}\Big[\int_{\mathbb{R}^2}\frac{\widetilde{c}^{\ \iota_1\iota_2\iota_3}_{k_1 k_2 k_3}e^{is\Psi^{\iota_1\iota_2\iota_3}}}{i\Psi^{\iota_1\iota_2\iota_3}}
\widehat{P_{k_1}w^{\iota_1}}(\xi+\eta)\partial_s\widehat{P_{k_2}w^{\iota_2}}(\xi+\sigma)
\widehat{P_{k_3}w^{\iota_3}}(-\xi-\eta-\sigma)d\eta d\sigma\Big]ds\\
&-\int_{t_1}^{t_2}e^{i\vartheta}\Big[\int_{\mathbb{R}^2}\frac{\widetilde{c}^{\ \iota_1\iota_2\iota_3}_{k_1 k_2 k_3}e^{is\Psi^{\iota_1\iota_2\iota_3}}}{i\Psi^{\iota_1\iota_2\iota_3}}
\widehat{P_{k_1}w^{\iota_1}}(\xi+\eta)\widehat{P_{k_2}w^{\iota_2}}(\xi+\sigma)
\partial_s\widehat{P_{k_3}w^{\iota_3}}(-\xi-\eta-\sigma)d\eta d\sigma\Big]ds.
\end{align*}
Hence, in order to establish this lemma, it suffices to prove that there exists $\delta_4>0$ such that
\begin{align}\label{sc5.47}
|K_1(t_1,\xi)|+|K_2(t_2,\xi)|+ |L_1(\xi)|+ |L_2(\xi)|\lesssim \epsilon_1^32^{-\delta_4 m}2^{-(N_1+10)k_+}
\end{align}
whenever  $|\xi|\sim 2^k$ and $\iota_1\iota_2\iota_3\in\{+--,+++,---\}$.

We first prove \eqref{sc5.47} for the case $\iota_1\iota_2\iota_3=+--$. It is easy to see
$$
|\Psi^{+--}|^{-1}\lesssim\langle\xi\rangle+\langle\xi+\eta\rangle+\langle\xi+\sigma\rangle
+\langle\xi+\eta+\sigma\rangle\lesssim \langle\xi+\eta\rangle+\langle\xi+\sigma\rangle
+\langle\xi+\eta+\sigma\rangle,
$$
then using \eqref{sc8} and Lemma \ref{L1norm},  we can see
\begin{align}\label{sc5.48}
\|\mathscr{F}^{-1}[(\Psi^{\iota_1\iota_2\iota_3})^{-1}\widetilde{c}^{\ \iota_1\iota_2\iota_3}_{k_1 k_2 k_3}]\|_{L^1(\mathbb{R}^2)}\lesssim 2^{8\max(k_1,k_2,k_3)_+}.
\end{align}
Applying Lemma \ref{multiplierlemma} with
\begin{align*}
\widehat{\alpha}(\eta)&:=e^{-it_j\langle\xi+\eta\rangle}\widehat{P_{k_1}w}(t_j,\xi+\eta),\\
\widehat{\beta}(\sigma)&:=e^{it_j\langle\xi+\sigma\rangle}\widehat{P_{k_2}\overline{w}}(t_j,\xi+\sigma),\\
\widehat{\gamma}(\zeta)&:=e^{it_j\langle-\xi+\zeta\rangle}\widehat{P_{k_3}\overline{w}}(t_j,-\xi+\zeta),
\end{align*}
for $j=1,2$, and using \eqref{sc5.48}, \eqref{L^2estimate}, \eqref{Linfinity of profile}, we can obtain, by estimating the lowest frequency component in $L^\infty$ and the other two components in $L^2$,
\begin{align*}
&|K_1(t_1,\xi)|+|K_2(t_2,\xi)|\\
&\quad\lesssim  2^{8\max(k_1,k_2,k_3)_+}\cdot \epsilon_1 2^{p_0m}2^{-(N-5)\max(k_1,k_2,k_3)_+}
\cdot\epsilon_1 2^{p_0 m}2^{-(N-5)\mathrm{med}(k_1,k_2,k_3)_+}\cdot \epsilon_12^{-m/2}\\
&\quad\lesssim \epsilon_1^32^{-m/2}2^{2p_0m}2^{-(N_1+10)k_+}.
\end{align*}
To estimate $L_1(\xi)$, note that
\begin{align}\label{sc5.49}
|\vartheta_s(s,\xi)|\lesssim c^*(\xi)\langle\xi\rangle^3(1+s)^{-1}|\widehat{w}(s,\xi)|^2
\lesssim \epsilon_1^22^{-m}2^{2k}2^{6k_+}2^{-(2N_1+20)k_+},\ |\xi|\sim 2^k,
\end{align}
then using Lemma \ref{multiplierlemma}, \eqref{sc5.48}, \eqref{sc5.49}, \eqref{L^2estimate} and \eqref{Linfinity of profile},  we obtain
\begin{align*}
|L_1(\xi)|\lesssim & 2^m \cdot \epsilon_1^22^{-m}2^{-(2N_1+12)k_+}\cdot 2^{8\max(k_1,k_2,k_3)_+}\cdot\epsilon_12^{p_0m}2^{-(N-5)\max(k_1,k_2,k_3)_+}\\
&\cdot\epsilon_12^{p_0m}2^{-(N-5)\mathrm{med}(k_1,k_2,k_3)_+}\cdot \epsilon_12^{-m/2}\\
\lesssim&\epsilon_1^52^{-m/2}2^{2p_0m}2^{-(N_1+10)k_+}.
\end{align*}
 For the term $L_2(\xi)$,  we use \eqref{L^2estimate}, \eqref{L^2 of partial time} to  get
 \begin{align*}
|L_2(\xi)|\lesssim &2^m \cdot 2^{8\max(k_1,k_2,k_3)_+}
\cdot\epsilon_1^52^{3p_0m}2^{-3m/2}\cdot(2^{-(N-7)\max(k_1,k_2,k_3)_+}2^{-(N-7)\mathrm{med}(k_1,k_2,k_3)_+}\\
&+2^{-(N-7)\min(k_1,k_2,k_3)_+}2^{-(N-7)\max(k_1,k_2,k_3)_+})
 \\
\lesssim& \epsilon_1^52^{-m/2}2^{3p_0m}2^{-(N_1+10)k_+}.
\end{align*}
Therefore, the estimate \eqref{sc5.47} is established for $\iota_1\iota_2\iota_3=+--$.

Note that
\begin{align*}
|\Psi^{+++}(\xi,\eta,\sigma)|&\geq \frac{2}{\langle\xi\rangle+\langle\xi+\eta\rangle+\langle\xi+\sigma\rangle
+\langle\xi+\eta+\sigma\rangle},\\
|\Psi^{---}(\xi,\eta,\sigma)|&\sim  \max(\langle\xi\rangle,\langle\xi+\eta\rangle,\langle\xi+\sigma\rangle
,\langle\xi+\eta+\sigma\rangle).
\end{align*}
Then we can apply the same argument as above to show the bound \eqref{sc5.47} in the case $\iota_1\iota_2\iota_3=+++$ and $--- $. For the sake of simplicity, we omit further details. This ends the proof of the lemma.
\hfill $\Box$

To complete the proof of Proposition \ref{scproposition}, we are left to prove \eqref{sc7}. That is, we are aiming to show
\begin{align*}
\left|\langle\xi\rangle^{N_1+10}\int_{t_1}^{t_2}e^{i\vartheta(s,\xi)}e^{is\langle\xi\rangle}
\widehat{\mathcal{N}_R}(s,\xi)ds\right|
\lesssim \epsilon_1^4(1+t_1)^{-\delta}
\end{align*}
for some $\delta>0$. Recall that the definitions of $\vartheta$ and $\mathcal{N}_R$  are given by \eqref{definition of theta} and  \eqref{sc-a1}, respectively. To prove this bound, we need the following lemma.

\begin{lem}\label{lemmascr1}
For any $t\in [0,T]$, there hold that
\begin{align}
 \|\langle\xi\rangle^{N-20}\widehat{\mathcal{N}_R}(\xi)\|_{L^\infty}&\lesssim \epsilon_1^4 (1+t)^{2p_0-1},\label{scr4}\\
\|\mathcal{N}_R\|_{H^{N_1+10}}&\lesssim\epsilon_1^4 (1+t)^{p_0-3/2},\label{scr5}\\
\|\Gamma\mathcal{N}_R\|_{H^{N_1-10}}&\lesssim\epsilon_1^4 (1+t)^{p_0-3/2},\label{scr6}\\
\|x\mathcal{N}_R\|_{H^{N_1-10}}&\lesssim\epsilon_1^4 (1+t)^{p_0-1/2}.\label{scr7}
\end{align}
\end{lem}

Recall  the  bounds for $g$ and $h$
\begin{align}\label{scr2}
\begin{split}
\|g\|_{H^{N-5}}+ \|h\|_{H^N} &\lesssim \epsilon_1(1+t)^{p_0},\\
\|\Gamma g\|_{H^{N_1-5}}+\|\Gamma h\|_{H^{N_1}}&\lesssim \epsilon_1(1+t)^{p_0},\\
\|g\|_{W^{N_1+10,\infty}}+\|h\|_{W^{N_1+10,\infty}}&\lesssim \epsilon_1(1+t)^{-1/2},
\end{split}
\end{align}
and the bounds for the difference $h-g$
\begin{align}\label{scr3}
\begin{split}
\|h-g\|_{H^{N-5}}&\lesssim \epsilon_1^2(1+t)^{p_0-1/2},\\
\|h-g\|_{W^{N_1+5,\infty}}&\lesssim \epsilon_1^2(1+t)^{-1},\\
\|\Gamma(h-g)\|_{H^{N_1-5}}&\lesssim \epsilon_1^2(1+t)^{p_0-1/2}.
\end{split}
\end{align}
The bounds \eqref{scr2} and \eqref{scr3} follow easily from \eqref{a-priori wwbound}, \eqref{ps-a1} and Lemma \ref{wwlemma1}.\newline

 {\noindent\emph{ Proof of Lemma \ref{lemmascr1}.}}
According to the definitions \eqref{ps4} and \eqref{sc-a1}, we see that in order to prove Lemma \ref{lemmascr1}, it suffices to show each term in $\mathcal{N}_R$ satisfies \eqref{scr4}--\eqref{scr7}. In this proof, we mainly concentrate on the term
$$
\mathcal{N}_R^{++++}:=\mathcal{O}[\mathcal{O}[h,q^{++}]h,b^{++}]h-\mathcal{O}[\mathcal{O}[g,q^{++}]g,b^{++}]g,
$$
and the treatments for the other terms are similar. Decompose this term as
$$
\mathcal{N}_R^{++++}=\mathcal{N}_{R1}^{++++}+\mathcal{N}_{R2}^{++++}+\mathcal{N}_{R3}^{++++}
$$
with
\begin{align*}
\mathcal{N}_{R1}^{++++}&:=\mathcal{O}[\mathcal{O}[h,q^{++}]h,b^{++}](h-g),\\
\mathcal{N}_{R2}^{++++}&:=\mathcal{O}[\mathcal{O}[h-g,q^{++}]h,b^{++}]g,\\
\mathcal{N}_{R3}^{++++}&:=\mathcal{O}[\mathcal{O}[g,q^{++}](h-g),b^{++}]g.
\end{align*}

We first show \eqref{scr4}. Using \eqref{ps-a2}, \eqref{ps1}, Lemma \ref{multiplierlemma2d} and the bounds \eqref{scr2}--\eqref{scr3}, we see
\begin{align*}
\|\langle\xi\rangle^{N-20}\mathscr{F}\mathcal{N}_{R1}^{++++}\|_{L^\infty}&\lesssim \|\mathcal{N}_{R1}^{++++}\|_{W^{N-20,1}}
\lesssim \|\mathcal{O}[h,q^{++}]h\|_{H^{N-15}}\|h-g\|_{H^{N-15}}\\
&\lesssim \|h\|_{H^N}\|h\|_{L^{\infty}}\|h-g\|_{H^{N-5}}\lesssim \epsilon_1^4 (1+t)^{2p_0-1},
\end{align*}
and
\begin{align*}
\|\langle\xi\rangle^{N-20}\mathscr{F}\mathcal{N}_{R2}^{++++}\|_{L^\infty}&
\lesssim  \|\mathcal{O}[h-g,q^{++}]h\|_{H^{N-15}}\|g\|_{H^{N-15}}\\
&\lesssim (\|h-g\|_{H^{N-5}}\|h\|_{L^{\infty}}+\|h-g\|_{L^{\infty}}\|h\|_{H^N})\|g\|_{H^{N-5}}\\
&\lesssim \epsilon_1^4 (1+t)^{2p_0-1}.
\end{align*}
The argument for the term $\mathcal{N}_{R3}^{++++}$ is similar as above. Hence, the bound \eqref{scr4} follows.

Similarly, we have
\begin{align*}
\|\mathcal{N}_{R1}^{++++}\|_{H^{N_1+10}}&\lesssim \|\mathcal{O}[h,q^{++}]h\|_{H^{N_1+15}}\|h-g\|_{L^{\infty}}+\|\mathcal{O}[h,q^{++}]h\|_{L^{\infty}}\|h-g\|_{H^{N_1+15}}\\
&\lesssim \|h\|_{H^N}\|h\|_{L^{\infty}}\|h-g\|_{{L^{\infty}}}+\|h\|_{W^{4,\infty}}^2\|h-g\|_{H^{N_1+15}}\\
&\lesssim \epsilon_1^4 (1+t)^{p_0-3/2},
\end{align*}
and
\begin{align*}
\|\mathcal{N}_{R2}^{++++}\|_{H^{N_1+10}}&\lesssim \|\mathcal{O}[h-g,q^{++}]h\|_{H^{N_1+15}}\|g\|_{L^{\infty}}+\|\mathcal{O}[h-g,q^{++}]h\|_{L^{\infty}}\|g\|_{H^{N_1+15}}\\
&\lesssim (\|h-g\|_{H^{N-5}}\|h\|_{L^{\infty}}+\|h-g\|_{L^{\infty}}\|h\|_{H^{N}})\|g\|_{{L^{\infty}}}\\
&\quad+\|h-g\|_{W^{4,\infty}}\|h\|_{W^{4,\infty}}\|g\|_{H^{N_1+15}}\\
&\lesssim \epsilon_1^4 (1+t)^{p_0-3/2}.
\end{align*}
Also, we can deal with the term $\|\mathcal{N}_{R3}^{++++}\|_{H^{N_1+10}}$ in a similar way. Combining these estimates yields \eqref{scr5} as desired.

Now we  prove the weighted estimate \eqref{scr6}. As \eqref{ps-a3}, we have
\begin{align*}
\Gamma\mathcal{N}_{R1}^{++++}=W_1+W_2+W_3
\end{align*}
with
\begin{align*}
\widehat{W_1}(\xi)&:=\mathscr{F}(\mathcal{O}[\Gamma\mathcal{O}[h,q^{++}]h,b^{++}](h-g))(\xi),\\
\widehat{W_2}(\xi)&:=\mathscr{F}(\mathcal{O}[\mathcal{O}[h,q^{++}]h,b^{++}]\Gamma(h-g))(\xi),\\
\widehat{W_3}(\xi)&:=\frac{i}{2\pi}\int_\mathbb{R}\partial_\xi b^{++}(\xi-\eta,\eta)\mathscr{F}(\partial_t\mathcal{O}[h,q^{++}]h)(\xi-\eta)\widehat{(h-g)}(\eta)d\eta\\
&\quad+\frac{i}{2\pi}\int_\mathbb{R}\partial_\xi b^{++} (\xi-\eta,\eta)\mathscr{F}(\mathcal{O}[h,q^{++}]h)(\xi-\eta)\widehat{(h-g)_t}(\eta)d\eta\\
&\quad+\frac{i}{2\pi}\int_\mathbb{R}\partial_\eta b^{++} (\xi-\eta,\eta)\mathscr{F}(\mathcal{O}[h,q^{++}]h)(\xi-\eta)\widehat{(h-g)_t}(\eta)d\eta.
\end{align*}
By expanding $\Gamma\mathcal{O}[h,q^{++}]h$ as  \eqref{ps-a3}, we can obtain
\begin{align*}
\|\Gamma\mathcal{O}[h,q^{++}]h\|_{H^{N_1-5}}\lesssim \epsilon_1^2(1+t)^{p_0-1/2}.
\end{align*}
Hence, using also \eqref{scr3}, we have
\begin{align*}
\|W_1\|_{H^{N_1-10}}&\lesssim \|\Gamma\mathcal{O}[h,q^{++}]h\|_{H^{N_1-5}}\|h-g\|_{L^\infty}+\|\Gamma\mathcal{O}[h,q^{++}]h\|_{L^{2}}\|h-g\|_{W^{N_1-5,\infty}}\\
&\lesssim \epsilon_1^4(1+t)^{p_0-3/2}.
\end{align*}
Similarly, from \eqref{scr2}--\eqref{scr3}, there holds
\begin{align*}
\|W_2\|_{H^{N_1-10}}&\lesssim \|\mathcal{O}[h,q^{++}]h\|_{W^{N_1-5,\infty}}\|\Gamma(h-g)\|_{L^2}+\|\mathcal{O}[h,q^{++}]h\|_{L^{\infty}}\|\Gamma(h-g)\|_{H^{N_1-5}}\\
&\lesssim \epsilon_1^4(1+t)^{p_0-3/2}.
\end{align*}
To estimate $W_3$, note that
\begin{align*}
&\|h_t\|_{H^{N-5}}+\|g_t\|_{H^{N-10}}\lesssim \epsilon_1(1+t)^{p_0},\\
&\|h_t\|_{W^{N_1,\infty}}+\|g_t\|_{W^{N_1,\infty}}\lesssim \epsilon_1(1+t)^{-1/2},\\
&\|(h-g)_t\|_{L^\infty}\lesssim \|h_t\|_{L^\infty}\|h\|_{W^{5,\infty}}+\|h_t\|_{W^{5,\infty}}\|h\|_{L^\infty}\lesssim \epsilon_1^2(1+t)^{p_0-1},\\
&\|(h-g)_t\|_{H^{N_1-5}}\lesssim \|h_t\|_{L^\infty}\|h\|_{H^{N_1}}+\|h_t\|_{H^{N_1}}\|h\|_{L^\infty}\lesssim \epsilon_1^2(1+t)^{p_0-1/2},
\end{align*}
which can be verified by  the equations \eqref{EP36}, \eqref{ps3} and the identity \eqref{ps-a1}, then
\begin{align*}
\|W_3\|_{H^{N_1-10}}&\lesssim\|\partial_t\mathcal{O}[h,q^{++}]h\|_{H^{N_1-5}}\|h-g\|_{L^\infty}
+\|\partial_t\mathcal{O}[h,q^{++}]h\|_{L^{\infty}}\|h-g\|_{H^{N_1-5}}\\
&\quad+\|\mathcal{O}[h,q^{++}]h\|_{H^{N_1-5}}\|(h-g)_t\|_{L^\infty}
+\|\mathcal{O}[h,q^{++}]h\|_{L^{\infty}}\|(h-g)_t\|_{H^{N_1-5}}\\
&\lesssim \epsilon_1^4(1+t)^{p_0-3/2}.
\end{align*}
Therefore, we conclude that
$$
\|\Gamma\mathcal{N}_{R1}^{++++}\|_{H^{N_1-10}}\lesssim \|W_1\|_{H^{N_1-10}}+\|W_2\|_{H^{N_1-10}}+\|W_3\|_{H^{N_1-10}}\lesssim \epsilon_1^4(1+t)^{p_0-3/2}.
$$
Moreover, we can estimate the $H^{N_1-10}$ norm of $\Gamma\mathcal{N}_{R2}^{++++}$ and $\Gamma\mathcal{N}_{R3}^{++++}$ in a similar way as above, and we omit further details for simplicity. Thus, the bound \eqref{scr6} is valid.

Finally, by  similar argument as the proof  of \eqref{ps13}, it is straightforward  to obtain the desired bound for $\|x\mathcal{N}_R\|_{H^{N_1-10}}$. This ends the proof of the lemma.
\hfill $\Box$

 Now, we end this subsection by presenting the proof of Proposition \ref{scpropostion2}.
\newline

{\noindent\emph{Proof of Proposition \ref{scpropostion2}.}} Denote
$$
\widehat{\mathcal{A}}(\xi):=\int_{t_1}^{t_2}e^{i\vartheta(s,\xi)}e^{is\langle\xi\rangle}\widehat{\mathcal{N}_R}(s,\xi)ds
=\widehat{\mathcal{A}}_1(\xi)+\widehat{\mathcal{A}}_2(\xi),
$$
where
\begin{align*}
\widehat{\mathcal{A}}_1(\xi)&
:=\int_{t_1}^{t_2}e^{i\vartheta(s,\xi)}e^{is\langle\xi\rangle}
\big(1-\varphi(\frac{\xi}{(1+s)^{p_0}})\big)\widehat{\mathcal{N}_R}(s,\xi)ds,\\
\widehat{\mathcal{A}}_2(\xi)&
:=\int_{t_1}^{t_2}e^{i\vartheta(s,\xi)}e^{is\langle\xi\rangle}\varphi(\frac{\xi}{(1+s)^{p_0}})\widehat{\mathcal{N}_R}(s,\xi)ds.
\end{align*}
By \eqref{scr4} and the fact $|\xi|\gtrsim (1+s)^{p_0}$, it is easy to see
\begin{align}
\|\langle\xi\rangle^{N_1+10}\widehat{\mathcal{A}}_1(\xi)\|_{L^\infty}&\lesssim \int_{t_1}^{t_2}\|\langle\xi\rangle^{N-20}\widehat{\mathcal{N}_R}(s,\xi)\|_{L^\infty} (1+s)^{-(N-N_1-30)p_0}ds\nonumber\\
&\lesssim \int_{t_1}^{t_2}\epsilon_1^4 (1+s)^{2p_0-1}
(1+s)^{-(N-N_1-30)p_0}ds\nonumber \\
&\lesssim \epsilon_1^4 (1+t_1)^{-(N-N_1-32)p_0}.\label{srrcf0}
\end{align}
For the term $\mathcal{A}_2$, we use
\begin{align}\label{srrcf1}
\|\langle\xi\rangle^{N_1+10}\widehat{\mathcal{A}_2}(\xi)\|_{L^\infty}\lesssim \|\mathcal{A}_2\|_{H^{N_1+10}}+ \|x\mathcal{A}_2\|_{H^{N_1+10}}.
\end{align}
The first term in the RHS of \eqref{srrcf1} can be estimated directly by   \eqref{scr5}
\begin{align}\label{srrcf2}
\|\mathcal{A}_2\|_{H^{N_1+10}}&\lesssim \int_{t_1}^{t_2}\|\mathcal{N}_R(s)\|_{H^{N_1+10}}ds\lesssim \int_{t_1}^{t_2}\epsilon_1^4 (1+s)^{p_0-3/2}ds\lesssim \epsilon_1^4 (1+t_1)^{p_0-1/2}.
\end{align}
To estimate the term $\|x\mathcal{A}_2\|_{H^{N_1+10}}$, we apply $\partial_\xi$ to $\widehat{\mathcal{A}}_2$. In view of  \eqref{ps9},
\begin{align*}
\langle\xi\rangle\partial_\xi(e^{it\langle\xi\rangle}\widehat{\mathcal{N}_R})&=e^{it\langle\xi\rangle}\widehat{\widetilde{\Gamma} \mathcal{N}_R}\\
&=e^{it\langle\xi\rangle}\widehat{\Gamma \mathcal{N}_R}-e^{it\langle\xi\rangle}[(\partial_t+i\langle\xi\rangle)\widehat{ x\mathcal{N}_R}]\\
&=e^{it\langle\xi\rangle}\widehat{\Gamma \mathcal{N}_R}-\partial_t[e^{it\langle\xi\rangle}\widehat{ x\mathcal{N}_R}],
\end{align*}
then $\partial_\xi\widehat{\mathcal{A}_2}$ can be   decomposed into
\begin{align*}
\partial_\xi\widehat{\mathcal{A}_2}(\xi)&=\int_{t_1}^{t_2}\partial_\xi[e^{i\vartheta(s,\xi)}\varphi(\frac{\xi}{(1+s)^{p_0}})] e^{is\langle\xi\rangle}\widehat{\mathcal{N}_R}(s,\xi)ds\\
&\quad+\int_{t_1}^{t_2}e^{i\vartheta(s,\xi)}\varphi(\frac{\xi}{(1+s)^{p_0}})\partial_\xi [e^{is\langle\xi\rangle}\widehat{\mathcal{N}_R}(s,\xi)]ds\nonumber\\
&=\widehat{\mathcal{A}}_{21}(\xi)+\widehat{\mathcal{A}}_{22}(\xi)
+\widehat{\mathcal{A}}_{23}(\xi).
\end{align*}
where
\begin{align*}
\widehat{\mathcal{A}}_{21}(\xi)&:=\int_{t_1}^{t_2}\partial_\xi[e^{i\vartheta(s,\xi)}\varphi(\frac{\xi}{(1+s)^{p_0}})] e^{is\langle\xi\rangle}\widehat{\mathcal{N}_R}(s,\xi)ds,\\
\widehat{\mathcal{A}}_{22}(\xi)&:=\int_{t_1}^{t_2}e^{i\vartheta(s,\xi)}\langle\xi\rangle^{-1} \varphi(\frac{\xi}{(1+s)^{p_0}})e^{is\langle\xi\rangle}\widehat{\Gamma\mathcal{N}_R}(s,\xi)ds,\\
\widehat{\mathcal{A}}_{23}(\xi)&:=-\int_{t_1}^{t_2}e^{i\vartheta(s,\xi)}\langle\xi\rangle^{-1}
\varphi(\frac{\xi}{(1+s)^{p_0}})\partial_s[ e^{is\langle\xi\rangle}\widehat{x\mathcal{N}_R}(s,\xi)]ds.
\end{align*}
Using \eqref{definition of theta}, \eqref{resonant2}, \eqref{sc-a2} and \eqref{a-priori wwbound}, we have
$$
|\partial_\xi\vartheta(t,\xi)|\lesssim \langle\xi\rangle^8\ln(1+t)\sup _{s\in[0,t]}(|\widehat{w}(s,\xi)|^2+|\partial_\xi\widehat{w}(s,\xi)|^2)\lesssim \epsilon_1^2(1+t)^{3p_0},
$$
hence, by \eqref{scr5},
\begin{align*}
\|\mathcal{A}_{21}\|_{H^{N_1+10}}\lesssim \epsilon_1^4\int_{t_1}^{t_2}(1+s)^{4p_0-3/2}ds\lesssim \epsilon_1^4(1+t_1)^{4p_0-1/2}.
\end{align*}
For the term $\mathcal{A}_{22}$, we obtain from  \eqref{scr6} that
\begin{align*}
\|\mathcal{A}_{22}\|_{H^{N_1+10}}&\lesssim \int_{t_1}^{t_2}(1+s)^{19p_0}\|\Gamma\mathcal{N}_R\|_{H^{N_1-10}}ds\\
&\lesssim \epsilon_1^4\int_{t_1}^{t_2}(1+s)^{20p_0-3/2}ds\lesssim \epsilon_1^4(1+t_1)^{20p_0-1/2}.
\end{align*}
To estimate $\|\mathcal{A}_{23}\|_{H^{N_1+10}}$, using integration by parts in time, the bound
$$
|\partial_s\vartheta(s,\xi)|\lesssim (1+s)^{-1} \langle\xi\rangle^8|\widehat{w}(s,\xi)|^2 \lesssim \epsilon_1^2(1+s)^{-1},
$$
and  \eqref{scr7}, we  obtain
$$
\|\mathcal{A}_{23}\|_{H^{N_1+10}}\lesssim \epsilon_1^4 (1+t_1)^{20p_0-1/2}.
$$
We finally  conclude that
\begin{align}\label{srrcf3}
\|x\mathcal{A}_2\|_{H^{N_1+10}}\lesssim \|\mathcal{A}_{21}\|_{H^{N_1+10}}+\|\mathcal{A}_{22}\|_{H^{N_1+10}}+\|\mathcal{A}_{23}\|_{H^{N_1+10}}\lesssim \epsilon_1^4 (1+t_1)^{20p_0-1/2}.
\end{align}
Therefore,  the desired bound \eqref{sc7} follows from \eqref{srrcf0}--\eqref{srrcf3}.
\hfill $\Box$

\subsection{Proof of \eqref{desired wwbound1}}

{\noindent \emph{Proof of \eqref{desired wwbound1}}.}
Using Bernstein's inequality, \eqref{xwbound} and \eqref{wn-5}, we have
\begin{align*}
\|xP_{\leq (1+t)^{1/240}}w\|_{H^{N_1+11}}&\lesssim(1+t)^{-1/240}\|P_{\leq (1+t)^{1/240}}w\|_{H^{N_1+11}}+ \|P_{\leq (1+t)^{1/240}}(xw)\|_{H^{N_1+11}}\\
&\lesssim \|w\|_{H^{N-5}}+(1+t)^{1/16}\|P_{\leq (1+t)^{1/240}}(xw)\|_{H^{N_1-4}}\\
&\lesssim (\epsilon_0+\epsilon_1^2)(1+t)^{1/16+p_0}.
\end{align*}
Then we deduce from  the linear estimate \eqref{linear estimate}, the bounds \eqref{desired wwbound3} and \eqref{wn-5} that
\begin{align}
\|P_{\leq (1+t)^{1/240}}g\|_{W^{N_1+10,\infty}}&\lesssim (1+t)^{-1/2}\|\langle\xi\rangle^{N_1+10}\widehat{w}(\xi)\|_{L^\infty}\nonumber\\
&\quad +(1+t)^{-5/8}(\|g\|_{H^{N_1+12}}+\|xP_{\leq (1+t)^{1/240}}w\|_{H^{N_1+11}})\nonumber\\
& \lesssim (1+t)^{-1/2}(\epsilon_0+\epsilon_1^2),\ \ \forall\ t\in[0,T].\label{decay-1v1}
\end{align}
On the other hand, by Bernstein's inequality and \eqref{wn-5}, there also holds
\begin{align}
\|P_{\geq (1+t)^{1/240}} g\|_{W^{N_1+10,\infty}}&\lesssim \|P_{\geq (1+t)^{1/240}} g\|_{H^{N_1 +11}} \nonumber\\
&\lesssim (1+t)^{-(N-N_1-11)/240}\|P_{\geq (1+t)^{1/240}} g\|_{H^{N}}\nonumber\\
& \lesssim (1+t)^{-1/2}(\epsilon_0+\epsilon_1^2).\label{decay-2v1}
\end{align}

Now, we conclude from  \eqref{decay-1v1} and \eqref{decay-2v1} that
$$
\|g(t)\|_{W^{N_1+10,\infty}}\lesssim\|P_{\leq (1+t)^{1/240}}g\|_{W^{N_1+10,\infty}}+\|P_{\geq (1+t)^{1/240}} g\|_{W^{N_1+10,\infty}}\lesssim (1+t)^{-1/2}(\epsilon_0+\epsilon_1^2).
$$
Moreover, if $\epsilon_1$ is small enough, \eqref{ps-a1} and  \eqref{ps12} lead to
$$\|h(t)\|_{W^{N_1+10,\infty}}\sim \|g(t)\|_{W^{N_1+10,\infty}}.$$
Therefore, \eqref{desired wwbound1} follows, and  this also completes the proof of Proposition \ref{decaye}.
\hfill $\Box$

Finally, combing Proposition \ref{eprop}, Propositions \ref{weprop}--\ref{energyxuprop}, Proposition \ref{decaye} and Lemma  \ref{linear dis est}, Theorem \ref{mainthm} follows by standard continuation argument.

\section*{Appendix A}

In this part, we prove the  linear dispersive estimate for Klein-Gordon operator.

\noindent \textbf{Lemma A.1}  \emph{There holds that}
\begin{align}\label{linear estimatev1}\tag{A.1}
\|e^{\pm it\langle\partial_x\rangle}f\|_{L^\infty}\lesssim (1+t)^{-1/2}\|\widehat{f} \|_{L^\infty}
+(1+t)^{-5/8}(\|f\|_{H^2}+\|xf\|_{H^1}),\ \ \ \forall\ t\geq 0.
\end{align}

{\noindent \emph{Proof.}} In this  proof, we only show \eqref{linear estimatev1} in the ``+'' case, since the discussion for the minus case is similar. Note that the estimate \eqref{linear estimatev1} is trivial if $0\leq t\leq 100$.  Denote
 $$J_k:=\frac{1}{2\pi}\left|\int_\mathbb{R}e^{ix\xi+it\langle\xi\rangle}\widehat{P_kf}(\xi)d\xi\right|,\ \ k\in \mathbb{Z},$$
then in order to prove this lemma, it suffices to show
$\sum_{k\in \mathbb{Z}}J_k \lesssim 1$
with $f$ satisfying
\begin{align}\label{LE2}\tag{A.2}
t^{-1/2}\|\widehat{f} \|_{L^\infty}
+t^{-5/8}(\|f\|_{H^2}+\|xf\|_{H^1}) \lesssim 1,\ \ t>100.
\end{align}
Now we divide this proof into four cases.

Case 1: $2^k\leq t^{-1/2}$. In this case, we use $\|\widehat{P_kf}\|_{L^\infty}\lesssim t^{1/2}$ to obtain
\begin{align*}
\sum_{2^k\leq  t^{-1/2}}J_k \lesssim \sum_{2^k\leq  t^{-1/2}} 2^k\|\widehat{P_kf}\|_{L^\infty}
\lesssim t^{1/2}\sum_{2^k\leq  t^{-1/2}} 2^k\lesssim 1.
\end{align*}
\indent Case 2: $2^k\geq t^{5/12}$. By using the estimate $2^{2k}\|\widehat{P_kf}\|_{L^2}\lesssim t^{5/8}$, we have
\begin{align*}
\sum_{2^k\geq  t^{5/12}}J_k \lesssim \sum_{2^k\geq  t^{5/12}} 2^{k/2}\|\widehat{P_kf}\|_{L^2}
\lesssim t^{5/8}\sum_{2^k\geq  t^{5/12}} 2^{k/2}2^{-2k}\lesssim 1.
\end{align*}
\indent Case 3: $t^{-1/2}\leq 2^k\leq t^{5/12}$ and $|x/t|\geq 1/2$.  Note that \eqref{LE2} implies
\begin{align}\label{LE4}\tag{A.3}
\|\widehat{P_kf}\|_{L^2}+2^{2k}\|\widehat{P_kf}\|_{L^2}+\|\partial_\xi\widehat{P_kf}\|_{L^2}
+2^k\|\partial_\xi\widehat{P_kf}\|_{L^2}\lesssim t^{5/8}.
\end{align}
Moreover, in this case we observe that
\begin{align}\label{LE5}\tag{A.4}
\|tf_x\|_{L^2}\leq 2\|xf_x\|_{L^2}\lesssim t^{5/8}\Rightarrow \|f_x\|_{L^2}\lesssim t^{-3/8}\Rightarrow 2^k\|\widehat{P_kf}\|_{L^2}\lesssim t^{-3/8}.
\end{align}
\indent Subcase 3--1: $|\xi|\geq 1/4$. In this subcase, thanks to \eqref{LE5}, it follows from the definition of $J_k$ that
\begin{align*}
\sum_{4^{-1}\leq 2^k\leq t^{5/12}}J_k\lesssim \sum_{4^{-1}\leq 2^k}\|\widehat{P_kf}\|_{L^2}2^{k/2}
\lesssim\sum_{{4^{-1}\leq 2^k }}t^{-3/8}2^{-k}2^{k/2}\lesssim 1.
\end{align*}

\indent Subcase 3--2: $|\xi|\leq 1/4$. Let $\Phi:=x\xi+t\langle\xi\rangle$, using integration by parts, we see that
\begin{align}
J_k&=\left|\int_\mathbb{R}\partial_{\xi}(e^{i\Phi})(i\partial_\xi\Phi)^{-1}\widehat{P_kf}(\xi)d\xi\right|\nonumber\\
&\leq  \left|\int_\mathbb{R}e^{i\Phi}(\partial_\xi\Phi)^{-2}\partial^2_\xi\Phi\widehat{P_kf}(\xi)d\xi\right|+
\left|\int_\mathbb{R}e^{i\Phi}(\partial_\xi\Phi)^{-1}\partial_\xi\widehat{P_kf}(\xi)d\xi\right|,\label{LE3}\tag{A.5}
\end{align}
where $\partial_\xi\Phi=t(xt^{-1}+\xi\langle\xi\rangle^{-1})$ and $\partial^2_\xi\Phi=t\langle\xi\rangle^{-3}$.
If $|\xi|\leq 1/4$, then $|\partial^2_\xi\Phi|\sim t$ and $|\partial_\xi\Phi|\geq t(|x/t|-|\xi|\langle\xi\rangle^{-1}) \geq t/4$ since $|\xi|\langle\xi\rangle^{-1}\leq 1/4$.  With the help of \eqref{LE4}, it follows from \eqref{LE3} that
\begin{align*}
\sum_{t^{-1/2}\leq 2^k\leq 4^{-1}}J_k\lesssim \sum_{ 2^k\leq 4^{-1} }(t^{-1}\|\widehat{P_kf}\|_{L^2}2^{k/2}+t^{-1}\|\partial_\xi\widehat{P_kf}\|_{L^2}2^{k/2})
\lesssim t^{-1}t^{5/8}\lesssim1.
\end{align*}

\indent Case 4: $t^{-1/2}\leq 2^k\leq t^{5/12}$ and $|x/t|\leq 1/2$.

 \indent Subcase 4--1:  $|\xi|\geq 2$. We see that $|\partial_\xi^2\Phi|\lesssim t|\xi|^{-3}$ and $|\partial_\xi\Phi|\gtrsim t$, so from \eqref{LE3}, there holds
\begin{align*}
\sum_{2\leq 2^k\leq t^{5/12}}J_k\lesssim \sum_{2\leq 2^k}(t^{-1}2^{-3k}\|\widehat{P_kf}\|_{L^2}2^{k/2}+t^{-1}2^{-k}2^k\|\partial_\xi\widehat{P_kf}\|_{L^2}2^{k/2})
\lesssim t^{-1}t^{5/8}\lesssim 1.
\end{align*}
\indent Subcase 4--2: $|\xi|\leq 2$. Let $\xi_0$ be the unique root of the equation $\partial_\xi\Phi=0$, i.e., $\xi_0=-\frac{x}{\sqrt{t^2-x^2}}$ and $|\xi_0|\leq 3^{-1/2}$. Then it is easy to see
\begin{align*}
\sum_{t^{-1/2}\leq 2^k\leq 2}J_k\lesssim \sum_{t^{-1/2}\leq 2^k\leq 2,\ l\geq l_0}\left|\int_\mathbb{R} e^{i\Phi}\widehat{P_kf}(\xi)\varphi_l^{(l_0)}(\xi-\xi_0)d\xi\right|=:\sum_{t^{-1/2}\leq 2^k\leq 2,\ l\geq l_0}J_{k,l},
\end{align*}
where $l_0$ is the smallest integer satisfying $2^{l_0}\geq t^{-1/2}$ and
\begin{equation*}
\varphi_l^{(l_0)}(\xi-\xi_0):=\left\{
\begin{array}{ll}
\varphi(|\xi-\xi_0|/2^l)-\varphi(|\xi-\xi_0|/2^{l-1}),&l\geq l_0+1,\\
\varphi(|\xi-\xi_0|/2^{l_0}),&l=l_0
\end{array}
\right.
\end{equation*}
with $\varphi$ the smooth  function given in Section 1. By this definition, $\partial_\xi\Phi$ vanishes in the integral domain of $J_{k,l_0}$, and  we estimate this term as
\begin{align*}
J_{k,l_0}\leq \|\widehat{P_kf}(\xi)\|_{L^\infty}\|\varphi_{l_0}^{(l_0)}(\xi-\xi_0)\|_{L^1}\lesssim t^{1/2} 2^{l_0}\lesssim1.
\end{align*}
For $l\geq l_0+1$, note that $|\partial^2_\xi\Phi|\sim t$ and $$|\partial_\xi\Phi|=|\partial_\xi\Phi(\xi)-\partial_\xi\Phi(\xi_0)|=|\partial^2_\xi\Phi(\xi^*)||\xi-\xi_0|\sim t2^l,$$
so integrating by parts in $\xi$ as \eqref{LE3}, we can obtain
\begin{align*}
&\sum_{t^{-1/2}\leq 2^k\leq 2,\ l> l_0}J_{k,l}\lesssim\sum_{l> l_0}\Big( t^{-1}2^{-2l}\|\widehat{P_kf}\|_{L^\infty}\|\varphi_{l}^{(l_0)}(\xi-\xi_0)\|_{L^1}
\\
&\qquad\qquad+t^{-1}2^{-l}\|\partial_\xi\widehat{P_kf}\|_{L^2}\|\varphi_{l}^{(l_0)}(\xi-\xi_0)\|_{L^2}
+t^{-1}2^{-l}\|\widehat{P_kf}\|_{L^\infty}\|\partial_\xi\varphi_{l}^{(l_0)}(\xi-\xi_0)\|_{L^1}\Big)\\
&\qquad\lesssim\sum_{l> l_0}( t^{-1}2^{-2l}t^{1/2}2^l+t^{-1}2^{-l}t^{5/8}2^{l/2}+t^{-1}2^{-l}t^{1/2})\lesssim1.
\end{align*}
This ends the proof of the lemma.
\hfill $\Box$

\section*{Appendix B}

In this appendix, we collect some analysis lemmas.

\noindent \textbf{Lemma B.1}
\emph{There holds}
\begin{align}
\|\mathcal{O}[f, M]V\|_{L^2(\mathbb{R})} \lesssim \| M(\xi, \eta- \xi)\|_{L_\eta^\infty H_\xi^1}\|f\|_{L^\infty (\mathbb{R})}\|V\|_{L^2(\mathbb{R})}.\label{5.5}\tag{B.1}
\end{align}

{\noindent \emph{Proof.}} Let $\mathscr{F}_x^\xi$ denote the Fourier transform from $x$ to $\xi$. Using  H\"{o}lder's inequality, we can see
\begin{align}
\left|\langle\mathcal{O}[f, M]V,W\rangle\right|&=  \frac{1}{(2\pi)^2}\Big|\int_{\mathbb{R}^2}  M(\xi, \eta- \xi)\widehat{f}(\xi)\widehat{V}(\eta - \xi) \overline{\widehat{W}}(\eta) d\xi d\eta \Big| \nonumber\\
&= \frac{1}{(2\pi)^2}\Big|\int_{\mathbb{R}}   \Big(\int_{\mathbb{R}}  M(\xi, \eta- \xi)\mathscr{F}_x^\xi \mathscr{F}_y^\eta (f(x+y)V(y)) d\xi \Big) \overline{\widehat{W}}(\eta)d\eta \Big|\nonumber\\
&= \frac{1}{2\pi}\Big| \int_{\mathbb{R}}   \Big(\int_{\mathbb{R}}  \mathscr{F}_\xi^{-1} M(\xi, \eta- \xi)(x) \mathscr{F}_y^\eta (f(-x+y)V(y))dx \Big) \overline{\widehat{W}}(\eta)d\eta \Big|\nonumber\\
&\lesssim \int_{\mathbb{R}}  \| \langle x \rangle \mathscr{F}_\xi^{-1} M(\xi, \eta- \xi)\|_{L_x^2} \cdot\|\langle x \rangle^{-1}\mathscr{F}_y^\eta (f(-x+y)V(y)) \|_{L_x^2}\cdot |\overline{\widehat{W}}(\eta) | d\eta \nonumber\\
&\lesssim   \| \langle x\rangle \mathscr{F}_\xi^{-1}  M(\xi, \eta- \xi)\|_{L_\eta^\infty L_x^2}\cdot \|\langle x \rangle^{-1} \mathscr{F}_y^\eta (f(-x+y)V(y)) \|_{L_\eta^2 L_x^2}\cdot \|\widehat{W}(\eta) \|_{L_\eta^2} \nonumber\\
&\lesssim  \|   M(\xi, \eta- \xi)\|_{L_\eta^\infty H_\xi^1} \cdot\|\langle x \rangle^{-1} f(-x+y)V(y) \|_{L_y^2 L_x^2}\cdot
\|W \|_{L^2}. \nonumber
\end{align}
Note that
\begin{align}
\|\langle x \rangle^{-1} f(-x+y)V(y) \|_{L_y^2 L_x^2} \lesssim \|f\|_{L^\infty} \|V\|_{L^2}, \nonumber
\end{align}
then the desired estimate \eqref{5.5} follows by duality argument.
$\hfill\Box$

\noindent \textbf{Lemma B.2} \emph{Let} $m(\xi,\eta)$ \emph{be a Fourier multiplier satisfying}
\begin{align}\label{B.2}\tag{B.2}
\|m\|_{L^2(\mathbb{R}^2)}+\|\partial_\xi^2 m\|_{L^2(\mathbb{R}^2)}+\|\partial_\eta^2 m\|_{L^2(\mathbb{R}^2)}\lesssim 1,
\end{align}
\emph{then for any} $p_0,\ p_1,\ p_2\in [1,+\infty]$ \emph{with} $p_0^{-1}=p_1^{-1}+p_2^{-1}$, \emph{\emph{we have}}
\begin{align}\label{B.3}\tag{B.3}
\|\mathcal{O}[f_1, m]f_2\|_{L^{p_0}(\mathbb{R})}
\lesssim
\|f_1\|_{L^{p_1}(\mathbb{R})}\|f_2\|_{L^{p_2}(\mathbb{R})}.
\end{align}

{\noindent \emph{Proof.}} Define
$$
K(x,y):=\mathscr{F}^{-1}[m(\xi,\eta)]=\frac{1}{(2\pi)^2}\int_{\mathbb{R}^2}e^{i(x\xi+y\eta)}m(\xi,\eta) d\xi d\eta.
$$
Note that
$$
\mathscr{F}_{\tilde{x}}^{\xi}[f_1(x+\tilde{x})]=e^{ix\xi}\widehat{f}_1(\xi),\ \ \ \ \mathscr{F}_{\tilde{y}}^{\eta}[f_2(x+\tilde{y})]=e^{ix\eta}\widehat{f}_2(\eta),
$$
where $\mathscr{F}_{\tilde{x}}^\xi$ is the Fourier transform from $\tilde{x}$ to $\xi$, then by \eqref{defofO},
\begin{align*}
(\mathcal{O}[f_1, m]f_2)(x)&=\frac{1}{(2\pi)^2}\int_{\mathbb{R}^2} e^{ix(\xi+\eta)}m(\xi,\eta)\widehat{f}_1(\xi)\widehat{f}_2(\eta)d\xi d\eta\\
&=\frac{1}{(2\pi)^4}\int_{\mathbb{R}^2} \mathscr{F}_{\tilde{x}}^{\xi}\mathscr{F}_{\tilde{y}}^{\eta}[K(\tilde{x},\tilde{y})]
\mathscr{F}_{\tilde{x}}^{\xi}\mathscr{F}_{\tilde{y}}^{\eta}[f_1(x+\tilde{x})f_2(x+\tilde{y})]d\xi d\eta\\
&=\frac{1}{(2\pi)^2}\int_{\mathbb{R}^2} K(\tilde{x},\tilde{y})
f_1(x-\tilde{x})f_2(x-\tilde{y})d\tilde{x} d\tilde{y},
\end{align*}
where we have used the identity $\langle\widehat{F},\widehat{G}\rangle=(2\pi)^2\langle F,G \rangle$ ($F,G: \mathbb{R}^2\rightarrow \mathbb{C}$) in the last step. Hence, by H\"{o}lder's inequality,
\begin{align}
\|\mathcal{O}[f_1, m]f_2\|_{L^{p_0}(\mathbb{R})}&\lesssim \int_{\mathbb{R}^2} |K(\tilde{x},\tilde{y})|\cdot
\|f_1(x-\tilde{x})\|_{L^{p_1}_x(\mathbb{R})}\|f_2(x-\tilde{y})\|_{L^{p_2}_x(\mathbb{R})}d\tilde{x} d\tilde{y}\nonumber\\
&\lesssim \|K\|_{L^1(\mathbb{R}^2)}\|f_1\|_{L^{p_1}(\mathbb{R})}\|f_2\|_{L^{p_2}(\mathbb{R})}\label{B.4}\tag{B.4}
\end{align}
with $\frac{1}{p_0}=\frac{1}{p_1}+\frac{1}{p_2}$. Moreover, using \eqref{B.2}, we have
\begin{align}
\|K(x,y)\|_{L^1(\mathbb{R}^2)}&\leq \|(1+x^2+y^2)^{-1}\|_{L^2(\mathbb{R}^2)}\|(1+x^2+y^2)K(x,y)\|_{L^2(\mathbb{R}^2)}
\nonumber\\
&\lesssim \|K(x,y)\|_{L^2(\mathbb{R}^2)}+\|x^2K(x,y)\|_{L^2(\mathbb{R}^2)}+\|y^2K(x,y)\|_{L^2(\mathbb{R}^2)}\nonumber\\
&\sim \|m(\xi,\eta)\|_{L^2(\mathbb{R}^2)}+\|\partial_{\xi}^2m(\xi,\eta)\|_{L^2(\mathbb{R}^2)}+
\|\partial_{\eta}^2m(\xi,\eta)\|_{L^2(\mathbb{R}^2)}\nonumber\\
&\lesssim 1.\label{B.5}\tag{B.5}
\end{align}
Therefore, the desired bound \eqref{B.3} follows from \eqref{B.4} and \eqref{B.5}.
$\hfill\Box$

\noindent \textbf{Lemma B.3} \emph{If} $m(\eta,\sigma)$ \emph{is a Fourier multiplier with} $\eta$ \emph{and} $\sigma$ \emph{localized in the size} $2^k$ \emph{and} $2^l$, \emph{respectively}, \emph{and satisfies}
\begin{align*}
|\partial_{\eta}^a\partial_{\sigma}^b m|\lesssim  A2^{-a k}2^{-b l}\  (\mathrm{resp.}\  A)
\end{align*}
\emph{for any} $a,b=0,1,2$, \emph{then we have}
\begin{align*}
\|\mathscr{F}^{-1}m\|_{L^1(\mathbb{R}^2)}\lesssim A\ (\mathrm{resp.}\  A2^k2^l).
\end{align*}

{\noindent \emph{Proof.}} Let $K(x,y):=\mathscr{F}^{-1}m$, namely,
$$
K(x,y)=(2\pi)^{-2}\int_{\mathbb{R}^2}e^{ix\eta}e^{iy\sigma}m(\eta,\sigma)d\eta d\sigma.
$$
We first assume
\begin{align}\label{B.6}\tag{B.6}
|\partial_{\eta}^a\partial_{\sigma}^b m|\lesssim  A2^{-a k}2^{-b l},\ \ a, b=0,1,2.
\end{align}
Using the localized property of $m$, we see that
\begin{align}\label{B.7}\tag{B.7}
|K(x,y)|\lesssim \|m\|_{L^\infty}2^k2^l\leq A2^k2^l,\ \ \forall\ (x,y)\in \mathbb{R}^2.
\end{align}
On the other hand, with integration by parts, it is easy to see
\begin{align}\label{B.8}\tag{B.8}
|K(x,y)|\lesssim x^{-a}y^{-b}\|\partial_\eta^a\partial^b_\sigma m\|_{L^\infty}2^k2^l,\ \ x\neq0\ \mathrm{and}\ y\neq0.
\end{align}
Let $\mathbb{R}^2=\Omega_1\cup \Omega_2\cup \Omega_3\cup \Omega_4$, where
\begin{align*}
&\Omega_1=\{(x,y);|x|\leq \alpha,\ |y|\leq \beta\},\ \ \Omega_2=\{(x,y);|x|\leq \alpha,\ |y|\geq \beta\},\\
&\Omega_3=\{(x,y);|x|\geq \alpha,\ |y|\leq \beta\},\ \ \Omega_4=\{(x,y);|x|\geq \alpha,\ |y|\geq \beta\}.
\end{align*}
Then using \eqref{B.7}, there holds
$$
\|K(x,y)\|_{L^1(\Omega_1)} \lesssim \alpha\beta A2^k2^l.
$$
Integrating by parts in $\sigma$ only and using \eqref{B.6}, \eqref{B.8} with $(a,b)=(0,2)$, we obtain
$$
\|K(x,y)\|_{L^1(\Omega_2)} \lesssim \alpha\beta^{-1}\|\partial^2_\sigma m\|_{L^\infty}2^k2^l\lesssim \alpha\beta^{-1}A2^k2^{-l}.
$$
Similarly, we can obtain
$$
\|K(x,y)\|_{L^1(\Omega_3)} \lesssim \alpha^{-1}\beta\|\partial^2_\eta m\|_{L^\infty}2^k2^l\lesssim \alpha^{-1}\beta A2^{-k}2^{l}.
$$
Also, with integration by parts in $\eta$ and $\sigma$, we have
$$
\|K(x,y)\|_{L^1(\Omega_4)} \lesssim \alpha^{-1}\beta^{-1}\|\partial^2_\eta\partial^2_\sigma m\|_{L^\infty}2^k2^l\lesssim \alpha^{-1}\beta^{-1}A2^{-k}2^{-l}.
$$
Now, we choose $\alpha,\beta$ satisfying $\alpha2^k=1$ and $\beta2^l=1$, then from the above four estimates, there holds
$$
\|\mathscr{F}^{-1}m\|_{L^1(\mathbb{R}^2)}=\|K\|_{L^1(\mathbb{R}^2)}\lesssim A.
$$

Next, we assume $|\partial_{\eta}^a\partial_{\sigma}^b m|\lesssim  A$ for any $a,b=0,1,2$. In this case, applying the same argument as above with $\alpha=\beta=1$, we can easily see that $\|K\|_{L^1(\mathbb{R}^2)}\lesssim A2^k2^l$. This ends the proof of the lemma.
$\hfill\Box$

\noindent \textbf{Lemma B.4} \emph{For any} $\lambda,\ \mu>0$ and $n\in \mathbb{N}$,  \emph{there holds that}
\begin{align}\label{B.9}\tag{B.9}
\int_{\mathbb{R}^2}e^{ i\lambda xy}\varphi(\mu^{-1}x)\varphi(\mu^{-1}y)dxdy=2\pi\lambda^{-1}+\lambda^{-1-n}\mu^{-2n}O(1),
\end{align}
\emph{where} $\varphi$ \emph{is the smooth radial function used in  the Littlewood-Paley decomposition}.  \emph{The implicit constant coming from the term} $O(1)$ \emph{depends only on } $n$ \emph{and} $\varphi$.

{\noindent \emph{Proof.}} We first set $\lambda=1$. A direct computation gives
\begin{align*}
\mathrm{LHS\ of\ }\eqref{B.9}&=\mu\int_{\mathbb{R}}\varphi(\mu^{-1}x)\widehat{\varphi}(-\mu x)dx
=\int_{\mathbb{R}}\varphi(\mu^{-2}x)\widehat{\varphi}(x)dx\\
&=\int_{\mathbb{R}}\varphi(0)\widehat{\varphi}(x)dx
+\int_{\mathbb{R}}[\varphi(\mu^{-2}x)-\varphi(0)]\widehat{\varphi}(x)dx\\
&=2\pi+\int_{\mathbb{R}}[\varphi(\mu^{-2}x)-\varphi(0)]\widehat{\varphi}(x)dx,
\end{align*}
since $\int_{\mathbb{R}}\widehat{\varphi}(x)dx=2\pi\varphi(0)=2\pi$. Using Taylor's expansion, we have
\begin{align*}
\varphi(\mu^{-2}x)&=\varphi(0)+\sum_{k=1}^{n-1}\frac{\varphi^{(k)}(0)}{k!}\left(\frac{x}{\mu^2}\right)^k
+\frac{\varphi^{(n)}(y)}{n!}\left(\frac{x}{\mu^2}\right)^n
=\varphi(0)+\frac{\varphi^{(n)}(y)}{n!}\left(\frac{x}{\mu^2}\right)^n,
\end{align*}
where $0<|y|<\mu^{-2}|x|$. Hence, there holds
\begin{align*}
\int_{\mathbb{R}}[\varphi(\mu^{-2}x)-\varphi(0)]\widehat{\varphi}(x)dx&
=(\mu^{2n}n!)^{-1}\int_{\mathbb{R}}x^n\widehat{\varphi}(x)\varphi^{(n)}(y)dx.
\end{align*}
Combining the above equalities, we obtain
\begin{align*}
\int_{\mathbb{R}\times \mathbb{R} }e^{-ixy}\varphi(\mu^{-1}x)\varphi(\mu^{-1}y)dxdy=2\pi+\mu^{-2n}O(1),
\end{align*}
and by transformation $\sqrt{\lambda}x\rightarrow x$, $\sqrt{\lambda}y\rightarrow y$, we thus get  \eqref{B.9} as desired. $\hfill\Box$

\section*{Acknowledgments}

L. Han and J. Zhang thank the Division of Applied Mathematics at Brown University for its hospitality, where the work was completed during their visits, supported by the China Scholarship Council. Y. Guo's research  was supported in part by NSFC grant 10828103, NSF grant DMS-0905255 and BICMR. L. Han's research  was supported by the Fundamental Research Funds for the Central Universities. J. Zhang's research was supported  by NSFC grant 11201185, 11471057.

\begin{center}

\end{center}

\noindent \sc{The Division of Applied Mathematics, Brown University}\\
\emph{E-mail address: yan\_guo@brown.edu}
\newline

\noindent {\sc Department of Mathematics and Physics,  North China Electric Power University}\\
\emph{E-mail address: hljmath@ncepu.edu.cn}
\newline

\noindent {\sc Department of Mathematics, Jiaxing University}\\
\emph{E-mail address: zjj@mail.zjxu.edu.cn}


\begin{thebibliography}{99}
\addcontentsline{toc}{section}{References}

  \bibitem{AD}Alazard, T., Delort, J. M.: Sobolev estimates
for two dimensional gravity water waves. arXiv:1307.3836v1 (2013)

\bibitem{Dafermous}  Dafermos, C. M.: \emph{Hyperbolic Conservation Laws in Continuum
Physics}. Grundlehren  der mathematischen
Wissenschaften, Volume 325, 2010.

\bibitem{Deng} Deng, Y., Ionescu, A. and Pausader, B. The Euler-Maxwell system for electrons: global solutions in 2D. Preprint(2014)

 \bibitem{GP} Germain, P.,  Masmoudi, N.:  Global existence for the Euler-Maxwell system. \emph{Ann. Sci. \'{E}c. Norm. Sup\'{e}r.}  \textbf{47},  469--503 (2014)

\bibitem{GMP} Germain, P., Masmoudi, N., Pausader, B.: Nonneutral global solutions for the electron
Euler-Poisson system in three dimensions. \emph{SIAM J. Math. Anal.} \textbf{45} (1), 267--278 (2013)

\bibitem{GMS} Germain, P., Masmoudi, N., Shatah, J.: Global solutions for the gravity water
waves equation in dimension 3. \emph{Ann. Math.} \textbf{175} (2), 691--754 (2012)

\bibitem{Guo} Guo, Y.: Smooth irrotational flows in the large to the Euler-Poisson system in $R^{3+1}$. \emph{Commun. Math. Phys.} \textbf{195}, 249--265 (1998)

\bibitem{GIP} Guo, Y., Ionescu, A., Pausader, B.: Global solutions of the Euler-Maxwell two-fluid system in 3D. arXiv:1303.1060v1 (2013)

\bibitem{GIP2} Guo, Y., Ionescu, A., Pausader, B.:  Global solutions of certain plasma fluid models in three-dimension. \emph{J. Math. Phys.} \textbf{55}, 123102  (2014)

\bibitem{GuoPausader}Guo, Y., Pausader, B.: Global smooth ion dynamics in the Euler-Poisson system.
\emph{Commun. Math. Phys.} \textbf{303}, 89--125 (2011)

\bibitem{GT} Gustafson, S., Nakanishi, K., Tsai, T. P.: Scattering theory for the Gross-Pitaevskii equation in three
dimensions. \emph{Commun. Contemp. Math.} \textbf{11}, 657--707 (2009)

\bibitem{HN1}Hayashi, N., Naumkin, P. I.: The initial value problem for the cubic nonlinear Klein-Gordon equation. \emph{Z. Angew. Math. Phys.} \textbf{59}, 1002--1028 (2008)

\bibitem{HN3}Hayashi, N., Naumkin, P. I.: Quadratic nonlinear Klein-Gordon equation in one dimension. \emph{J. Math. Phys.} \textbf{53}, 103711 (2012)

\bibitem{Hormander} H\"{o}rmander, L.: Lectures on nonlinear hyperbolic differential equations. \emph{Math\'{e}matiques
 Applications} \textbf{26}, Springer-Verlag, Berlin, (1997)

\bibitem{HZG} Han, L., Zhang, J., Guo, B.: Global smooth solution for a kind of two-fluid system
in plasmas. \emph{J. Differential Equations} \textbf{252}, 3453--3481 (2012)

\bibitem{IP}Ionescu, A., Pausader, B.: The Euler-Poisson system in 2D: global stability
of the constant equilibrium solution. \emph{Int. Math. Res. Notices} \textbf{2013} (4), 761--826 (2013).

\bibitem{IP2}Ionescu, A., Pausader, B.: Global solutions of quasilinear systems of Klein-Gordon equations in 3D. \emph{J. Eur. Math. Soc.} \textbf{16}, 2355--2431 (2014)

\bibitem{IPu1}Ionescu, A., Pusateri, F.: Nonlinear fractional Schr\"{o}dinger equations in one dimension. \emph{J. Funct. Anal.} \textbf{266}, 139--176 (2014)

\bibitem{IPu2}Ionescu, A., Pusateri, F.: Global solutions  for the gravity water waves system in 2D. \emph{Invent. Math.}
DOI 10.1007/s00222-014-0521-4

\bibitem{IPu3}Ionescu, A., Pusateri, F.: Global analysis of a model for capillary water waves in 2D. arXiv:1406.6042v1 (2014)

\bibitem{IPu4}Ionescu, A., Pusateri, F.: Global regularity for 2D water waves with surface tension. arXiv:1408.4428v1 (2014)

\bibitem{Jackson} Jackson, J. D.: \emph{Classical Electrodynamics}.  John Wiley \& Sons Inc, 1962.

\bibitem{Jang}Jang, J.: The two-dimensional Euler-Poisson system with
spherical symmetry. \emph{J. Math. Phys.} \textbf{53}, 023701  (2012)

\bibitem{JLZ}Jang, J., Li, D., Zhang, X.: Smooth global solutions for
the two-dimensional Euler-Poisson system. \emph{Forum Math.} \textbf{26}, 645--701 (2014)

\bibitem{LW}Li, D., Wu, Y.: The Cauchy problem for the two dimensional
Euler-Poisson system. \emph{J. Eur. Math. Soc.} \textbf{10}, 2211--2266 (2014)

\bibitem{Shatah} Shatah, J.:  Normal forms and quadratic nonlinear Klein-Gordon equations. \emph{Comm. Pure  Appl. Math.}  \textbf{38} (5), 685--696 (1985)

\bibitem{ST} Sideris, T.: Formation of singularities in three-dimensional compressible fluids. \emph{Commun. Math. Phys.} \textbf{101}, 475--485 (1985)
    


\end{thebibliography}
\end{document}